\documentclass[leqno, 11pt]{amsart}
\usepackage{amsmath, amsthm, amssymb, latexsym}
 \newlength{\auxwidth}
 \newlength{\auxheight}
\newcounter{theorem}
\newcounter{proposition}
\newcounter{lemma}
\newcounter{corollary}
\newcounter{conjecture}
\newtheorem{definition}{Definition}[section]
\newtheorem{theorem}[definition]{Theorem}
\newtheorem{proposition}[definition]{Proposition}
\newtheorem{lemma}[definition]{Lemma}

\theoremstyle{remark}
\newtheorem{example}[definition]{Example}
\newtheorem{remark}[definition]{Remark}

\numberwithin{equation}{section}

\newcommand{\op}[1]{\operatorname{#1}}
\newcommand{\acou}[2]{\ensuremath{\langle #1 , #2 \rangle}}

\newcommand{\brak}[1]{\ensuremath{\langle\! #1\!\rangle}}

\newcommand{\Tr}{\ensuremath{\op{Tr}}}

\newcommand{\Tra}{\ensuremath{\op{Trace}}}

\newcommand{\Hol}{\op{Hol}}

\newcommand{\C}{\ensuremath{\mathbb{C}}} 
\newcommand{\bH}{\ensuremath{\mathbb{H}}} 
\newcommand{\N}{\ensuremath{\mathbb{N}}} 
\newcommand{\R}{\ensuremath{\mathbb{R}}}

\newcommand{\Rd}{\ensuremath{\R^{d+1}}}

\newcommand{\Rdo}{\R^{d+1}\!\setminus\! 0}

\newcommand{\URd}{U\times\R^{d+1}}

\newcommand{\URdo}{U\times(\R^{d+1}\!\setminus\! 0)}

\newcommand{\Ca}[1]{\ensuremath{\mathcal{#1}}}

\newcommand{\cD}{\ensuremath{\mathcal{D}}}
\newcommand{\cE}{\Ca{E}}
\newcommand{\cF}{\Ca{F}}

\newcommand{\cH}{\ensuremath{\mathcal{H}}}

\newcommand{\cK}{\ensuremath{\mathcal{K}}}
\newcommand{\cL}{\ensuremath{\mathcal{L}}}

\newcommand{\cS}{\ensuremath{\mathcal{S}}}

\newcommand{\cU}{\ensuremath{\mathcal{U}}}

\newcommand{\fg}{\ensuremath{\mathfrak{g}}}
\newcommand{\fh}{\ensuremath{\mathfrak{h}}}

\newcommand{\psivdo}{$\Psi_{H}$DO}
\newcommand{\psivdos}{$\Psi_{H}$DO's}

\newcommand{\pvdo}{\ensuremath{\Psi_{H}}}

\newcommand{\psido}{$\Psi$DO} 
\newcommand{\psidos}{$\Psi$DO's} 
\newcommand{\psinf}{\ensuremath{\Psi^{-\infty}}}

\newcommand{\reg}{{\text{reg}}}

\newcommand{\xiy}{{\xi\rightarrow y}}

\newcommand{\yxi}{{y\rightarrow\xi}}

\newcommand{\supp}{\op{supp}}

\newcommand{\rk}{\op{rk}}
\newcommand{\im}{\op{im}}
 
\newcommand{\dom}{\op{dom}}

\newcommand{\End}{\ensuremath{\op{End}}}

\newcommand{\hotimes}{\hat\otimes}

\renewcommand{\Box}{\square}

\newcommand{\Sp}{\op{Sp}}

\newcommand{\dbarb}{\bar\partial_{b}}
 \newtheorem*{acknowledgements}{Acknowledgements}

\begin{document}
\title{Intrinsic notion of principal symbol\\ for the Heisenberg calculus} 

\author{Rapha\"el Ponge}

\address{Department of Mathematics, Ohio State University, Columbus, USA.}
\email{ponge@math.ohio-state.edu}

    \begin{abstract}
        In this paper we define an intrinsic notion of principal for the Hypoelliptic calculus on Heisenberg manifolds. More precisely, 
        the principal symbol of a \psivdo\ appears as a homogeneous section over the linear dual of the 
        tangent Lie algebra bundle of the manifold. This definition is an important step towards using global $K$-theoretic tools in the Heisenberg 
        setting, such as those involved in the elliptic setting for proving the Atiyah-Singer index theorem or the regularity of the eta invariant. 
        On the other hand, the intrinsic definition of the principal symbol enables us to give an intrinsic sense to the model operator of \psivdo\ at 
        point and to give a definitive proof that the Heisenberg calculus is modelled at each point by the calculus of left-invariant \psidos\ on the 
        tangent group at the point. This also allows us to define an intrinsic Rockland condition for \psivdos\ which is shown to be  
        equivalent to the invertibility of the principal symbol, provided that the Levi form has constant rank. Furthermore, we review the 
        main hypoellipticity results on Heisenberg manifolds in terms of the results of the paper. In particular, we complete the treatment of the 
        Kohn Laplacian of~\cite{BG:CHM} and establish that for the horizontal sublaplacian the invertibility of the principal symbol is equivalent to 
        some condition on the Levi form, called condition $X(k)$. Incidentally, this paper provides us with a relatively up-to-date overview of the main 
        facts about the Heisenberg calculus. 
    \end{abstract}

    \maketitle 
    
\section{Introduction} 
The aim of this paper is to define a global and intrinsic notion of the principal symbol for the hypoelliptic calculus on Heisenberg 
manifolds, along with some applications. On the way, we complete the treatment of the Kohn Laplacian of~\cite{BG:CHM} and gives a complete 
treatment for the horizontal sublaplacian on a Heisenberg manifold acting on forms. 

Recall that a Heisenberg manifold $(M,H)$ consists of a manifold $M$ together 
with a distinguished hyperplane bundle $H\subset TM$. This definition 
 includes as main examples the Heisenberg group and its quotients by cocompact lattices, (codimension 1) foliations, CR and contact manifolds, as well as the 
 confolations of~\cite{ET:C}. Moreover, in this context the natural operators such as H{\"o}rmander's sum of squares, the Kohn Laplacian, the horizontal 
 sublaplacian or the contact Laplacian of~\cite{Ru:FDVC}, while they may be hypoelliptic, are definitely not elliptic. 
 Thus, the standard elliptic calculus cannot be used efficiently  to study these operators.  
 
 The relevant substitute to the elliptic calculus is the Heisenberg calculus of Beals-Greiner~\cite{BG:CHM} and Taylor~\cite{Ta:NCMA}, which extends 
 seminal works of Boutet de Monvel~\cite{BdM:HODCRPDO}, Dynin~(\cite{Dy:POHG}, \cite{Dy:APOHSC}) and  Folland-Stein~\cite{FS:EDdbarbCAHG}
 (see also~\cite{BdMGH:POPDCM}, \cite{CGGP:POGD}, \cite{EMM:HAITH}, \cite{Gr:HPOPCD}, \cite{Ho:CHPODC}).  The idea in the 
 Heisenberg calculus, which goes back to Elias Stein, is to build a pseudodifferential calculus on Heisenberg manifolds
  modelled on that of left-invariant pseudodifferential operators on nilpotent groups. 
  This stems from the fact that the relevant tangent structure for a Heisenberg manifold $(M,H)$ is that of a bundle $GM$ of two-step nilpotent Lie 
 groups (see~\cite{BG:CHM}, \cite{Be:TSSRG}, \cite{EMM:HAITH}, \cite{FS:EDdbarbCAHG}, \cite{Gr:CCSSW}, \cite{Po:Pacific1}, \cite{Ro:INA}).  

  In the original monographs~\cite{BG:CHM} and~\cite{Ta:NCMA} the principal symbol is only defined in local coordinates, so the definition 
\emph{a priori} depends on the choice of these coordinates. In the special case of a contact manifold, an intrinsic definition have been given in the 
unpublished book~\cite{EMM:HAITH}, as a section over a bundle of jets of vector fields representing the tangent group bundle of the contact manifold. 
This approach is similar to that of Melrose~\cite{Me:APSIT} in the setting of the $b$-calculus for manifolds with boundary.   
% using an alternative description of the Heisenberg calculus.
% similar that of Melrose's $b$-calculus. 

In this paper we give an intrinsic definition of the principal symbol, valid for an arbitrary Heisenberg manifold,
in terms of the description of the tangent Lie group bundle in~\cite{Po:Pacific1}. 
As a consequence, we can reformulate in a global fashion previously known criterions for existence of parametrices in the Heisenberg calculus. In particular, we 
can define a Rockland condition  for operators in the Heisenberg calculus in a fully intrinsic way and show that this condition allows to determine the existence of 
a parametrix when the Levi form of the Heisenberg manifold has constant rank. 
%  
% 
% This 
% gives us som
% 
% about existence
% 
% 
% 
% As it would be explained in 
% 
% group bundle of
% 
% In this paper that thanks to the results of~\cite{Po:Pacific1}, we can define a global and intrinisic principal symbol for general Heisenberg manifolds 
% by remaining within the original approach 
%  of~\cite{BG:CHM}. 
% There are two main motivations for carrying out this project. First, as shown in this paper, this allows us to 
% have a coordinate free, Rockland-type criterion for the existence of parametrices within the Heisenberg calculus 
% (see Theorem~\ref{thm:Intro.Rockland-Parametrix} below).

More importantly, since our approach of the principal symbol connects nicely with the construction of the tangent groupoid of a Heisenberg manifold 
in~\cite{Po:Pacific1}, this presumably allows us to make use of global $K$-theoretic arguments in the Heisenberg setting, as those involved in the proof of the (full) 
Atiyah-Singer index theorem~(\cite{AS:IEO1}, \cite{AS:IEO3}) and the regularity of the eta invariant for general selfadjoint elliptic \psidos~(\cite{APS:SARG3}, 
\cite{Gi:ITHEASIT}) or, equivalently, the vanishing of the noncommutative residue of a \psido\ projection~(see~\cite{Wo:LISA}, \cite{BL:OEICNLBVP}, 
\cite{Po:JFA1}). More precisely, this paper can also be seen as a step towards extending the aforementioned results to the Heisenberg setting. 

On the other hand, for strictly pseudoconvex CR manifolds this should allows us to recover and to extend in a rather simple way a recent result of 
Hirachi~\cite{Hi:LSSKGISPD} on the invariance of the integral of the logarithmic singularity of the Szeg\"o kernel (see also~\cite{BdM:LTTP}), in connection 
with the program of Fefferman~\cite{Fe:PITCA} for the Szeg\"o and the Bergman kernels on strictly pseudoconvex complex domains. 
% 
% In particular, in case of strictly 
% 
% 
% From this point of view, this paper is a sequel to~\cite{Po:Pacific1} where an analogue for Heisenberg manifolds of Connes' tangent groupoid of a manifold, 
% which is used 
% in~\cite{Co:NCG} to give an alternative proof of the index theorem of Atiyah-Singer. Thus our approach to the principal symbol of a \psivdo\
% should be useful to extend to the Heisenberg setting Connes' proof of the Atiyah-Singer index theorem (see also~\cite{VE:PhD} for 
% related results on this point). 
% 
% On the other hand, we can make use of the results of this paper to give a more general proof of recent results of Hirachi~\cite{Hi:LSSKGISPD} and 
% Boutet de Monvel~\cite{BdM:LTTP} on the contact invariance of the integral of the logarithmic singularity of the Szeg\"o kernel on contact manifolds, in 
% connection with the homotopy invariance of the noncommutative residue of a (bounded) projection in the Heisenberg calculus (see~\cite{Po:MRL1}). 
% 
% It should also be mentioned that the approach to the principal symbol and to the Heisenberg calculus in this paper should allow us to construct a hypoelliptic 
% pseudodifferential calculus on Carnot-Carath\'eodory manifolds which are equiregular in the sense of~\cite{Gr:CCSSW}. 
% Following is a more detailed outline of the results of the paper.

\subsection{Intrinsic notion of principal symbol for the Heisenberg calculus} 
Given a Heisenberg a manifold $(M,H)$ and a vector bundle $\cE$ we let $\fg^{*}M$ denote the linear dual of the Lie algebra bundle associated to the 
tangent Lie group bundle. For $m \in \C$ we let $\pvdo^{m}(M,\cE)$ denote the class of \psivdos\ 
of order $m$ acting on section of $\cE$ and, letting $\pi:\fg^{*}M\rightarrow M$ be the canonical projection, we define
$S_{m}(M,\cE)$ as the space of sections $p_{m}\in C^{\infty}(\fg^{*}M\setminus 0,\End \pi_{*}\cE)$ such that for any $\lambda>0$ we have
% which are homogeneous with 
% respect to the grading of $\fg^{*}M$ induced by that of $GM$, that is, 
\begin{equation}
    p_{m}(\lambda^{2}\xi_{0},\lambda\xi')=\lambda^{m}p_{m}(\xi_{0},\xi'), \qquad \xi_{0}\in  (TM/H)^{*}, \quad \xi'\in H^{*}.
\end{equation}

The principal symbol of an operator $P\in \pvdo^{m}(M,\cE)$  is obtained as an element of $S_{m}(\fg^{*}M,\cE)$ as follows. 
As shown in~\cite{Po:Pacific1} the tangent Lie group bundle $GM$ of a Heisenberg manifold $(M,H)$  can be described as the bundle $(TM/H)\oplus H$ 
endowed with a law group encoded by the Levi form $\cL: H\times H \rightarrow TM/H$ such that, for sections $X$ and $Y$ of $H$, we have
 \begin{equation}
   \cL(X,Y)=[X,Y] \bmod H.
     \label{eq:Intro.Levi-form}
\end{equation}

It is also shown in~\cite{Po:Pacific1} that at a point $x \in M$ this approach to $GM$ is equivalent to the extrinsic one in~\cite{BG:CHM} in terms of the 
Lie group of a nilpotent Lie algebra of jets of vector fields at point of a local chart. This is made by means of a special kind of  
privileged coordinates at $x$, called Heisenberg coordinates. 

As a byproduct of the equivalence between the two approaches, we obtain that in Heisenberg coordinates a Heisenberg diffeomorphism $\phi:(M, H) 
\rightarrow (M',H')$ is well approximated by the induced isomorphism $\phi_{H}':GM\rightarrow GM'$ between the tangent Lie group bundles 
(see~\cite{Po:Pacific1}). This allows us to carry out the proof of the invariance by Heisenberg diffeomorphisms of the Heisenberg calculus in Heisenberg 
coordinates, rather than in privileged coordinates as in~\cite{BG:CHM} (see Section~\ref{sec.Invariance}). As a consequence we can prove: 

\begin{theorem}\label{prop:Intro.principal-symbol}
    Let $P:C^{\infty}_{c}(M,\cE)\rightarrow C^{\infty}(M,\cE)$ be a \psivdo\ of order $m$. Then there exists a unique symbol a unique symbol 
    $\sigma_{m}(P)\in S_{m}(\fg^{*}M, \cE)$ such that, for any $a \in M$, the symbol $\sigma_{m}(P)(a,.)$ agrees in trivializing Heisenberg coordinates centered at 
    $a$ with the principal symbol of $P$ at $x=0$. 
\end{theorem}

We call the symbol $\sigma_{m}(P)(x,\xi)$ the principal symbol of $P$.  In order to distinguish it from the other definition of the principal symbol 
in local coordinates, we will sometimes refer to it as the \emph{global} principal symbol, while the other principal symbol will be called the 
\emph{local} principal symbol. 

Let us now describe the main properties of the global principal symbol. First, we have: 

\begin{proposition}\label{prop:Intro.surjectivity-principal-symbol-map}
    For every $m\in \C$ the principal symbol map $\sigma_{m}: \pvdo^{m}(M,\cE) \rightarrow S_{m}(\fg^{*}M,\cE)$ gives rise to a linear isomorphism
   $\pvdo^{m}(M,\cE)/\pvdo^{m-1}(M,\cE)\stackrel{\sim}{\longrightarrow}  S_{m}(\fg^{*}M,\cE)$. 
\end{proposition}

If $P\in \pvdo^{m}(M,\cE)$ has principal symbol $\sigma_{m}(P)$ then for any $a\in M$ we define the model operator $P^{a}$ of $P$  at $a$ as the operator 
as the left-invariant \psivdo\ on $G_{a}M$ with symbol $\sigma_{m}(a,.)$, that is, the left-convolution operator with 
the inverse Fourier transform of $\sigma_{m}(a,.)$ (see Definition~\ref{def:PsiHDO.model-operator}). 

On the other hand, for any $a \in M$ the product law on $G_{a}M$ defines a product for symbols,
\begin{equation}
    *^{a}:S_{m_{1}}(\fg^{*}_{a}M)\times S_{m_{2}}(\fg^{*}_{a}M) \longrightarrow S_{m_{1}+m_{2}}(\fg^{*}_{a}M).
\end{equation}
This product depends smoothly enough on $a$ to give rise to a product 
\begin{equation}
    *:S_{m_{1}}(\fg^{*}M,\cE)\times S_{m_{2}}(\fg^{*}M,\cE) \longrightarrow S_{m_{1}+m_{2}}(\fg^{*}M,\cE),
\end{equation}
such that for $p_{j}\in S_{m_{j}}(\fg^{*}M,\cE)$, $j=1,2$, we have 
\begin{equation}
      p_{m_{1}}*p_{m_{2}}(a,\xi)=[p_{m_{1}}(a,.)*^{a}p_{m_{2}}(a,.)](\xi) \qquad \forall (a,\xi) \in \fg^{*}M\setminus 0.
\end{equation}
This comes from the fact that in local coordinates the above product is nicely related to the product of local homogeneous symbols of~\cite{BG:CHM}. 
As a consequence we get: 

\begin{proposition}\label{prop:Intro.composition}
 For $j=1,2$ let $P_{j}\in \pvdo^{m_{j}}(M,\cE)$ and suppose that $P_{1}$ or $P_{2}$ 
    is properly supported.\smallskip 
   
    1) We have $\sigma_{m_{1}+m_{2}}(P_{1}P_{2})=\sigma_{m_{1}}(P)*\sigma_{m_{2}}(P)$.\smallskip
    
    2) At every point $a\in M$ the model operator of $P_{1}P_{2}$ is $(P_{1}P_{2})^{a}=P^{a}_{1}P_{2}^{a}$.
\end{proposition}

Next, it is shown in~\cite{BG:CHM} that the transpose of a \psivdo\ is again a \psivdo. Thanks to the results of~\cite{Po:Pacific1}, 
we can prove a version of this result in Heisenberg coordinates (see Section~\ref{sec:transpose}). As a consequence, 
we identify the principal symbol of transpose, for we get:

\begin{proposition}
  Let $P \in \pvdo^{m}(M,\cE)$ have principal symbol $\sigma_{m}(P)$. Then:\smallskip
 
 1) The principal symbol of the transpose $P^{t}$ is $\sigma_{m}(P^{t})(x,\xi)= \sigma_{m}(x,-\xi)^{t}$;\smallskip 
  
  2) If $P^{a}$ is the model operator of $P$ at $a$, then the model operator of $P^{t}$ at $a$ is the transpose operator 
  $(P^{a})^{t}: \cS_{0}(G_{x}M,\cE_{x}^{*})\rightarrow \cS_{0}(G_{x}M,\cE_{x}^{*})$.
\end{proposition}

Assume now that $M$ is endowed with a positive density and $\cE$ with a Hermitian metric, and let $L^{2}(M,\cE)$ be the associated 
$L^{2}$-Hilbert space. Then we have:

\begin{proposition}
    Let $P \in \pvdo^{m}(M,\cE)$ have principal symbol $\sigma_{m}(P)$. Then:\smallskip
   
 1) The principal symbol of the adjoint $P^{*}$ is $\sigma_{\bar{m}}(P^{*})(x,\xi)=\sigma_{m}(P)(x,\xi)^{*}$.\smallskip 
  
  2) If $P^{a}$ denotes the model operator of $P$ at $a\in M$ then the model operator of $P^{*}$ at $a$ is  
  the adjoint  $(P^{a})^{*}$ of $P^{a}$. 
 \end{proposition} 

\subsection{Rockland condition and hypoellipticity}
It is shown in~\cite{BG:CHM} that, in local coordinates, the invertibility of the local principal symbol is equivalent to the existence of a \psivdo\ 
parametrix. Thanks to Proposition~\ref{prop:Intro.composition} and  the relationship in local coordinates between the local and global principal symbols, we can 
reformulate this result in a global fashion. More precisely, we have:  %using the global principal symbol: 

\begin{proposition}\label{thm:Intro.hypoellipticity}
   Let $P:C^{\infty}_{c}(M,\cE)\rightarrow C^{\infty}(M,\cE)$ be a \psivdo\ of order $m$. Then the following are equivalent:\smallskip 

   1) The principal symbol $\sigma_{m}(P)$ of $P$ is invertible with respect to the convolution product for homogeneous 
    symbols;\smallskip 
    
    2) The operator $P$ admits a parametrix $Q$ in $\pvdo^{-m}(M,\cE)$.\smallskip 
 
 \noindent Furthermore, if 1) and 2) hold then $P$ is hypoelliptic with gain of $\frac{k}{2}$-derivatives.
\end{proposition}

In fact, it may be difficult to determine whether the principal symbol of a \psivdo\ is invertible because in general the convolution product for symbols is not 
the pointwise product of symbols. In particular, in the Heisenberg setting, the product of 
principal symbols is neither commutative, nor microlocal.

Nevertheless, to determine the invertibility of the principal symbol we can make use of a general 
 representation theoretic criterion, the Rockland condition: a \psivdo-operator $P$ satisfies the Rockland 
 condition at a point $a\in M$ when its model operator $P^{a}$ satisfies the Rockland condition on $G_{a}M$, i.e., 
 for every nontrivial irreducible unitary representation $\pi$ of $G_{a}M$ some (unbounded) operator $\pi_{P^{a}}$ 
acting on the  representation space of $\pi$ is injective on the space of smooth vectors of $\pi$. 

It is well known that a left-invariant homogeneous \psido\ on a nilpotent group is hypoelliptic if, and only if, it satisfies the Rockland condition, 
and it further admits an inverse if, and only if, together with its transpose it satisfies the Rockland condition (e.g.~\cite{HN:HGNRN3}, \cite{HN:COHHIGGLG},
\cite{CGGP:POGD}). 
Thus if $P$ is a \psivdo\ such that $P$ and $P^{t}$ satisfies 
the Rockland condition at every point then, for any $x \in M$, then the model operator of $P^{x}$ is invertible. However, whether this inverse 
depends smoothly enough with respect to $x$ to yield an inverse for the principal symbol is a more delicate issue. 

We show here that this occurs when the Levi form~(\ref{eq:Intro.Levi-form}) has constant rank, say $2n$. In this case the tangent Lie group is a fiber 
bundle with typical fiber $\bH^{2n+1}\times \R^{d-2n}$. This allows us to make use of results of Christ~\emph{et al.}~\cite{CGGP:POGD} 
about families of \psidos\ on a \emph{fixed} nilpotent graded group to get: 

\begin{proposition}\label{thm:Intro.Rockland-Parametrix}
  Let $P$ be a \psivdo\ of order $m$, $\Re m\geq 0$, and  suppose that the Levi form~(\ref{eq:Intro.Levi-form}) has constant rank. 
  Then the following are equivalent:\smallskip 
    
    (i) $P$ and $P^{t}$ satisfy the Rockland condition at every point of $M$;\smallskip 
    
  (ii) $P$ and $P^{*}$ satisfy the Rockland condition at every point of $M$;\smallskip 
    
    (iii) The principal symbol of $P$ is invertible.\smallskip 
    
   \noindent  In particular,  if (i) and (ii) holds then $P$ admits a \psivdo\ parametrix and is hypoelliptic with loss of $\frac{1}{2}\Re m$ 
   derivatives.
\end{proposition}

In particular, when $P$ is selfadfoint and the Levi form has constant rank the Rockland condition for $P$ is equivalent to the invertibility of the 
principal symbol of $P$. 

Finally, if $2n$ denotes the rank of the Levi form at $a$ then we have $G_{a}M\simeq \bH^{2n+1}\times \R^{d-2n}$. Therefore, the irreducible 
unitary representations of $G_{a}M$ are two kinds, one dimensional representations on $\C$ and infinite dimensional representations on 
$L^{2}(\R^{n})$. In the former case the Rockland condition corresponds to the invertibility of the principal symbol along $H^{*}$, while in the latter 
case it is enough to look at the representations coming from that of $\bH^{2n+1}$ with Planck constants $\pm 1$ (see Section~\ref{sec:hypoellipticity}). 

\subsection{Hypoellipticity criterions for sublaplacians}
It is interesting to look at the previous results in the case of sublaplacians, as this covers several important examples such as H\"ormander's sum of 
squares, the Kohn Laplacian or the horizontal sublaplacian. 

Here by sublaplacian we mean a differential operator 
$\Delta:C^{\infty}(M,\cE)\rightarrow C^{\infty}(M,\cE)$ such that near every point $a \in M$ we can write $\Delta$ in the form,
\begin{equation}
    \Delta=-(X_{1}^{2}+\ldots+X_{d}^{2})-i\mu(x)X_{0}+\op{O}_{H}(1),
     \label{eq:Intro.sublaplacian}
\end{equation}
where  $X_{0},X_{1},\ldots,X_{d}$ is a $H$-frame of $TM$, so that $X_{1},\ldots,X_{d}$ span $H$, the term $\mu(x)$ is a smooth section of $\End \cE$ 
and the notation $\op{O}_{H}(1)$ stands for a differential operator of Heisenberg order~$\leq 1$. 

In this case the Rockland condition and the invertibility of the reduces to the following. Let $L(x)=(L_{jk}(x))$ be the matrix of 
$\cL$ with respect to the $H$-frame $X_{0},X_{1},\ldots,X_{d}$, so that for $j,k=1,\ldots,d$ we have 
\begin{equation}
    \cL(X_{j},X_{k})=[X_{j},X_{k}]=L_{jk}(x)X_{0} \quad \bmod H.
%     \label{eq:¥}
\end{equation}
Let $2n$ be the rank of $\cL_{a}$ and $L(a)$, let  $\lambda_{1},\ldots,\lambda_{d}$ denoting the eigenvalues of $L(a)$ and 
consider the condition, 
 \begin{equation}
    \Sp \mu(a) \cap \Lambda_{a}=\emptyset,\\
     \label{eq:Intro.sublaplacian.condition}
 \end{equation}
 where the singular set $\Lambda_{a}$ is defined as follows, 
 \begin{gather} 
 \Lambda_{a}=     (-\infty, -\frac12 \Tra |L(a)|]\cup [\frac12 \Tra 
   |L(a)|,\infty) \qquad \text{if $2n<d$},\\
    \Lambda_{a}=\{\pm(\frac12 \Tra |L(a)|+2\alpha_{j}|\lambda_{j}|); \alpha_{j}\in \N^{d}\}\qquad \text{if $2n=d$}.
%    \frac{1}{2}\sum_{j=1}^{2n}|\lambda_{j}|(1+2\N) 
 \end{gather}
As is turns out this condition makes sense independently of the choice of the $H$-frame and is the relevant condition to look at in the case of 
 a sublaplacian. More precisely, we have: 

\begin{proposition}\label{prop:Intro.sublaplacian.Rockland-bundle}
    1) The condition~(\ref{eq:Sublaplacian.condition}) makes sense intrinsically for any $a \in M$.\smallskip
    
    2) At every point $a\in M$ the Rockland conditions for $\Delta$ and $\Delta^{t}$ are equivalent to~(\ref{eq:Intro.sublaplacian.condition}).\smallskip
    
    3) The principal symbol of $\Delta$ is invertible if, and only if, the condition~(\ref{eq:Intro.sublaplacian.condition}) holds at every point of $M$. 
    Moreover, when the latter occurs $\Delta$ admits a parametrix in $\pvdo^{-2}(M,\cE)$ and is hypoelliptic with loss of one derivative.
\end{proposition}

The above result is proved in~\cite{BG:CHM} in the case of scalar sublaplacians, for which the coefficient $\mu(x)$ in~(\ref{eq:Intro.sublaplacian}) 
is a complex-valued function, but the general case is not dealt with in~\cite{BG:CHM}.  However, it is necessary to extend the result to sublaplacians acting on sections 
of vector bundles in order to deal with the Kohn Laplacian and the horizontal sublaplacian acting on forms. In particular, 
Proposition~\ref{prop:Intro.sublaplacian.Rockland-bundle} allows us to complete the treatment of the Kohn Laplacian in~\cite{BG:CHM} (see below).
 
 \subsection{Main examples of hypoelliptic operators on Heisenberg manifolds}
We devote a section where we explain how the Heisenberg calculus, including the results of this paper, provides us with a unified framework to deal with the 
hypoellipticity of the main geometric operators on a Heisenberg manifold and to recover well-known results proved using various different approaches. 

\subsubsection*{(a) H\"ormander's sum of squares} A H\"ormander's sum of squares is an operator of the form $\Delta=-(X_{1}^{2}+\ldots+X_{m}^{2})$ 
where $X_{1},\ldots,X_{m}$ are vector fields on $M$. When $X_{1},\ldots,X_{m}$ span $H$ we get a sublaplacian and 
Proposition~\ref{prop:Intro.sublaplacian.Rockland-bundle}  allows us to 
recover, in this special case, the celebrated result of H\"ormander~\cite{Ho:HSODE} about the hypoellipticity of sum of squares under the bracket 
condition. 

\subsubsection*{(b) Kohn Laplacian} The Kohn Laplacian is the Laplacian associated to the tangential Cauchy-Riemann complex, or $\dbarb$-complex, on 
a CR manifold~(\cite{KR:EHFBCM}, \cite{Ko:BCM}). It was shown by Kohn~\cite{Ko:BCM} that under the condition $Y(q)$ the Kohn Laplacian acting on 
$(p,q)$-forms is hypoelliptic with loss of one derivative. 

It was proved by Beals-Greiner~\cite{BG:CHM} that for the Kohn Laplacian acting on $(p,q)$-forms the 
condition~(\ref{eq:Intro.sublaplacian.condition})  reduces to the condition 
$Y(q)$, so we may apply Proposition~\ref{prop:Intro.sublaplacian.Rockland-bundle} to recover Kohn's result. This allows us to complete the treatment  of the Kohn 
in~\cite{BG:CHM}, because the initial argument there is not quite complete (see Remark~\ref{rem:Examples.Boxb}). 

\subsubsection*{(c) Horizontal sublaplacian} The horizontal sublaplacian $\Delta_{b}$ on a Heisenberg manifold $(M^{d+1},H)$ can be seen as a horizontal 
Laplacian acting on the horizontal forms, that is on the sections of $\Lambda^{*}_{\C}H^{*}$. This operator was first introduced by Tanaka~\cite{Ta:DGSSPCM} 
for strictly pseudoconvex CR manifolds, but versions of this operator acting on functions 
were independently defined by Greenleaf~\cite{Gr:FESPM} and Lee~\cite{Le:FMPHI}. 

While $\Delta_{b}$ acting on functions can be seen as a sum of squares up to lower order terms, it 
seems that little has been done concerning the hypoellipticity of $\Delta_{b}$ acting on forms on a general Heisenberg manifolds, except in the 
contact case (see~\cite{Ta:DGSSPCM}, \cite{Ru:FDVC}). 

We show here that the invertibility of the principal symbol of $\Delta_{b}$ acting on the sections of $\Lambda^{k}_{\C}H^{*}$  
reduces to a condition involving $k$ and the Levi form of $(M,H)$ only. More precisely, for $a \in M$ let $2n$ be the rank of the Levi form $\cL$ at 
$a$. We say that the condition $X(k)$ is satisfied at $a$ when 
\begin{equation}
     k\not \in\{n,n+1,\ldots,d-n\}.
%     \label{eq:¥}
\end{equation}
We then prove that this condition is the condition~(\ref{eq:Intro.sublaplacian.condition}) for  $\Delta_{b}$ acting on $\Lambda^{k}_{\C}H^{*}$. Thus,  
using Proposition~\ref{prop:Intro.sublaplacian.Rockland-bundle} we get: 

\begin{proposition}\label{prop:Intro.horizontal-sublaplacian}
    Let $\Delta_{b}:C^{\infty}(M,\Lambda^{k}_{\C}H^{*})\rightarrow C^{\infty}(M,\Lambda^{k}_{\C}H^{*})$ be the horizontal  sublaplacian acting on 
    horizontal forms of degree $k$.\smallskip  
   
   1) At a point $x\in M$ the Rockland condition for $\Delta_{b}$ is equivalent to the condition~$X(k)$.\smallskip
   
   2) The principal symbol of $\Delta_{b}$ is invertible if, and only if, the condition $X(k)$ is satisfied at every point. 
%    In particular, when the latter 
%    occurs $\Delta_{b}$ admits a parametrix in $\pvdo^{-2}(M,\Lambda^{k}_{\C}H^{*})$ and is hypoelliptic with loss of 1 derivative.
\end{proposition}

\subsubsection*{(d) Contact Laplacian} Given a contact manifold $(M^{2n+1},\theta)$ the contact Laplacian is associated to the contact complex 
defined by Rumin~\cite{Ru:FDVC}. Unlike the previous examples this not a sublaplacian and it is even of order 4 on contact forms of middle degrees. 

It has been shown by Rumin~\cite{Ru:FDVC} that the contact Laplacian satisfies the Rockland condition on forms of any degree. Rumin then used results 
of Helffer-Nourrigat~\cite{HN:HMOPCV} to deduce that the contact Laplacian was hypoelliptic maximal. Alternatively, we may use 
Proposition~\ref{thm:Intro.Rockland-Parametrix}   to deduce 
that on contact forms of any degree the contact Laplacian has an invertible principal symbol, hence admits a parametrix in the Heisenberg calculus and 
is hypoelliptic. 

\subsection{Organization of the paper} 
The paper is organized as follows.  In Section~\ref{sec.Heisenberg} we recall the main definitions and examples concerning Heisenberg manifolds and their tangent 
Lie group bundles. In Section~\ref{sec:PsiHDO} we give a detailed overview of the Heisenberg 
calculus of~\cite{BG:CHM} and~\cite{Ta:NCMA}, following closely the exposition of~\cite{BG:CHM}. 

The section~\ref{sec:principal-symbol} is devoted to the definitions and the 
main properties of the 
principal symbols and model operators of \psivdos. In Section~\ref{sec:hypoellipticity} we study the relationships between invertibility of the principal symbol, 
Rockland condition and hypoellipticity and, in particular, we prove Theorem~\ref{thm:Intro.Rockland-Parametrix}. 

In Section~\ref{sec:sublaplacian} we deal with sublaplacians and in particular extend the results of~\cite{BG:CHM} to sublaplacians acting on sections of vector 
bundles. In Section~\ref{sec:Examples} we deal with the invertibility of the principal symbols of  main geometric 
operators on Heisenberg manifolds: H\"ormander's sum of squares, Kohn Laplacian, horizontal sublaplacian  and contact Laplacian. In particular, we 
prove Proposition~\ref{prop:Intro.horizontal-sublaplacian}.

The last two sections are devoted to the proofs of Proposition~\ref{prop:PsiHDO.invariance} 
and Proposition~\ref{prop:PsiHDO.transpose-chart}. These are versions in Heisenberg coordinates of the invariance of the Heisenberg calculus 
by  Heisenberg diffeomorphisms and transposition. The latter were proved in~\cite{BG:CHM} in a less precise form, but in this paper we need the  
precise version in Heisenberg coordinates in order to show that the definition of the principal symbol makes sense intrinsically and 
to determine the principal symbols of the transpose and the adjoint of a \psivdo. 

\begin{acknowledgements} 
The author is grateful  to the hospitality of the Mathematics Departments of Princeton 
 University and Harvard University where most part of this paper was written. The research of the author was partially supported 
 by the NSF grant DMS 0409005. 
\end{acknowledgements}

\section{Heisenberg manifolds and their tangent Lie group bundles} 
\label{sec.Heisenberg}
In this section we recall the main facts about Heisenberg manifolds and their tangent Lie group bundle. The exposition here follows closely that of~\cite{Po:Pacific1}.

\begin{definition}
   1) A Heisenberg manifold is a smooth manifold $M$ equipped with a distinguished hyperplane bundle $H \subset TM$.\smallskip 
   
   2) A Heisenberg diffeomorphism $\phi$ from a Heisenberg manifold $(M,H)$ onto another Heisenberg manifold 
   $(M,H')$ is a diffeomorphism $\phi:M\rightarrow M'$ such that $\phi_{*}H = H'$. 
\end{definition}

\begin{definition}
   Let $(M^{d+1},H)$ be a Heisenberg manifold. Then:\smallskip 
   
   1) A (local) $H$-frame for $TM$ is a (local) frame $X_{0}, X_{1}, \ldots, X_{d}$ of $TM$ so that $X_{1}, \ldots, X_{d}$ span~$H$.\smallskip  
   
   2) A local Heisenberg chart is a  local chart together with a local $H$-frame of $TM$ over its domain.
\end{definition}

The main examples of Heisenberg manifolds are the following.\smallskip 

\emph{(a) Heisenberg group}. The $(2n+1)$-dimensional Heisenberg group
$\bH^{2n+1}$ is $\R^{2n+1}=\R \times \R^{2n}$ equipped with the 
group law, 
\begin{equation}
    x.y=(x_{0}+y_{0}+\sum_{1\leq j\leq n}(x_{n+j}y_{j}-x_{j}y_{n+j}),x_{1}+y_{1},\ldots,x_{2n}+y_{2n}).  
\end{equation}
A left-invariant basis for its Lie algebra $\fh^{2n+1}$ is 
provided by the vector-fields, 
\begin{equation}
    X_{0}=\frac{\partial}{\partial x_{0}}, \quad X_{j}=\frac{\partial}{\partial x_{j}}+x_{n+j}\frac{\partial}{\partial 
    x_{0}}, \quad X_{n+j}=\frac{\partial}{\partial x_{n+j}}-x_{j}\frac{\partial}{\partial 
    x_{0}}, \quad 1\leq j\leq n,
     \label{eq:Examples.Heisenberg-left-invariant-basis}
\end{equation}
which  for $j,k=1,\ldots,n$ and $k\neq j$ satisfy the Heisenberg relations,
\begin{equation}
    [X_{j},X_{n+k}]=-2\delta_{jk}X_{0}, \qquad [X_{0},X_{j}]=[X_{j},X_{k}]=[X_{n+j},X_{n+k}]=0.
     \label{eq:Examples.Heisenberg-relations}
\end{equation}
In particular, the subbundle spanned by the vector fields 
$X_{1},\ldots,X_{2n}$ defines a left-invariant Heisenberg structure on 
$\bH^{2n+1}$.\smallskip

 \emph{(b) Foliations.} A (smooth) foliation is a manifold $M$ together with a subbundle $\cF \subset TM$ 
which is integrable in the Froebenius' sense, i.e., we have
$[\cF,\cF]\subset \cF$. Thus, any codimension 1 foliation is a Heisenberg manifold.\smallskip  

 \emph{(c) Contact manifolds}. 
Opposite to foliations are contact manifolds: a \emph{contact
structure} on a manifold $M^{2n+1}$ is given by a global non-vanishing $1$-form $\theta$ on $M$ such that
$d\theta$ is non-degenerate on $H=\ker \theta$. In particular, $(M,H)$ is a Heisenberg manifold. In fact, by
Darboux's theorem any contact manifold $(M^{2n+1},\theta)$ is locally
contact-diffeomorphic to the Heisenberg group $\bH^{2n+1}$ equipped with its standard contact
form $\theta^{0}= dx_{0}+\sum_{j=1}^{n}(x_{j}dx_{n+j}-x_{n+j}dx_{j})$.\smallskip

 \emph{(d) Confoliations}. According to Elyashberg-Thurston~\cite{ET:C} a \emph{confoliation structure} on an oriented manifold
$M^{2n+1}$ is given by a global non-vanishing $1$-form $\theta$ on $M$ such that
$(d\theta)^{n}\wedge \theta\geq 0$. In particular, when $(d\theta)^{n}
\wedge \theta=0$ (resp.~$(d\theta)^{n}\wedge \theta>0$) we are in presence of a foliation (resp.~a contact structure). In any case the hyperplane bundle 
$H=\ker \theta$ defines a Heisenberg structure on $M$.\smallskip

 \emph{(e) CR manifolds.} A CR
structure on an orientable manifold $M^{2n+1}$ is given by a rank $n$
complex subbundle $T_{1,0}\subset T_{\C}M$ which is integrable in  Froebenius' sense and such that 
$T_{1,0}\cap T_{0,1}=\{0\}$, where $T_{0,1}=\overline{T_{1,0}}$. 
Equivalently, the subbundle $H=\Re (T_{1,0}\otimes T_{0,1})$ has the 
structure of a complex bundle of (real) dimension $2n$. In
particular, the pair $(M,H)$ forms a Heisenberg manifold. 

Moreover, since $M$ is orientable and $H$ is orientable by means of its complex structure, the normal bundle $TM/H$ is orientable, hence admits a global 
nonvanishing section $T$. Let $\theta$ be the global section of $T^{*}M/H^{*}$ such that $\theta(T)=1$ and $\theta$ annihilates $H$. Then Kohn's Levi 
form is the form $L_{\theta}$ on $T_{1,0}$ such that, for sections $Z$ and $W$ of $T_{1,0}$, we have
\begin{equation}
    L_{\theta}(Z,W)=-id\theta(Z,\bar{W})=i\theta([Z,\bar{W}]). 
     \label{eq:Heisenberg.Kohn-Levi-form}
\end{equation}

We say that $M$ is strictly pseudoconcex (resp.~nondegenerate, $\kappa$-strictly pseudoconvex) when for some choice of $\theta$  the Levi form 
$L_{\theta}$  is everywhere positive definite (resp.~is everywhere non-degenerate, has everywhere signature $(n-\kappa,\kappa,0)$). In particular, 
when $M$ is nondegenerate the 1-form $\theta$ is a contact form on $M$. 

The main example of a CR manifold is that of the (smooth) boundary $M=\partial D$ of a complex domain $D \subset \C^{n}$. In particular, 
when $D$ is strongly pseudoconvex (or strongly pseudoconcave) $M$ is strictly pseudoconvex.

\subsection{The tangent Lie group bundle}
A simple description of the tangent Lie group bundle of a Heisenberg manifold $(M^{d+1},H)$ can be given as follows.

\begin{lemma}[\cite{Po:Pacific1}]
The Lie bracket of vector fields induces 2-form,
\begin{equation}
    \cL: H\times H \longrightarrow TM/H,
     \label{eq:Bundle.Levi-form1}
\end{equation}
such that, for any $a\in M$ and any sections $X$ and $Y$ of $H$ near $a$, we have
\begin{equation}
    \cL_{a}(X(a),Y(a)) = [X,Y](a) \quad \bmod H_{a}.
     \label{eq:Bundle.Levi-form2}
\end{equation}
\end{lemma}

\begin{definition}
 The $2$-form  $\cL$ is called the Levi form of $(M,H)$.
\end{definition}

The Levi form $\cL$ allows us to define a bundle $\fg M$ of graded Lie algebras  by endowing $(TM/H)\oplus H$ 
with the smooth fields of Lie Brackets and gradings such that
\begin{equation}
    [X_{0}+X',Y_{0}+Y']_{a}=\cL_{a}(X',Y') \qquad \text{and} \qquad t.(X_{0}+X')=t^{2}X_{0}+tX', \quad t \in \R,
    \label{eq:Heisenberg.intrinsic-Lie-algebra-structure}
\end{equation}
for $a\in M$ and $X_{0}$, $Y_{0}$ in $T_{a}M/H_{a}$ and $X'$, $Y'$ in $H_{a}$. 

\begin{definition}
    The bundle $\fg M$ is called the tangent Lie algebra bundle of $M$.
\end{definition}

As we can easily check $\fg M$ is a bundle of $2$-step nilpotent Lie algebras which contains the normal bundle $TM/H$ in its center.
Therefore, its associated 
graded Lie group bundle $GM$ can be described as follows. As a bundle $GM$ is $(TM/H)\oplus H$ and the exponential 
map is the identity. In particular, the grading of $GM$ is as in~(\ref{eq:Heisenberg.intrinsic-Lie-algebra-structure}). 
Moreover, since  $\fg M$ is 
2-step nilpotent the Campbell-Hausdorff formula gives 
\begin{equation}
    (\exp X)(\exp Y)= \exp(X+Y+\frac{1}{2}[X,Y]), \qquad \text{$X$, $Y$ sections of $\fg M$}.
\end{equation}
From this we deduce that the product on $GM$ is such that  
\begin{equation}
    (X_{0}+X').(Y_{0}+X')=X_{0}+Y_{0}+\frac{1}{2}\cL(X',Y')+X'+Y',    
    \label{eq:Bundle.Lie-group-law}
\end{equation}
for sections $X_{0}$, $Y_{0}$ of $TM/H$  and sections $X'$, $Y'$ of $H$.

\begin{definition}
    The bundle  $GM$ is called the tangent Lie group bundle of $M$. 
\end{definition}

In fact, the fibers of $GM$ are classified by the Levi form $\cL$ as follows.

\begin{proposition}[\cite{Po:Pacific1}]\label{prop:Bundle.intrinsic.fiber-structure}
  1) Let $a\in M$. Then $\cL_{a}$ has rank $2n$ if, and only if, as a 
  graded Lie group $G_{a}M$ is isomorphic to $\bH^{2n+1}\times \R^{d-2n}$.\smallskip 
  
  2) The Levi form $\cL$ has constant rank $2n$ if, and only if, $GM$ is  a fiber bundle with typical fiber 
  $\bH^{2n+1}\times \R^{d-2n}$.
\end{proposition}

Now, let $\phi:(M,H)\rightarrow (M',H')$ be a Heisenberg diffeomorphism from $(M, H)$ onto another Heisenberg manifold 
$(M',H')$. Since $\phi_{*}H=H'$ we see that $\phi'$ induces a smooth vector bundle isomorphism 
$\overline{\phi}:TM/H\rightarrow TM'/H'$. 

\begin{definition}
We let  $\phi_{H}':(TM/H)\oplus 
  H \rightarrow (TM'/H')\oplus H'$ denote the vector bundle isomorphism such that
    \begin{equation}
    \phi'_{H}(a)(X_{0}+X')=\overline{\phi}'(a)X_{0}+\phi'(a)X',
     \label{eq:Bundle.Intrinsic.Phi'H}
\end{equation}
for any $a\in M$ and any $X_{0}\in T_{a}/H_{a}$ and $X'\in H_{a}$.
\end{definition}

\begin{proposition}[\cite{Po:Pacific1}]\label{prop:Bundle.Intrinsic.Isomorphism}
The vector bundle isomorphism  $\phi'_{H}$ is an isomorphism of graded Lie group bundles from $GM$ onto $GM'$. In particular, the Lie group bundle isomorphism 
class of $GM$ depends only on the Heisenberg diffeomorphism class of $(M,H)$.  
\end{proposition}

\subsection{Heisenberg coordinates and nilpotent approximation of vector fields}
It is interesting to relate the intrinsic description of $GM$ above with the more extrinsic description of~\cite{BG:CHM} (see also~\cite{Be:TSSRG}, 
\cite{EMM:HAITH}, \cite{FS:EDdbarbCAHG},  \cite{Gr:CCSSW}, \cite{Ro:INA}) in terms of the Lie group 
associated to a nilpotent Lie algebra of model vector fields. 

First, let $a\in M$ and let us describe $\fg_{a}M$ as the graded Lie algebra of left-invariant vector fields on $G_{a}M$  
by identifying any $X \in \fg_{a}M$ with the left-invariant vector fields $L_{X}$ on $G_{a}M$ given by 
\begin{equation}
    L_{X}f(x)= \frac{d}{dt}f[(t\exp X).x]_{|_{t=0}}= \frac{d}{dt}f[(tX).x]_{|_{t=0}}, \qquad f \in C^{\infty}(G_{a}M).
\end{equation}
This allows us to associate to any vector fields $X$ near $a$ a unique left-invariant vector fields $X^{a}$ on $G_{a}M$ 
such that 
\begin{equation}
    X^{a}= \left\{ 
    \begin{array}{ll}
        L_{X_{0}(a)} & \text{if $X(a)\not \in H_{a}$},  \\
        L_{X(a)} & \text{otherwise,} 
    \end{array}\right.
     \label{eq:Bundle.intrinsic.model-vector-fields}
\end{equation}
where $X_{0}(a)$ denotes the class of $X(a)$ modulo $H_{a}$. 

\begin{definition}
    The left-invariant vector fields $X^{a}$ is called the model vector fields of $X$ at $a$.
\end{definition}

Let us look at the above construction in terms of a $H$-frame $X_{0},\ldots,X_{d}$ near 
$a$, i.e.~of a local trivialization of the vector bundle $(TM/H)\oplus H$. For $j,k=1,\ldots,d$ we let 
\begin{equation}
    \cL(X_{j},X_{k})=[X_{j},X_{k}]X_{0}=L_{jk}X_{0} \quad \bmod H.
\end{equation}
With respect to the coordinate system $(x_{0},\ldots,x_{d})\rightarrow x_{0}X_{0}(a)+\ldots+x_{d}X_{d}(a)$ we can 
write the product law of $G_{a}M$ as 
\begin{equation}
    x.y=(x_{0}+\frac{1}{2}\sum_{j,k=1}^{d}L_{jk}x_{j}x_{k},x_{1},\ldots,x_{d}).
     \label{eq:Heisenberg.productGmM-coordinates}
\end{equation}
Then the vector fields $X_{j}^{a}$, $j=0,\ldots,d$, in~(\ref{eq:Bundle.intrinsic.model-vector-fields}) 
are the left-invariant vector fields corresponding to the vectors $e_{j}$, $j=0,\ldots,d$, of the canonical basis 
of $\Rd$, i.e., we have
\begin{equation}
    X_{0}^{a}=\frac{\partial}{\partial x_{0}} \quad \text{and}  \quad X_{j}^{a}=\frac{\partial}{\partial x_{j}} 
    -\frac{1}{2}\sum_{k=1}^{d}L_{jk}x_{k}\frac{\partial}{\partial x_{0}}, \quad 1\leq j\leq d.
     \label{eq:Heisenberg.Xjm.coordinates}
\end{equation}
In particular, for $j,k=1,\ldots,d$ we have the relations, 
\begin{equation}
    [X_{j}^{a},X_{k}^{a}]=L_{jk}(a)X_{0}^{a}, \qquad [X_{j}^{a},X_{0}^{a}]=0.
     \label{eq:Heisenberg.constant-structures.Gm}
\end{equation}

Now, let $\kappa:\dom \kappa \rightarrow U$ be a Heisenberg chart near $a=\kappa^{-1}(u)$ and let 
$X_{0},\ldots,X_{d}$ be the associated $H$-frame of $TU$.  
Then there exists a unique affine coordinate change $x \rightarrow \psi_{u}(x)$ such that 
$\psi_{u}(u)=0$ and $\psi_{u*}X_{j}(0)=\frac{\partial}{\partial x_{j}}$ for 
$j=0,1,\ldots,d$. Indeed, if for $j=1,\ldots,d$ we set $X_{j}(x)=\sum_{k=0}^{d}B_{jk}(x)\frac{\partial}{\partial x_{k}}$ then 
we have
\begin{equation}
    \psi_{u}(x)=A(u)(x-u), \qquad A(u)=(B(u)^{t})^{-1}.
\end{equation}

\begin{definition}\label{def:Heisenberg.extrinsic.u-coordinates}
1) The coordinates provided by $\psi_{u}$ are called the privileged coordinates at $u$ 
with respect to the $H$-frame $X_{0},\ldots,X_{d}$. 

2) The map $\psi_{u}$ is called the privileged-coordinate map with respect to the $H$-frame $X_{0},\ldots,X_{d}$.
\end{definition}
\begin{remark}
    The privileged coordinates at $u$ are called $u$-coordinates in~\cite{BG:CHM}, but they correspond to the privileged coordinates 
    of~\cite{Be:TSSRG} and \cite{Gr:CCSSW} in the special case of a Heisenberg manifold. 
\end{remark}

Next, on $\Rd$ we consider the dilations 
\begin{equation}
    \delta_{t}(x)=t.x=(t^{2}x_{0},tx_{1}, \ldots, tx_{d}), \qquad t \in \R,
    \label{eq:Heisenberg.dilations}
\end{equation}
with respect to which $\frac{\partial}{\partial_{x_{0}}}$ is 
homogeneous of degree $-2$ and $\frac{\partial}{\partial_{x_{1}}},\ldots,\frac{\partial}{\partial_{x_{d}}}$ 
is homogeneous of degree~$-1$. 

Since in the privileged coordinates at $u$ we have $X_{j}(0)=\frac{\partial}{\partial x_{j}}$ we can write
\begin{equation}
    X_{j}= \frac{\partial}{\partial_{x_{j}}}+ \sum_{k=0}^{d} a_{jk}(x) \frac{\partial}{\partial_{x_{k}}}, 
    \qquad j=0,1,\ldots d,
\end{equation}
where the $a_{jk}$'s are smooth functions such that $a_{jk}(0)=0$. Therefore, we may define
\begin{gather}
    X_{0}^{(u)}= \lim_{t\rightarrow 0} t^{2}\delta_{t}^{*}X_{0}= \frac{\partial}{\partial_{x_{0}}},
    \label{eq:Heisenberg.X0u}\\
     X_{j}^{(u)}= \lim_{t\rightarrow 0} t^{-1}\delta_{t}^{*}X_{j}= 
     \frac{\partial}{\partial_{x_{j}}}+\sum_{k=1}^{d}b_{jk}x_{k} \frac{\partial}{\partial_{x_{0}}}, \quad 
     j=1,\ldots,d, \label{eq:Heisenberg.Xju}     
\end{gather}
where for $j,k=1,\ldots,d$ we have set $b_{jk}= \partial_{x_{k}}a_{j0}(0)$. 

 Observe that $X_{0}^{(u)}$ is homogeneous of degree $-2$ and $X_{1}^{(u)},\ldots,X_{d}^{(u)}$ are homogeneous of degree $-1$.
 Moreover, for $j,k=1,\ldots,d$ we have 
\begin{equation}
    [X_{j}^{(u)},X_{0}^{(u)}]=0 \quad \text{and} \quad [X_{j}^{(u)},X_{0}^{(u)}]=(b_{kj}-b_{jk})X_{0}^{(u)}.
     \label{eq:Heisenberg.constant-structures.Gu1}
\end{equation}
Thus, the linear space spanned by $X_{0}^{(u)},X_{1}^{(u)}, \ldots, X_{d}^{(u)}$ is a graded 2-step nilpotent 
Lie algebra $\fg^{(u)}$. In particular, $\fg^{(u)}$ is the Lie algebra of left-invariant vector fields over the graded Lie group $G^{(u)}$ 
consisting of $\Rd$ equipped with the grading~(\ref{eq:Heisenberg.dilations}) and the group law,
\begin{equation}
    x.y=(x_{0}+\sum_{j,k=1}^{d}b_{kj}x_{j}x_{k},x_{1},\ldots,x_{d}).
\end{equation}

Now, if near $a$ we let $\cL(X_{j},X_{k})=[X_{j},X_{k}]=L_{jk}(x)X_{0}\bmod H$, then we get
\begin{equation}
    [X_{j}^{(u)},X_{k}^{(u)}]=\lim_{t\rightarrow 0}[t\delta_{t}^{*}X_{j},t\delta_{t}^{*}X_{k}] = 
    \lim_{t\rightarrow 0} t^{2}\delta_{t}^{*}(L_{jk}(\circ \kappa^{-1}(x))X_{0})=L_{jk}(a)X_{0}^{(u)}.
     \label{eq:Heisenberg.constant-structures.Gu2}
\end{equation}
Comparing this with~(\ref{eq:Heisenberg.constant-structures.Gm}) and~(\ref{eq:Heisenberg.constant-structures.Gu1}) 
then shows that $\fg^{(u)}$ has the same the constant structures as those of 
$\fg_{a}M$, hence is isomorphic to it. Consequently, the Lie groups $G^{(u)}$ and $G_{a}M$ are isomorphic. 
In fact, as shown in~\cite{BG:CHM} and~\cite{Po:Pacific1}, an explicit isomorphism is given by 
\begin{equation}
     \phi_{u}(x_{0},\ldots,x_{d})= (x_{0}-\frac{1}{4}\sum_{j,k=1}^{d}(b_{jk}+b_{kj})x_{j}x_{k},x_{1},\ldots,x_{d}).
     \label{eq:Bundle.Extrinsic.Phiu}
\end{equation}

\begin{definition}\label{def:Bundle.extrinsic.normal-coordinates}
Let $\varepsilon_{u}=\phi_{u}\circ \psi_{u}$. Then:\smallskip 

1) The new coordinates provided by $\varepsilon_{u}$  are called Heisenberg 
coordinates at $u$ with respect to the $H$-frame $X_{0},\ldots,X_{d}$.\smallskip  

2) The map $\varepsilon_{u}$ is called the $u$-Heisenberg coordinate map.
\end{definition}

\begin{remark}
       The Heisenberg coordinates at $u$ have been also considered in~\cite{BG:CHM} as a technical tool 
       for inverting the principal symbol of a hypoelliptic sublaplacian.
\end{remark}

  Next, by~\cite[Lem.~1.17]{Po:Pacific1} we have
\begin{equation}
    \phi_{*}X_{0}^{(u)}=\frac{\partial}{\partial x_{0}}=X_{0}^{a} \quad \text{and} \quad 
    \phi_{*}X_{j}^{(u)}=\frac{\partial}{\partial x_{j}}-\frac{1}{2}\sum_{k=1}^{d}L_{jk}x_{k}\frac{\partial}{\partial x_{0}}=X_{j}^{a}, \quad 
    j=1,\ldots,d. 
    \label{eq:Heisenberg.Xu-Xm}
\end{equation}
Since $\phi_{u}$ commutes with the Heisenberg dilations~(\ref{eq:Heisenberg.dilations}),  
using~(\ref{eq:Heisenberg.X0u})--(\ref{eq:Heisenberg.Xju}) we get
\begin{equation}
    \lim_{t\rightarrow 0} t^{2}\delta_{t}^{*}\phi_{u*}X_{0}^{(u)}=X^{a}_{0} \quad \text{and} \quad 
    \lim_{t\rightarrow 0} t\delta_{t}^{*}\phi_{u*}X_{j}^{(u)}=X^{a}_{j}, \quad j=1,\ldots,d.
\end{equation}

In fact, as shown in~\cite{Po:Pacific1}, in Heisenberg coordinates at $a$ 
for any vector fields $X$ as $t\rightarrow 0$  we have
\begin{equation}
   \delta_{t}^{*}X=  \left\{ 
   \begin{array}{ll}
       t^{-2}X^{a} +\op{O}(t^{-1})& \text{if $X(a)\in H_{a}$},\\
       t^{-1}X^{a} +\op{O}(1) & \text{otherwise}. 
   \end{array}\right. 
   \label{eq:Bundle.Extrinsic.approximation-normal}
\end{equation}
Therefore, we obtain:

\begin{proposition}[\cite{Po:Pacific1}]\label{prop:Bundle.equivalent-descriptions}
    In the Heisenberg coordinates centered at $m=\kappa^{-1}(u)$ the tangent Lie group $G_{a}M$ coincides with $G^{(u)}$ and for any vector fields 
    $X$ the model vector fields $X^{a}$ approximates $X$ near $a$ in the sense of~(\ref{eq:Bundle.Extrinsic.approximation-normal}).
\end{proposition}

One consequence of the equivalence between the two approaches to $GM$ is a tangent approximation for Heisenberg diffeomorphisms as follows. 

Let  $\phi:(M,H)\rightarrow (M',H')$ be a Heisenberg diffeomorphism
from $(M,H)$ to another Heisenberg manifold $(M',H')$. We also endow $\Rd$ with the pseudo-norm,
\begin{equation}
    \|x\|= (x_{0}^{2}+(x_{1}^{2}+\ldots+x_{d}^{2})^{2})^{1/4}, \qquad x\in \Rd,
\end{equation}
so that, for any $x \in \Rd$ and any $t \in \R$, we have 
\begin{equation}
    \|t.x\|=|t|\, \|x\| . 
     \label{eq:Bundle.homogeneity-pseudonorm}
\end{equation}

\begin{proposition}[{\cite[Prop.~2.21]{Po:Pacific1}}]\label{prop:Heisenberg.diffeo}
   Let $a\in M$ and set $m'=\phi(a)$. Then, in Heisenberg coordinates at $a$ and at $a'$   
   the diffeomorphism $\phi(x)$ has a behavior near 
   $x=0$ of the form 
   \begin{equation}
       \phi(x)= \phi_{H}'(0)x+(\op{O}(\|x\|^{3}), \op{O}(\|x\|^{2}),\ldots,\op{O}(\|x\|^{2})).
        \label{eq:Bundle.Approximation-diffeo}
   \end{equation}
   In particular, there are no terms of the form $x_{j}x_{k}$, $1\leq j,k\leq d$, in the Taylor expansion of $\phi_{0}(x)$ 
   at $x=0$.
\end{proposition}
   \begin{remark}
  An asymptotics similar to~(\ref{eq:Bundle.Approximation-diffeo}) is given  
  in~\cite[Prop.~5.20]{Be:TSSRG} in privileged coordinates at $u$ and 
   $u'=\kappa_{1}(a')$, 
  but the leading term there is only a Lie algebra isomorphism from $\fg^{(u)}$ onto $\fg^{(u')}$. This is only 
  in Heisenberg coordinates that we recover the Lie group isomorphism $\phi'_{H}(a)$ as the leading term of the asymptotics.
\end{remark}
   
% \begin{remark}
%     An interesting application  of Proposition~\ref{prop:Heisenberg.diffeo} in~\cite{Po:Pacific1} is the construction of the tangent groupoid 
%     $\cG_{H}M$ of $(M,H)$ as the differentiable groupoid encoding the smooth deformation of  $M\times M$ to $GM$. 
%    This groupoid is the analogue in the Heisenberg setting of Connes' tangent groupoid of a manifold~(\cite[II.5]{Co:NCG}, \cite{HS:MKOEFFTK}) and its shows 
%    that $GM$ is tangent to $a$ in a differentiable  fashion (compare~\cite{Be:TSSRG}, \cite{Gr:CCSSW}). 
% \end{remark}

\section{Hypoelliptic calculus on Heisenberg manifolds}
\label{sec:PsiHDO}
The Heisenberg calculus is the relevant pseudodifferential tool to study hypoelliptic operators  on Heisenberg 
manifolds. It was independently invented by Beals-Greiner~\cite{BG:CHM} and Taylor~\cite{Ta:NCMA}, extending previous 
works of Boutet de Monvel~\cite{BdM:HODCRPDO}, Folland-Stein~\cite{FS:EDdbarbCAHG} and Dynin~(\cite{Dy:POHG}, \cite{Dy:APOHSC}) 
 (see also~\cite{BdMGH:POPDCM}, \cite{CGGP:POGD}, \cite{EMM:HAITH}, \cite{Gr:HPOPCD}, \cite{Ho:CHPODC}, \cite{RS:HDONG}). 

 The idea in the Heisenberg calculus is to have a pseudodifferential calculus on a Heisenberg manifold $(M,H)$ which is modeled at any point 
 $a\in M$ by the calculus of left-invariant pseudodifferential operators on the tangent group $G_{a}M$. 

 \subsection{Left-invariant pseudodifferential operators}
 Let $(M^{d+1},H)$ be a Heisenberg manifold and let $G=G_{a}M$ be the tangent Lie group of $M$ at a point $a \in M$. We recall here the main facts about 
left-invariant pseudodifferential operators on $G$ (see also~\cite{BG:CHM}, \cite{CGGP:POGD}, \cite{Ta:NCMA}). 

 Recall that for any finite dimensional  vector space $E$ the Schwartz class $\cS(E)$ is a Fr\'echet space and the Fourier transform is the continuous 
 isomorphism of $\cS(E)$ onto $\cS(E^{*})$ given by
 \begin{equation}
     \hat{f}(\xi)=\int_{E}e^{i\acou{\xi}{x}}f(x)dx, \qquad f \in \cS(E), \quad \xi \in E^{*},
 \end{equation}
 where $dx$ denotes the Lebesgue measure of $E$. 
 \begin{definition}
     $\cS_{0}(E)$ is the closed subspace of $\cS(E)$ consisting of $f \in \cS(E)$ such that for any differential operator $P$ on $E^{*}$ we have $(P\hat{f})(0)=0$. 
 \end{definition}
 
 Since $G$ has the same underlying set as that of its Lie algebra $\fg=\fg_{x}M$ we can let $\cS(G)$ and $\cS_{0}(G)$ denote the Fr\'echet spaces 
 $\cS(E)$ and $\cS_{0}(E)$ associated to the underlying linear space $E$ of $\fg$ 
 (notice that the Lebesgue measure of $E$ coincides with the Haar measure of $G$ since $G$ is nilpotent). 
 
 Next, for  $\lambda \in \R$ and $\xi=\xi_{0}+\xi'$ in $\fg^{*}= (T_{a}^{*}M/H^{*}_{a})\oplus H_{a}$ we let 
  \begin{equation}
     \lambda.\xi=  \lambda.(\xi_{0}+\xi')=\lambda^{2}\xi_{0}+\lambda \xi'.
      \label{eq:PsiHDO.Heisenberg-dilation-fg*}
 \end{equation}

 \begin{definition}
     $S_{m}(\fg^{*})$, $m \in \C$, is the space of functions $p \in C^{\infty}(\fg^{*}\setminus 0)$ which are homogeneous of degree $m$, in the sense 
     that, for any $\lambda>0$ we have 
\begin{equation}
         p(\lambda.\xi)=\lambda^{m}p(\xi) \qquad \forall \xi\in \fg^{*}\setminus 0. 
\end{equation}
    In addition $S_{m}(\fg^{*})$ is endowed with the Fr\'echet space topology induced from that of $C^{\infty}(\fg^{*}\setminus 0)$. 
 \end{definition}

 Note that the image $\hat{\cS}_{0}(G)$ of $\cS(G)$ under the Fourier transform consists of functions $v\in \cS(\fg^{*})$ such that, given any norm 
 $|.|$ on $G$,  near $\xi=0$ we have $|g(\xi)|=\op{O}(|\xi|^{N})$ for any integer $N\geq 0$. Thus, any $p \in S_{m}(\fg^{*})$ defines an element of
 $\hat{\cS}_{0}(\fg^{*})'$ by letting 
 \begin{equation}
     \acou{p}{g}= \int_{\fg^{*}}p(\xi)g(\xi)d\xi, \qquad g \in \hat{\cS}_{0}(\fg^{*}).
 \end{equation}
 This allows us to define the inverse Fourier transform of $p$ as the element $\check{p}\in \cS_{0}(G)'$ such that 
 \begin{equation}
     \acou{\check{p}}{f}=\acou{p}{\check{f}} \qquad \forall f \in \cS_{0}(G).
       \label{eq:PsiHDO.inverse-Fourier-transform-symbol}
 \end{equation}

\begin{proposition}[\cite{BG:CHM}, \cite{CGGP:POGD}]\label{prop:PsiHDO.convolution-symbols-group}
    1) For any $p \in S_{m}(\fg^{*})$ the left-convolution by $\check{p}$, 
\begin{equation}
   \check{p}*f(x):=\acou{\check{p}(y)}{f(x.y^{-1})}, \qquad f \in \cS_{0}(G), 
     \label{eq:PsiHDO.convolution-operator}
\end{equation}
defines a continuous endomorphism of $\cS_{0}(G)$.\smallskip 

   2) There is a continuous bilinear product,
 \begin{equation}
     *:S_{m_{1}}(\fg^{*})\times S_{m_{2}}(\fg^{*}) \longrightarrow S_{m_{1}+m_{2}}(\fg^{*}),
 \end{equation} 
such that, for any  $p_{1}\in S_{m_{1}}(\fg^{*})$ and $p_{2}\in S_{m_{2}}(\fg^{*})$, the composition of the left-convolution operators by 
$\check{p_{1}}$ and $\check{p_{2}}$ is the left-convolution operator by $(p_{1}*p_{2})^{\vee}$, that is, we have
\begin{equation}
    \check{p_{1}}*(\check{p_{2}}*f)=(p_{1}*p_{2})^{\vee}*f  \qquad \forall f \in \cS_{0}(G) .
\end{equation}
\end{proposition}

Let us also mention that if $p \in S_{m}(\fg^{*})$ then the convolution operator $Pu=\check{p}*f$ is a pseudodifferential operator. 
Indeed, let $X_{0}(a),\ldots,X_{d}(a)$ be a (linear) basis of $\fg$ so that $X_{0}(a)$ is in $T_{a}M/H_{a}$ and 
$X_{1}(a),\ldots,X_{d}(a)$ span $H_{a}$. For $j=0,\ldots,d$ let $X_{j}^{a}$ be the left-invariant vector fields on $G$ such that 
$X^{w}_{j|_{x=0}}=X_{j}(a)$. The basis 
$X_{0}(a),\ldots,X_{d}(a)$ yields a linear isomorphism $\fg\simeq \Rd$, hence a global chart of $G$. In this chart $p$ is a 
homogeneous symbol on $\Rd\setminus 0$ with respect to the dilations 
\begin{equation}
    \lambda.x=(\lambda^{2}x_{0},\lambda x_{1},\ldots,\lambda x_{d}), \qquad x\in \Rd, \quad \lambda>0.
     \label{eq:PsiHDO.Heisenberg-dilations-Rd}
\end{equation}

Similarly, each vector fields $\frac{1}{i}X_{j}^{a}$, $j=0,\ldots,d$, corresponds to a vector fields on $\Rd$ whose symbol is denoted 
$\sigma_{j}^{a}(x,\xi)$. Then, setting $\sigma=(\sigma_{0},\ldots,\sigma_{d})$, it can be shown that in the above chart the operator $P$ is 
given by
\begin{equation}
    Pf(x)=\int_{\Rd} e^{ix.\xi}p(\sigma^{a}(x,\xi))\hat{f}(\xi), \qquad f \in \cS_{0}(\Rd).
     \label{eq:PsiHDO.PsiDO-convolution}
\end{equation}
In other words $P$ is the pseudodifferential operator $p(-iX^{a}):=p(\sigma^{a}(x,D))$ acting on $\cS_{0}(\Rd)$.

\subsection{$\mathbf{\Psi_{H}}$DO's on an open subset of $\Rd$} 
Let $U$ be an open subset of $\Rd$ together with a hyperplane bundle $H \subset TU$ and a $H$-frame $X_{0},X_{1},\ldots,X_{d}$ of $TU$. 
Then the class of \psivdos\ on $U$  is a class of pseudodifferential operators modelled on that of homogeneous convolution operators on the fibers of $GU$.

\begin{definition} $S_{m}(\URd)$, $m\in\C$, is the space of symbols 
    $p(x,\xi)\in C^{\infty}(U\times\Rdo)$ that are homogeneous of degree $m$ with respect to the $\xi$-variable, that is, 
    \begin{equation}
        p(x,\lambda.\xi)=\lambda^m p(x,\xi) \qquad \text{for any $\lambda>0$},
    \end{equation}
    where $\xi \rightarrow \lambda.\xi$ denotes the Heisenberg dilation~(\ref{eq:PsiHDO.Heisenberg-dilations-Rd}).
\end{definition}

Observe that the homogeneity of $p\in S_{m}(\URd)$ implies that, for any compact $K \subset U$, it satisfies the estimates
\begin{equation}
     | \partial^\alpha_{x}\partial^\beta_{\xi}p(x,\xi)|\leq C_{K\alpha\beta}\|\xi\|^{\Re m-\brak\beta}, \qquad x\in K, \quad \xi \neq 0,
     \label{eq:PsiVDO.estimates-homogeneous-symbols}
\end{equation}
where $ \|\xi\|=(|\xi_{0}|^{2}+|\xi_{1}|^{4}+\ldots +|\xi_{d}|^{4})^{1/4}$ and $\brak\alpha = 2\alpha_{0}+\alpha_{1}+\ldots+ 
\alpha_{d}$. 

\begin{definition}$S^m(\URd)$,  $m\in\C$, consists of symbols  $p\in C^{\infty}(\URd)$ with
an asymptotic expansion $ p \sim \sum_{j\geq 0} p_{m-j}$, $p_{k}\in S_{k}(\URd)$, in the sense that, for any integer $N$ and 
for any compact $K \subset U$, we have
\begin{equation}
    | \partial^\alpha_{x}\partial^\beta_{\xi}(p-\sum_{j<N}p_{m-j})(x,\xi)| \leq 
    C_{\alpha\beta NK}\|\xi\|^{\Re m-\brak\beta -N}, \qquad  x\in K, \quad \|\xi \| \geq 1.
    \label{eq:PsiVDO.asymptotic-expansion-symbols}
\end{equation}
\end{definition}

Next, for $j=0,\ldots,d$ let  $\sigma_{j}(x,\xi)$ denote the symbol of $\frac{1}{i}X_{j}$ (in the 
classical sense) and set  $\sigma=(\sigma_{0},\ldots,\sigma_{d})$. For any $p \in S^{m}(\URd)$ it can be shown that the symbol 
$p_{\sigma}(x,\xi):=p(x,\sigma(x,\xi))$ is in the H\"ormander class of symbols of type $(\frac{1}{2}, \frac{1}{2})$ (see~\cite[Prop.~10.22]{BG:CHM}). 
Therefore, we define a continuous linear operator from $C^{\infty}_{c}(U)$ to $C^{\infty}(U)$ by letting 
    \begin{equation}
          p(x,-iX)f(x)= (2\pi)^{-(d+1)} \int e^{ix.\xi} p(x,\sigma(x,\xi))\hat{f}(\xi)d\xi,
    \qquad f\in C^{\infty}_{c}(U).
    \end{equation}

    In the sequel we let $\Psi^{-\infty}(U)$ denotes the class of smoothing operators, i.e.~of operators given by smooth kernels. 
    
\begin{definition}
   $\pvdo^{m}(U)$, $m\in \C$, consists of operators $P:C^{\infty}_{c}(U)\rightarrow C^{\infty}(U)$ of the form
\begin{equation}
         P= p(x,-iX)+R,
\end{equation}
 with $p$ in $S^{m}(\URd)$, called the symbol of $P$, and with $R$ in $\Psi^{-\infty}(U)$.
\end{definition}

    The above definition of the symbol of $P$ differs from that of~\cite{BG:CHM}, since 
    there the authors defined it to be $p_{\sigma}(x,\xi)=p(x,\sigma(x,\xi))$. Note also that $p$ is unique 
    modulo $S^{-\infty}(\URd)$. 

\begin{lemma}\label{lem:PsiHDO.asymptotic-completeness}
   For $j=0,1,\ldots$ let $p_{m-j}\in S_{m-j}(\URd)$. Then there exists $P\in \pvdo^{m}(U)$ with symbol $p\sim \sum_{j\geq 0}p_{m-j}$. Moreover, the 
   operator $P$ is unique modulo smoothing operators. 
\end{lemma}
   
    The class $\pvdo^{m}(U)$ does not depend on the choice of the $H$-frame $X_{0}, \ldots, X_{d}$ (see~\cite[Prop.~10.46]{BG:CHM}). Moreover, 
since it is contained in the class 
of \psidos\ of type $(\frac{1}{2},\frac{1}{2})$ we get:
 
\begin{proposition}\label{prop:PsiHDO.Sobolev-regularity}
    Let $P$ be a \psivdo\ of order $m$ on $U$.\smallskip 
    
    1)  $P$ extends to a continuous linear mapping from 
    $\cE'(U)$ to $\cD'(U)$ and has a distribution kernel which is smooth off the diagonal. \smallskip 
    
    3) Let $k=\Re m$ if $\Re m\geq 0$ and $k=\frac{1}{2}\Re m$ otherwise. Then for any $s \in \R$ the operator $P\in \pvdo^{m}(U)$ extends
    to a continuos mapping from $L^{2}_{s,\op{comp}}(U)$ to $L^{2}_{s-k,\op{loc}}(U)$. 
\end{proposition}

\subsection{Composition of $\mathbf{\Psi_{H}}$DO's}
Recall that there is no symbolic calculus for \psidos\ of type $(\frac{1}{2},\frac{1}{2})$ since the product of two such \psidos\ needs not be 
again a \psido\ of type $(\frac{1}{2},\frac{1}{2})$. However, the fact that the \psivdos\ are modelled on left-invariant pseudodifferential operators 
allows us  to construct a symbolic calculus for \psivdos.

First, for $j=0,\ldots,d$ let $X_{j}^{(x)}$ be the leading homogeneous part of $X_{j}$ in privileged coordinates centered at $x$ 
defined according to~(\ref{eq:Heisenberg.X0u})--(\ref{eq:Heisenberg.Xju}). These vectors span a nilpotent Lie algebra of left-invariant vector fields
on a nilpotent graded Lie group $G^{x}$ which corresponds 
to $G_{x}U$ by pulling back the latter from the Heisenberg coordinates at $x$ to the privileged coordinates at $x$. 

As alluded to above the product law of $G^{(x)}$ defines a convolution product for symbols, 
\begin{equation}
    *^{(x)}: S_{m_{1}}(\Rd) \times S_{m_{2}}(\Rd) \longrightarrow S_{m_{1}+m_{2}}(\Rd).
     \label{eq:PsiHDO.convolution-symbol-pointwise}
\end{equation}
such that, with the notations of~(\ref{eq:PsiHDO.PsiDO-convolution}), on $\cL(\cS_{0}(\Rd))$ we have  
\begin{equation}
    p_{1}(-iX^{(x)})p_{2}(-iX^{(x)})=(p_{1}*^{(x)}p_{2})(-iX^{(x)}) \qquad \forall p_{j}\in S_{m_{j}}(\Rd).
\end{equation}

As it turns out the product $*^{(x)}$ depends smoothly on $x$ (see~\cite[Prop.~13.33]{BG:CHM}). Therefore, we get a continuous bilinear product,
\begin{gather}
    *: S_{m_{1}}(\URd) \times S_{m_{2}}(\URd) 
        \rightarrow S_{m_{1}+m_{2}} (\URd),\\   
        p_{1}*p_{2}(x,\xi)=(p_{1}(x,.)*^{(x)}p_{2}(x,.))(\xi), \qquad p_{j}\in S_{m_{j}}(\URd).
 \label{eq:PsiHDO.convolution-symbols-URd}
\end{gather}

\begin{proposition}[{\cite[Thm.~14.7]{BG:CHM}}] \label{prop:PsiHDO.composition}
    For $j=1,2$ let $P_{j}\in \pvdo^{m_{j}}(U)$ have  symbol $p_{j}\sim \sum_{k\geq 0} p_{j,m_{j}-k}$ and assume  
    that one of these operators is properly supported. Then the operator  
$P=P_{1}P_{2}$ is a \psivdo\ of order $m_{1}+m_{2}$ and has symbol  $p\sim \sum_{k\geq 0} p_{m_{1}+m_{2}-k}$, with  

\begin{equation}
     p_{m_{1}+m_{2}-k}(x,\xi) = \sum_{k_{1}+k_{2}\leq k} \sum_{\alpha,\beta,\gamma,\delta}^{(k-k_{1}-k_{2})}
            h_{\alpha\beta\gamma\delta} (x)  (D_{\xi}^\delta p_{1,m_{1}-k_{1}})* (\xi^\gamma 
            \partial_{x}^\alpha \partial_{\xi}^\beta p_{2,m_{2}-k_{2}})(x,\xi),    
\end{equation}
where $\underset{\alpha\beta\gamma\delta}{\overset{\scriptstyle{(l)}}{\sum}}$ denotes the sum over the indices such that 
$|\beta|=|\gamma|$ and $|\alpha|+|\beta| \leq \brak\beta -\brak\gamma+\brak\delta = l$, and the functions 
$h_{\alpha\beta\gamma\delta}(x)$'s are  polynomials in the derivatives of the coefficients of 
the vector fields $X_{0},\ldots,X_{d}$.
\end{proposition}

% \begin{remark}
% It follows from~(\ref{eq:PsiHDO.convolution-symbols-URd}) that for any $x \in U$ the $x$-symbol $p_{1,m_{1}}*p_{2,m_{2}}(x,.)$ 
% depends only on $p_{m_{1}}(x,.)$ and $p_{m_{2}}(x,.)$. However, the value of  $(p_{1,m_{1}}*p_{2,m_{2}})(x,\xi)$ at $(x,\xi)\in \URd$ depends on all 
% the $\eta$-values $p_{m_{1}}(x,\eta)$ and $p_{m_{2}}(x,\eta)$ as $\eta$ ranges over $\Rd$. Thus we may localize the product of Heisenberg symbols with 
% respect to $x$, but not respect to $(x,\xi)$, that is, the product of \psivdos\ is local, but is not microlocal. 
% \end{remark}

\subsection{The distribution kernels of $\mathbf{\Psi_{H}}$DO's}
An important fact about \psidos\ is their characterization in terms of their distribution kernels.  

 First, we extend the notion of homogeneity of functions to distributions. For $K\in \cS'(\Rd)$ and for $\lambda >0$ we let $K_{\lambda}$ denote 
  the element of $\cS'(\Rd)$ such that
     \begin{equation}
           \acou{K_{\lambda}}{f}=\lambda^{-(d+2)} \acou{K(x)}{f(\lambda^{-1}.x)} \quad \forall f\in\cS(\Rd). 
            \label{eq:PsiHDO.homogeneity-K-m}
      \end{equation}
In the sequel we will also use the notation $K(\lambda.x)$ for denoting $K_{\lambda}(x)$. We then say that $K$ is homogeneous of degree $m$, $m\in\C$, 
when $K_{\lambda}=\lambda^m K$ for any $\lambda>0$.

\begin{definition}\label{def:PsiHDO.regular-distributions}
  $\cS'_{\reg}(\Rd)$  consists of tempered distributions on $\Rd$ which are smooth outside the 
  origin. We equip it with the weakest  topology such that  the inclusions of $\cS'_{\reg}(\Rd)$  into $\cS'(\Rd)$ and $C^{\infty}(\Rdo)$ are continuous. 
\end{definition}

\begin{definition}
  $\cK_{m}(\URd)$, $m\in\C$, consists of distributions $K(x,y)$ in $C^\infty(U)\hotimes\cS'_{\reg}(\Rd)$ such that  for some functions $c_{\alpha}(x) \in 
   C^{\infty}(U)$,  $\brak\alpha=m$, we have  
    \begin{equation}
        K(x,\lambda.y)= \lambda^m K(x,y) + \lambda^m\log\lambda
                \sum_{\brak\alpha=m}c_{\alpha}(x)y^\alpha \qquad \text{for any $\lambda>0$}.
    \end{equation}
\end{definition}
 
The interest of considering the distribution class $\cK_{m}(\URd)$ stems from:

\begin{lemma}[{\cite[Prop.~15.24]{BG:CHM}},~{\cite[Lem.~I.4]{CM:LIFNCG}}]\label{prop:PsiHDO.Sm-Km}
   1) Any  $p\in S_{m}(\URd)$ agrees on $\URdo$ with  a distribution $\tau(x,\xi)\in C^{\infty}(U)\hotimes \cS'(\Rd)$ 
    such that $\check \tau_{\xiy}$ is in $\cK_{\hat m}(\URd)$, $\hat m=-(m+d+2)$. 
    
    2) If $K(x,y)$ is  in $\cK_{\hat m}(\URd)$ then the restriction of 
    $\hat{K}_{\yxi }(x,\xi)$ to $\URdo$ belongs to $S_{m}(\URd)$. 
\end{lemma}

This result is a consequence of the solution to the problem of extending a homogeneous function $p\in C^{\infty}(\Rd\setminus 0)$ into a 
homogeneous distribution on $\Rd$ and of the fact that for $\tau \in \cS'(\Rd)$ we have 
\begin{equation}
    (\hat{\tau})_{\lambda} = |\lambda|^{-(d+2)}(\tau_{\lambda^{-1}})^{\wedge}\qquad \forall \lambda \in \R\setminus 0. 
     \label{eq:PsiHDO.dilations-Fourier-transform}
\end{equation}
In particular, if $\tau$ is homogeneous of degree $m$ then $\hat{\tau}$ is homogeneous of degree $-(m+d+2)$.

The relevant class of kernels for the Heisenberg calculus is the following.
\begin{definition}\label{def:PsiHDO.K*}
$\cK^{m}(\URd)$, $m\in \C$, consists of distributions $K\in \cD'(\URd)$ with an asymptotic expansion 
     $K\sim \sum_{j\geq0}K_{m+j}$,  $K_{l}\in \cK_{l}(\URd)$, 
 in the sense that,  for any integer $N$, as soon as $J$ is large enough we have 
   \begin{equation}
K-\sum_{j\leq J}K_{m+j}\in  C^{N}(\URd). 
     \label{eq:PsiHDO.asymptotics-kernel}
 \end{equation} 
\end{definition}

Since under the Fourier transform the asymptotic expansion~(\ref{eq:PsiVDO.asymptotic-expansion-symbols}) 
for symbols corresponds to that for distributions in~(\ref{eq:PsiHDO.asymptotics-kernel}),  
using Lemma~\ref{prop:PsiHDO.Sm-Km} we get: 

\begin{lemma}[{\cite[pp.~133--134]{BG:CHM}}]\label{lem:PsiHDO.characterization.Km}
    Let $K \in \cD'(\URd)$. Then the following are equivalent:\smallskip
    
    (i) The distribution $K$ belongs to $\cK^{m}(\URd)$;\smallskip
   
   (ii) We can put $K$ into the form
  \begin{equation}
      K(x,y)=\check{p}_{\xiy}(x,y)+R(x,y),
  \end{equation}
 for some $p\in S^{\hat{m}}(\URd)$, $\hat{m}=-(m+d+2)$, and some $R\in 
    C^{\infty}(\URd)$.\smallskip 
    
    Moreover, if (i) and (ii) holds and we expand $K \sim \sum_{j\geq 0} K_{m+j}$, $K_{l}\in \cK_{l}(\URd)$, 
    then we have $p\sim \sum_{j\geq 0}p_{\hat{m}-j}$ where $p_{\hat{m}-j} \in S_{\hat{m}-j}(\URd)$ is the restriction to $\URdo$ of 
    $(K_{m+j})^{\wedge}_{\yxi}$. 
\end{lemma}

Next, for $x \in U$ let $\psi_{x}$ denote the affine change to the privileged coordinates at $x$ and let us write  
$(A_{x}^{t})^{-1}\xi=\sigma(x,\xi)$ with  $A_{x}\in \op{GL}_{d+1}(\R)$. 
Since $\psi_{x}(x)=0$ and $\psi_{x*}X_{j}=\frac{\partial}{\partial 
y_{j}}$ at $y=0$ for $j=0,\ldots,d$,  one checks that $\psi_{x}(y)=A_{x}(y-x)$. 

Let $p\in S^{m}(\URd)$. As $p(x,-iX)=p_{\sigma}(x,D)$ with $p_{\sigma}(x,\xi)=p(x,\sigma(x,\xi))=p(x,(A_{x}^{t})^{-1}\xi)$ 
the distribution kernel 
$k_{p(x,-iX)}(x,y)$ of $p(x,-iX)$ is represented by the oscillating integrals
\begin{equation}
   (2\pi)^{-(d+1)} \int e^{i(x-y).\xi} p(x,(A_{x}^{t})^{-1}\xi)d\xi =   (2\pi)^{-(d+1)}|A_{x}|\int e^{iA_{x}(x-y).\xi} p(x,\xi)d\xi.
\end{equation}
Since $\psi_{x}(y)=A_{x}(y-x)$ we deduce that 
\begin{equation}
    k_{p(x,-iX)}(x,y)=|\psi_{x}'| \check{p}_{\xiy}(x,-\psi_{x}(y)).
    \label{eq:PsiHDO.kernel-quantization-symbol-psiy}
\end{equation}
Combining this with Lemma~\ref{lem:PsiHDO.characterization.Km} then gives: 

\begin{proposition}[{\cite[Thms.~15.39, 15.49]{BG:CHM}}]\label{prop:PsiVDO.characterisation-kernel1}
 Let $P:C_{c}^\infty(U)\rightarrow C^\infty(U)$ be a continuous linear operator with distribution kernel $k_{P}(x,y)$. Then the following are 
 equivalent:\smallskip  
 
 (i) $P$ is a \psivdo\ of order $m$, $m\in \C$.\smallskip 
 
 (ii) There exist $K\in \cK^{\hat{m}}(\URd)$, $\hat{m}=-(m+d+2)$, and $R \in C^{\infty}(U\times U)$ such that 
 \begin{equation}
     k_{P}(x,y)=|\psi_{x}'|K(x,-\psi_{x}(y)) +R(x,y) .
      \label{eq:PsiHDO.characterization-kernel.privilegedx}
 \end{equation}
 Furthermore, if (i) and (ii) hold and $K\sim \sum_{j\geq 0}K_{\hat{m}+j}$, 
$K_{l}\in \cK_{l}(\URd)$, then $P$ has symbol $p\sim \sum_{j\geq 0} p_{m-j}$, $p_{l}\in S_{l}(\URd)$, where 
$p_{m-j}$ is the restriction to $\URdo$ of 
    $(K_{m+j})^{\wedge}_{\yxi}$. 
\end{proposition}

In the sequel we will need a version of Proposition~\ref{prop:PsiVDO.characterisation-kernel1} in Heisenberg coordinates.  
To this end let $\varepsilon_{x}$ denote the coordinate change to the Heisenberg coordinates at $x$ and set $\phi_{x}=\varepsilon_{x}\circ \psi_{x}^{-1}$. 
Recall that $\phi_{x}$ is  a Lie group isomorphism from $G^{(x)}$ to 
$G_{x}U$ such that $\phi_{x}(\lambda.y)=\lambda.\phi_{x}(y)$ for any $\lambda \in \R$. Moreover, using~(\ref{eq:Bundle.Extrinsic.Phiu}) 
one can check that $|\phi_{x}'|=1$ and $\phi_{x}^{-1}(y)=-\phi_{x}(-y)$. Therefore,  
from~(\ref{eq:PsiHDO.kernel-quantization-symbol-psiy}) we see that we can put $ k_{p(x,-iX)}(x,y)$ into the form
\begin{equation}
    k_{p(x,-iX)}(x,y)=|\varepsilon_{x}'| K_{P}(x,-\varepsilon_{x}(y)), \quad K_{P}(x,y)=  
    \check{p}_{\xiy}(x,-\phi_{x}(-y))=\check{p}_{\xiy}(x,\phi_{x}^{-1}(y)).
    \label{eq:PsiHDO.kernel-quantization-symbol}
\end{equation}

In fact, the coordinate changes $\phi_{x}$ give rise to an action on distributions on $\URd$ given by
   \begin{equation}
     K(x,y) \longrightarrow \phi_{x}^{*}K(x,y), \qquad \phi_{x}^{*}K(x,y)=K(x,\phi_{x}^{-1}(y)).
       \label{eq:PsiHD.isomorphism-cK}
   \end{equation}
  Since $\phi_{x}$ depends smoothly on $x$, this action induces a continuous linear isomorphisms of $C^{N}(\URd)$, $N \geq 0$, and $C^{\infty}(\URd)$ 
  onto themselves. As $\phi_{x}(y)$ is 
  polynomial in $y$ in such way that $\phi_{x}(0)=0$  and $\phi_{x}(\lambda.y)=\lambda.\phi_{x}(y)$ for every $\lambda \in \R$, we deduce that the above 
  action also yields a continuous linear isomorphism of $C^{\infty}(U)\hotimes \cS_{\reg}'(\Rd)$ onto itself and, for every $\lambda>0$,  we have 
   \begin{equation}
       (\phi_{x}^{*}K)(x,\lambda.y)= \phi_{x}^{*}[K(x,\lambda.y)], \qquad K \in \cD'(\URd).
        \label{eq:PsiHDO.homogeneity.phix*}
   \end{equation}

   Furthermore, as $\phi_{x}(y)$ is polynomial in $y$ we see that for every $\alpha \in \N^{d+1}$ we can write 
$\phi_{x}(y)^{\alpha}$ in the form $\phi_{x}(y)^{\alpha}=\sum_{\brak\beta=\brak \alpha}d_{\alpha\beta}(x)y^{\beta}$ with $d_{\alpha\beta}\in C^{\infty}(\URd)$.  
   It then follows that, for every $m \in \C$,  the map $K(x,y) \rightarrow \phi_{x}^{*}K(x,y)$ induces a linear isomorphisms of $\cK_{m}(\URd)$ and 
   $\cK^{m}(\URd)$ onto themselves. Combining this with~(\ref{eq:PsiHDO.kernel-quantization-symbol}) and 
   Proposition~\ref{prop:PsiVDO.characterisation-kernel1} then gives:

\begin{proposition}\label{prop:PsiVDO.characterisation-kernel2}
 Let $P:C_{c}^\infty(U)\rightarrow C^\infty(U)$ be a continuous linear operator with distribution kernel $k_{P}(x,y)$. Then the following are 
 equivalent:\smallskip  
 
 (i) $P$ is a \psivdo\ of order $m$, $m\in \C$.\smallskip 
 
 (ii) There exist $K_{P}\in \cK^{\hat{m}}(\URd)$, $\hat{m}=-(m+d+2)$, and $R \in C^{\infty}(U\times U)$ such that 
 \begin{equation}
     k_{P}(x,y)=|\varepsilon_{x}'|K_{P}(x,-\varepsilon_{x}(y)) +R(x,y) .
      \label{eq:PsiHDO.characterization-kernel.Heisenberg}
 \end{equation}
Furthermore, if (i) and (ii) hold and $K_{P}\sim \sum_{j\geq 0}K_{P,\hat{m}+j}$,  
$K_{l}\in \cK_{l}(\URd)$, then $P$ has symbol $p\sim \sum_{j\geq 0} p_{m-j}$, $p_{l}\in S_{l}(\URd)$, where 
$p_{m-j}$ is the restriction to $\URdo$ of 
$[K_{P, \hat{m}+j}(x,\phi_{x}^{-1}(y))]^{\wedge}_{\yxi }$.
\end{proposition}

\begin{remark}\label{rem:PsiHDO.principal-symbol-at-x}
 Let $a\in U$. Then~(\ref{eq:PsiHDO.characterization-kernel.Heisenberg}) 
 shows that the distribution kernel of $\tilde{P}=(\varepsilon_{a})_{*}P$ at $x=0$ is 
 \begin{equation}
     k_{\tilde{P}}(0,y)=|\varepsilon_{a}'|^{-1} k_{P}(\varepsilon_{a}^{-1}(0),\varepsilon_{a}^{-1}(y))=K_{P}(a,-y).
      \label{eq:PsiHDO.characterization-kernel.Heisenberg-origin}
 \end{equation}
Moreover, as we are in Heisenberg coordinates already, we have $\psi_{0}=\varepsilon_{0}=\phi_{0}=\op{id}$. Thus, in the 
 form~(\ref{eq:PsiHDO.characterization-kernel.Heisenberg}) for $\tilde{P}$ we have $K_{\tilde{P}}(0,y)=K_{P}(a,y)$. 
Therefore, if we let 
 $p_{m}(x,\xi)$ denote the principal symbol of $P$ and let  $K_{P,\hat{m}}\in \cK_{\hat{m}}(\URd)$ denote the leading kernel of $K_{P}$,  
 then by Proposition~\ref{prop:PsiVDO.characterisation-kernel2} we have 
 \begin{equation}
     p_{m}(0,\xi)=[K_{P,\hat{m}}]^{\wedge}_{\yxi }(a,\xi).
 \end{equation}
This shows that $[K_{P,\hat{m}}]^{\wedge}_{\yxi }(a,\xi)$ is the principal symbol of $P$ at $x=0$ in Heisenberg coordinates centered at $a$.
\end{remark}

\subsection{$\mathbf{\Psi_{H}}$DO's on a general Heisenberg manifold}
Let $(M^{d+1},H)$ be a Heisenberg manifold. As alluded to before the \psivdos\ on an subset of $\Rd$ are \psidos\ of type 
$(\frac{1}{2},\frac{1}{2})$. However, the latter don't make sense on a general manifold, for their class is not preserved by an arbitrary change of chart. Nevertheless, 
when dealing with \psivdos\ this issue is resolved if we restrict ourselves to changes of Heisenberg charts. Indeed, we have:

 \begin{proposition}\label{prop:PsiHDO.invariance}
     Let $U$ (resp.~$\tilde{U}$) be an open subset of $\Rd$ together with a hyperplane bundle $H\subset TU$ (resp.~$\tilde{H}\subset T\tilde{U}$) and a 
    $H$-frame of $TU$ 
    (resp.~a $\tilde{H}$-frame of $T\tilde{U}$). Let $\phi:(U,H)\rightarrow (\tilde{U},\tilde{H})$ be a Heisenberg diffeormorphism and let $\tilde{P}\in 
    \Psi_{\tilde{H}}^{m}(\tilde{U})$.\smallskip 
    
   1) The operator $P=\phi^{*}\tilde{P}$ is a \psivdo\ of order $m$ on $U$.\smallskip 
    
   2) If the distribution kernel of $\tilde{P}$ is of the form~(\ref{eq:PsiHDO.characterization-kernel.Heisenberg}) 
   with $K_{\tilde{P}}(\tilde{x},\tilde{y})\in \cK^{\hat{m}}(\tilde{U}\times \Rd)$ then the distribution kernel of 
   $P$ can be written in the form~(\ref{eq:PsiHDO.characterization-kernel.Heisenberg}) with $K_{P}(x,y)\in \cK^{\hat{m}}(\URd)$ such that 
   \begin{equation}
       K_{P}(x,y) \sim \sum_{\brak\beta\geq \frac{3}{2}\brak\alpha} \frac{1}{\alpha!\beta!} 
       a_{\alpha\beta}(x)y^{\beta}(\partial_{\tilde{y}}^{\alpha}K_{\tilde{P}})(\phi(x),\phi_{H}'(x)y),
        \label{eq:PsiHDO.asymptotic-expansion-KP}
   \end{equation}
   where we have let $a_{\alpha\beta}(x)=\partial^{\beta}_{y}[|\partial_{y}(\tilde{\varepsilon}_{\phi(x)}\circ \phi\circ \varepsilon_{x}^{-1})(y)| 
   (\tilde{\varepsilon}_{\phi(x)}\circ \phi\circ \varepsilon_{x}^{-1}(y)-\phi_{H}'(x)y)^{\alpha}]_{|_{y=0}}$ and $\tilde{\varepsilon}_{\tilde{x}}$ 
   denote the change to the Heisenberg coordinates at $\tilde{x}\in \tilde{U}$. In particular, we have 
   \begin{equation}
       K_{P}(x,y)=|\phi_{H}'(x)| K_{\tilde{P}}(\phi(x),\phi_{H}'(x)y) \qquad \bmod \cK^{\hat{m}+1}(\URd).
       \label{eq:PsiHDO.asymptotic-expansion-KP-principal}
   \end{equation}
\end{proposition}

\begin{remark}
    The version of the above statement in~\cite{BG:CHM} does not contain the asymptotics~(\ref{eq:PsiHDO.asymptotic-expansion-KP-principal}),  
    which will be crucial for giving a global definition of the 
    principal symbol of a \psivdo\ in the next section. For this reason give a detailed proof of the above version in Section~\ref{sec.Invariance}. 
    This proof will also be useful in~\cite{Po:CPDE1} and~\cite{Po:CPDE2} for generalizing Proposition~\ref{prop:PsiHDO.invariance} 
    to holomorphic families of \psivdos\ and to \psivdos\ with parameter. 
\end{remark}

As a consequence of Proposition~\ref{prop:PsiHDO.invariance} we can define \psivdos\ on $M$  acting on the sections of 
a vector bundle $\cE$ over $M$.
\begin{definition}
    $\pvdo^{m}(M,\cE)$,  $m\in \C$,  consists of continuous operators $P:C^{\infty}_{c}(M,\cE) \rightarrow 
    C^{\infty}(M,\cE)$ such that:\smallskip 
    
    (i) The distribution kernel of $P$ is smooth off the diagonal;\smallskip 
    
    (ii) For any trivialization $\tau:\cE_{|_{U}}\rightarrow U\times \C^{r}$ over a 
     local Heisenberg chart $\kappa:U \rightarrow V\subset \Rd$ the operator $\kappa_{*}\tau_{*}(P_{|_{U}})$ belongs to
     $\pvdo^{m}(V, \C^{r}):=\pvdo^{m}(V)\otimes \End \C^{r}$. 
\end{definition}

All the previous properties of \psidos\  on an open subset of $\Rd$ hold \emph{mutatis standis} for \psivdos\ on $M$ acting on sections of $\cE$. 

\subsection{Transposes and adjoints of $\mathbf{\Psi_{H}}$DO's}
Let us now look at the transpose and adjoints of \psivdos. First, given a Heisenberg chart $U \subset \Rd$ we have:

\begin{proposition}\label{prop:PsiHDO.transpose-chart}
    Let $P\in \pvdo^{m}(U)$. Then:\smallskip 
    
    1) The transpose operator $P^{t}$ is a \psivdo\ of order $m$ on $U$.\smallskip 
    
    2) If we write the distribution kernel of 
    $P$ in the form~(\ref{eq:PsiHDO.characterization-kernel.Heisenberg}) 
   with $K_{P}\in \cK^{\hat{m}}(\URd)$ then $P^{t}$ can be written in the form~(\ref{eq:PsiHDO.characterization-kernel.Heisenberg}) with 
   $K_{P^{t}}\in  \cK^{\hat{m}}(\URd)$ such that 
   \begin{equation}
       K_{P^{t}}(x,y) \sim  \sum_{\frac{3}{2}\brak\alpha \leq \brak \beta} \sum_{|\gamma|\leq |\delta| \leq 2|\gamma|} 
       a_{\alpha\beta\gamma\delta}(x) y^{\beta+\delta}  
       (\partial^{\gamma}_{x}\partial_{y}^{\alpha}K_{P})(x,-y), 
        \label{eq:PsiHDO.transpose-expansion-kernel}
   \end{equation}
   where 
   $a_{\alpha\beta\gamma\delta}(x)=\frac{|\varepsilon_{x}^{-1}|}{\alpha!\beta!\gamma!\delta!}
   [\partial_{y}^{\beta}(|\varepsilon_{\varepsilon_{x}^{-1}(-y)}'|(y-\varepsilon_{\varepsilon_{x}^{-1}(y)}(x))^{\alpha})
   \partial_{y}^{\delta}(\varepsilon_{x}^{-1}(-y)-x)^{\gamma}](x,0)$. In particular, 
   \begin{equation}
       K_{P^{t}}(x,y)=K_{P}(x,-y) \quad \bmod \cK^{\hat m+1}(\URd). 
        \label{eq:PsiHDO.transpose-principal-kernels}
   \end{equation}
\end{proposition}

 \begin{remark}
     The asymptotic expansion~(\ref{eq:PsiHDO.transpose-expansion-kernel}) is not stated in~\cite{BG:CHM}, 
     but we need it in order to determine the global principal symbol of the transpose of a 
     \psivdo\ (see next section). A detailed proof of Proposition~\ref{prop:PsiHDO.transpose-chart}
     can be found in Section~\ref{sec:transpose}. %for sake of completeness and of further generalizations in~\cite{Po:CPDE1} and \cite{Po:CPDE2}. 
 \end{remark}

Using this result, or its version in~\cite{BG:CHM}, we obtain: 

\begin{proposition}[{\cite[Thm.~17.4]{BG:CHM}}]\label{prop:PsiHDO.transpose-adjoint}
  Let $P:C^{\infty}(M,\cE)\rightarrow C^{\infty}(M,\cE)$ be a \psivdo\ of order $m$. Then:
 
  1) The transpose operator $P^{t}:\cE'(M,\cE^{*})\rightarrow \cD'(M,\cE^{*})$ is a \psivdo\ of order $m$;\smallskip 
  
 2) If $M$ is endowed with a smooth positive density and $\cE$ with a Hermitian metric then the adjoint operator $P^{*}:C^{\infty}(M,\cE)\rightarrow 
 C^{\infty}(M,\cE)$ is a \psivdo\ of order $\overline{m}$. 
\end{proposition}

\section{Principal symbol and model operators.}\label{sec:principal-symbol}
In this section we define the principal symbols and model operators of \psivdos\ and check their main properties.

\subsection{Principal symbol and model operators}
% Let us now give an intrinsic definition of the principal symbol of  a \psivdo\ and of its model 
% operator at a point. 
% 
Let $\fg^{*}M$ be the dual bundle of $\fg M$ with canonical projection $\pi:\fg^{*}M \rightarrow M$.
%be the canonical projection of $\fg^{*}M$ onto $M$. 

\begin{definition}
    $S_{m}(\fg^{*}M,\cE)$, $m\in \C$, is the space of sections $p\in C^{\infty}(\fg^{*}M\setminus 0, \End \pi^{*}\cE)$ which are homogeneous of 
    degree $m$ in the sense that, for any $\lambda>0$, we have 
    \begin{equation}
        p(x,\lambda.\xi)=\lambda^{m}p(x,\xi) \qquad \forall (x,\xi) \in \fg^{*}M\setminus 0,
    \end{equation}
    where $\xi \rightarrow \lambda.\xi$ denotes the dilation~(\ref{eq:PsiHDO.Heisenberg-dilation-fg*}). 
\end{definition}

Let $P\in \pvdo^{m}(M)$ and for $j=1,2$ let $\kappa_{j}$ be a Heisenberg chart with domain 
$V_{j}\subset M$ and let $\phi:U_{1}\rightarrow U_{2}$ be the corresponding transition map, where we have let $U_{j}=\kappa_{j}(V_{1}\cap V_{2})\subset \Rd$.

Let us first assume that $\cE$ is the trivial line bundle, so that $P$ is a scalar operator. For $j=1,2$ we let 
$P_{j}:=\kappa_{j*}(P_{|_{V_{1}\cap V_{2}}})$, so that $P_{1}=\phi^{*} P_{2}$. Since $P_{j}$ belongs to $\pvdo^{m}(U_{j})$ its distribution kernel is of 
the form~(\ref{eq:PsiHDO.characterization-kernel.Heisenberg}) with $K_{P_{j}}\in \cK^{\hat m}(U_{j}\times \Rd)$. 
Moreover, by Proposition~\ref{prop:PsiHDO.invariance} we have 
\begin{equation}
    K_{P_{1}}(x,y)=|\phi_{H}'(x)| K_{P_{2}}(\phi(x),\phi_{H}'(x)y) \qquad \bmod \cK^{\hat{m}+1}(U_{1}\times\Rd).
\end{equation}
Therefore, if we let $K_{P_{j},\hat{m}}\in \cK_{\hat m}(U_{j}\times\Rd)$ be the leading kernel of $K_{P_{j}}$ then we get
\begin{equation}
    K_{P_{1},\hat m}(x,y)=|\phi_{H}'(x)|K_{P_{2},\hat m}(\phi(x),\phi_{H}'(x)y).
     \label{eq:PsiHDO.principal-kernel-transition}
\end{equation}

Next, for $j=1,2$ we define 
\begin{equation}
    p_{j,m}(x,\xi)=[K_{P_{j},\hat{m}}]^{\wedge}_{\yxi }(x,\xi), \qquad (x,\xi)\in U_{j}\times \Rdo. 
     \label{eq:PsiHDO.global-principal-symbol}
\end{equation}
By  Remark~\ref{rem:PsiHDO.principal-symbol-at-x} for any $a \in U_{j}$ the symbol $p_{j}(a,.)$ 
yields in Heisenberg coordinates centered at $a$ the principal symbol of $P_{j}$ at $x=0$.  Moreover, 
since $\phi_{H}'(a)$ is a linear map, from~(\ref{eq:PsiHDO.principal-kernel-transition}) we get
\begin{equation}
    p_{1,m}(x,\xi)= p_{2,m}(\phi(x),[\phi_{H}'(x)]^{-1t}\xi).
\end{equation}
This shows that $p_{m}: =\kappa_{1}^{*}p_{1,m}$ is an element of $S_{m}(\fg^{*}(V_{1}\cap V_{1}))$ which is independent of the choice of the chart 
$\kappa_{1}$. Since $S_{m}(\fg^{*}M)$ is a sheaf this gives rise to a uniquely defined symbol $p_{m}\in S_{m}(\fg^{*}M)$.

When $\cE$ is a general vector bundle, the above construction can be carried out similarly,  so that we obtain: 

\begin{theorem}\label{prop:PsiHDO.principal-symbol}
    For any $P \in \pvdo^{m}(M,\cE)$ there is a unique symbol $\sigma_{m}(P)\in S_{m}(\fg^{*}M, \cE)$ such that if in a local trivializing Heisenberg 
    chart $U\subset \Rd$ we let $K_{P,\hat{m}}(x,y)\in \cK_{\hat{m}}(\URd)$ be the leading kernel for the kernel $K_{P}(x,y)$  in the 
    form~(\ref{eq:PsiHDO.characterization-kernel.Heisenberg})  
    for $P$,  then we have
    \begin{equation}
       \sigma_{m}(P)(x,\xi)=[K_{P,\hat{m}}]^{\wedge}_{\yxi }(x,\xi), \qquad (x,\xi)\in U\times \Rdo. 
    \end{equation}
    
    Equivalently, for any $x_{0} \in M$ the symbol $\sigma_{m}(P)(x_{0},.)$ agrees in trivializing Heisenberg coordinates centered at $x_{0}$ with 
  the principal symbol of $P$ at $x=0$. 
\end{theorem}

\begin{definition}\label{def:PsiHDO.principal-symbol}
   For $P\in \pvdo^{m}(M,\cE)$ the symbol $\sigma_{m}(P)\in S_{m}(\fg^{*}M,\cE)$ provided by Theorem~\ref{prop:PsiHDO.principal-symbol}
     is called the principal symbol of $P$. 
\end{definition}

\begin{remark}
    Since we have two notions of principal symbol we shall distinguish between them by saying that  $\sigma_{m}(P)$ is the global 
    principal symbol of $P$ and that in a local trivializing chart the principal symbol $p_{m}$  of $P$ in the sense of~(\ref{eq:PsiVDO.asymptotic-expansion-symbols}) 
    is the local principal symbol of $P$ in this chart. 
   \end{remark}
   
    In a local Heisenberg chart $U\subset \Rd$ the global symbol $\sigma_{m}(P)$ and the local principal symbol $p_{m}$ of $P\in \pvdo^{m}(U)$  
   can be easily related to each other. Indeed, by Proposition~\ref{prop:PsiVDO.characterisation-kernel2} we have 
    \begin{equation}
       p_{m}(x,\xi)= [K_{P,\hat{m}}(x,\phi_{x}^{-1}(y))]^{\wedge}_{\yxi }(x,\xi),
    \end{equation}
    where $K_{P,\hat{m}}$ denotes the leading kernel for the kernel $K_{P}$ in the form~(\ref{eq:PsiHDO.characterization-kernel.Heisenberg}) for $P$. 
    By combining this with the definition~(\ref{eq:PsiHDO.global-principal-symbol}) 
    of $\sigma_{m}(P)$ we thus get 
    \begin{gather}
         p_{m}(x,\xi)= (\hat{\phi}_{x}^{*}\sigma_{m}(P))(x,\xi), 
           \label{eq:PsiHDO.local-global-principal-symbol} \\
       (\hat{\phi}_{x}^{*}\sigma_{m}(P))= 
       [[\sigma_{m}(P)]^{\vee}_{\xiy}(x,\phi_{x}^{-1}(y))]^{\wedge}_{\yxi }=  [\phi_{x}^{*}[\sigma_{m}(P)]^{\vee}_{\xiy}]^{\wedge}_{\yxi },
       \label{eq:PsiHDO.isomorphism-symbols}
     \end{gather}
    where $\phi_{x}^{*}$ is the isomorphism map~(\ref{eq:PsiHD.isomorphism-cK}).
    In particular, since the latter is a linear isomorphism of $\cK_{m}(\URd)$ onto itself, 
    we see that the map  $p\rightarrow \hat{\phi}_{x}^{*}p$ is a linear isomorphism of $S_{m}(\URd)$ onto itself. 

\begin{example}
    Let $X_{0},\ldots,X_{d}$ be a local $H$-frame of $TM$ near a point $a\in M$. In any Heisenberg chart associated with this frame the Heisenberg 
    symbol of $X_{j}$ is $\frac{1}{i}\xi_{j}$. In particular, this is true in Heisenberg coordinates centered at $a$. Thus the (global) principal symbol 
    of $X_{j}$ is equal to $\frac{1}{i}\xi_{j}$ in the local trivialization of $\fg^{*}M\setminus 0$ defined by the frame $X_{0},\ldots,X_{d}$. 
    More generally, for any differential $P=\sum_{\brak \alpha\leq m} a_{\alpha}(x) X^{\alpha}$ on $M$ we have
    \begin{equation}
        \sigma_{m}(P)(x,\xi)=\sum_{\brak \alpha\leq m} a_{\alpha}(x) i^{-|\alpha|}\xi^{\alpha}.
         \label{eq:Pincipal.example}
    \end{equation}
    Thus, for differential operators the global and local principal symbols agree in suitable coordinates. Alternatively, this result follows from the fact 
    that the isomorphism~(\ref{eq:PsiHD.isomorphism-cK}) induces the identity map on distributions $K(x,y)$ supported in $U\times \{y=0\}$.
\end{example}    
    
\begin{proposition}\label{prop:PsiHDO.surjectivity-principal-symbol-map}
    For every $m\in \C$ the principal symbol map $\sigma_{m}: \pvdo^{m}(M,\cE) \rightarrow S_{m}(\fg^{*}M,\cE)$ gives rise to a linear isomorphism
   $\pvdo^{m}(M,\cE)/\pvdo^{m-1}(M,\cE)\stackrel{\sim}{\longrightarrow}  S_{m}(\fg^{*}M,\cE)$. 
\end{proposition}
\begin{proof}
  By construction the principal symbol of $P\in \pvdo^{m}(M,\cE)$ depends only on his principal part in local coordinates and vanishes everywhere 
  if, and only if, the order of $P$ is~$\leq m-1$. Therefore, the kernel of the principal symbol map 
  $\sigma_{m}$ is $\pvdo^{m-1}(M,\cE)$, so $\sigma_{m}$ induces an injective linear map  
  $\pvdo^{m}(M,\cE)/\pvdo^{m-1}(M,\cE)\rightarrow  S_{m}(\fg^{*}M,\cE)$. 
  
  To complete the proof it is enough to show that $\sigma_{m}$ is surjective. To this end consider a symbol 
  $p_{m}(x,\xi)\in S_{m}(\fg^{*}M,\cE)$ and let $(\varphi_{i})_{i\in I}$ be 
  a partition of the unity subordinated to an open covering $(U_{i})_{i\in I}$ of $M$ by
 domains of Heisenberg charts $\kappa_{i}:U_{i}\rightarrow V_{i}$ over which there are trivializations $\tau_{i}:\cE_{|_{U_{i}}}\rightarrow 
  U_{i}\times \C^{r}$. For each index $i$ let $\psi_{i}\in C^{\infty}(U_{i})$ be such that $\psi_{i}=1$ near 
  $\supp \varphi_{i}$ and set
\begin{equation}
     p_{m}^{(i)}(x,\xi)=(1-\chi(\xi))(\hat{\phi}_{i,x}^{*}\kappa_{i*}\tau_{i*}p_{m|_{\fg^{*}U_{i}\setminus 0}})(x,\xi) \in S_{m}(V_{i})\otimes \End 
  \C^{r},
\end{equation}
 where $\chi\in C^{\infty}(\Rd)$ is such that $\chi=1$ near the origin and $\hat{\phi}_{i,x}^{*}$ denotes the 
 isomorphism~(\ref{eq:PsiHDO.isomorphism-symbols})  with respect to the 
 chart $V_{i}$. Then we define a   a \psivdo\ of order $m$ by letting
  \begin{equation}
      P=\sum \varphi_{i}[\tau_{i}^{*}\kappa_{i}^{*} p_{m}^{(i)}(x,-iX)]\psi_{i}.
  \end{equation}
  
    For for every index $i$  the local principal symbol 
    of $\varphi_{i}[\tau_{i}^{*}\kappa_{i}^{*} p_{m}^{(i)}(x,-iX)]\psi_{i}$ in the chart $V_{i}$  
    is $\varphi_{i}\circ\kappa_{i}^{-1}(\hat{\phi}_{i,x}^{*}\kappa_{i*}\tau_{i*}p_{m|_{\fg^{*}U_{i}\setminus 0}})$, so 
    by~(\ref{eq:PsiHDO.local-global-principal-symbol})  its global principal is 
    $\varphi_{i}\circ\kappa_{i}^{-1}(\kappa_{i*}\tau_{i*}p_{m|_{\fg^{*}U_{i}\setminus 0}})$, which pulls back to $\varphi_{i}p_{m}$ 
    on $U_{i}$. It follows that the global principal symbol of $P$ is $\sigma_{m}(P)= \sum_{i} \varphi_{i} p_{m}=p_{m}$. 
    This proves the surjectivity of  the map $\sigma_{m}$, so the proof is now complete.
\end{proof}
  
  Next, granted the above definition of the principal symbol, we can define the model operator at a point as follows. 
  
  \begin{definition}\label{def:PsiHDO.model-operator}
    Let $P\in \pvdo^{m}(M,\cE)$ have (global) principal symbol $\sigma_{m}(P)$.  Then the model operator of $P$ at $a\in M$ is the left-invariant 
    \psivdo-operator $P^{a}:\cS_{0}(G_{a}M,\cE_{a})\rightarrow 
    \cS_{0}(G_{a}M,\cE_{x})$ with symbol $\sigma_{m}(P)^{\vee}_{\xiy}(a,.)$, i.e.,
    \begin{equation}
        P^{a}f(x)=\acou{\sigma_{m}(P)^{\vee}_{\xiy}(a,y)}{f(x.y^{-1})}, \qquad f \in 
        \cS_{0}(G_{a}M,\cE_{a}).   
    \end{equation}
\end{definition}

Consider a local trivializing 
chart $U \subset \Rd$ near $a$ and let us relate the model operator $P^{a}$ on $G_{a}M$ to the operator $P^{(a)}=\tilde{p}_{m}^{a}(-iX^{(a)})$ on $G^{(a)}$ 
defined 
using the local principal symbol $\tilde{p}_{m}(x,\xi)$ of $P$ in this chart. Using~(\ref{eq:PsiHDO.convolution-operator}) 
and~(\ref{eq:PsiHDO.local-global-principal-symbol}) for $f\in \cS_{0}(\Rd)$ we get 
\begin{equation}
    P^{(a)}f(y)=\acou{(p_{m}^{a})^{\vee}(z)}{f(y.z^{-1})}=\acou{(\sigma_{m}(P)^{\vee}_{\xiy}(x,\phi_{a}^{-1}(y))}{f(y.z^{-1})}.
\end{equation}
Since $|\phi_{a}'|=1$ and $\phi_{a}$ is a Lie group isomorphism from $G^{(a)}$ onto $G_{a}M$ we obtain
\begin{equation}
      P^{(a)}f(y)= \acou{(\sigma_{m}(P)^{\vee}_{\xiy}(x,y)}{f\circ \phi_{a}^{-1}(y.\phi_{a}(x)^{-1})}=(\phi_{a}^{*}P^{a})f(y).
\end{equation}
In particular, we have
\begin{equation}
    P^{a}=(\phi_{a})_{*}p_{m}^{a}(-iX^{(a)}).
    \label{eq:PsiHDO.global-local-model-operator}
\end{equation}

\subsection{Composition of principal symbols and model operators}
Let us now look at the composition of principal symbols. To this end for $a \in M$ we 
let $*^{a}:S_{m_{1}}(\Rd)\times S_{m_{2}}(\Rd)\rightarrow S_{m_{1}+m_{2}}(\Rd)$ be the 
convolution product for symbols defined by the product law of $G_{a}M$ under the identification $G_{a}M\simeq \Rd$ provided by a $H$-frame 
$X_{0},\ldots,X_{d}$ of $TM$ near $a$, that is, 
\begin{equation}
    (p_{m_{1}}*^{a}p_{m_{j}})(-iX^{a})= p_{m_{1}}(-iX^{a})\circ p_{m_{2}}(-iX^{a}), \qquad p_{m_{j}}\in S_{m_{j}}(\Rd). 
     \label{eq:PsiHDO.global-convolution-symbols}
\end{equation}

Let $U \subset \Rd$ be a local trivializing Heisenberg chart
chart  near $a$ and for $j=1,2$ let $P_{j}\in \pvdo^{m_{j}}(U)$ have (global) principal symbol $\sigma_{m_{j}}(P_{j})$. 
Recall that under the trivialization of $GU$ provided 
by the $H$-frame $X_{0},\ldots,X_{d}$ we have $P^{a}_{j}=\sigma(P_{j})(x, -iX^{a})$. Thus, 
\begin{equation}
    [\sigma_{m_{j}}(P_{j})(x,,)*^{a}\sigma_{m_{j}}(P_{j})(x,.)](-iX^{a})=P_{1}^{a}P_{2}^{a}.
     \label{eq:PsiHDO.convolution-symbol-product-model-operators}
\end{equation}

On the other hand, using~(\ref{eq:PsiHDO.local-global-principal-symbol}) and~(\ref{eq:PsiHDO.global-local-model-operator}) 
we see that $ \hat{\phi}_{a}^{*}[p_{m_{1}}*^{a}p_{m_{2}}](-iX^{a})$ is equal to
\begin{multline}
   \phi_{a}^{*}[p_{m_{1}}(-iX^{a})\circ p_{m_{2}}(-iX^{a})] = 
    \phi_{a}^{*}[p_{m_{1}}(-iX^{a})]\circ  \phi_{a}^{*}[p_{m_{2}}(-iX^{a})] \\ = 
    (\hat{\phi}_{a}^{*}p_{m_{1}})(-iX^{(a)})\circ (\hat{\phi}_{a}^{*}p_{m_{2}})(-iX^{(a)}) =  
     [(\hat{\phi}_{a}^{*}p_{m_{1}})*^{(a)}(\hat{\phi}_{a}^{*}p_{m_{2}})](-iX^{(a)}).
\end{multline}
Hence we have
\begin{equation}
    p_{m_{1}}*^{a}p_{m_{2}}=(\hat{\phi}_{a})_{*}[(\hat{\phi}_{a}^{*}p_{m_{1}})*^{(a)}(\hat{\phi}_{a}^{*}p_{m_{2}})] \qquad \forall p_{m_{j}}\in 
    S_{m_{j}}(\Rd),
    \label{eq:PsiHDO.global-local-convolution-symbols}
\end{equation}
where $(\hat{\phi}_{a})_{*}$ denotes the inverse of $\hat{\phi}_{a}^{*}$. Since $\hat{\phi}_{a}^{*}$, its inverse and $*^{(a)}$ depend smoothly on 
$a$, we deduce that that so does $*^{a}$.  Therefore, we obtain: 

\begin{proposition}
  The group laws on the fibers of $GM$ give rise to a convolution product,
    \begin{gather}
        *:S_{m_{1}}(\fg^{*}M,\cE)\times S_{m_{2}}(\fg^{*}M,\cE) \longrightarrow S_{m_{1}+m_{2}}(\fg^{*}M,\cE),\\
%     \end{gather}
% such that for symbols $p_{m_{j}}\in S_{m_{j}}(\fg^{*}M,\cE)$, $j=1,2$, we have
%     \begin{gather}
        p_{m_{1}}*p_{m_{2}}(x,\xi)=[p_{m_{1}}(x,.)*^{x}p_{m_{2}}(x,.)](\xi), \qquad p_{m_{j}}\in S_{m_{j}}(\fg^{*}M,\cE),
    \end{gather}
 where $*^{x}$ denote the convolution product for symbols on $G_{x}M$.
\end{proposition}

Notice that~(\ref{eq:PsiHDO.global-local-convolution-symbols}) shows that, under the 
relation~(\ref{eq:PsiHDO.local-global-principal-symbol}) between local and global principal symbols, the convolution 
product~(\ref{eq:PsiHDO.global-convolution-symbols})  for 
global principal symbols corresponds to the convolution product~(\ref{eq:PsiHDO.convolution-symbols-URd}) 
for local principal symbols. Since by Proposition~\ref{prop:PsiHDO.composition} the latter yields 
the \emph{local} principal symbol of the product of two \psivdos\ in a local chart, we deduce that the convolution 
product~(\ref{eq:PsiHDO.global-convolution-symbols}) yields the 
\emph{global} principal  symbol of the product two \psivdos. 

Moreover,  by~(\ref{eq:PsiHDO.convolution-symbol-product-model-operators}) 
the global convolution product~(\ref{eq:PsiHDO.global-convolution-symbols}) corresponds to the product of model 
operators, so the model operator of a product of two \psivdos\ is equal to the product of the model operators. We have thus proved:

\begin{proposition}\label{prop:PsiHDO.composition2}
    For $j=1,2$ let $P_{j}\in \pvdo^{m_{j}}(M,\cE)$ and suppose that $P_{1}$ or $P_{2}$ 
    is properly supported.\smallskip 
   
    1) We have $\sigma_{m_{1}+m_{2}}(P_{1}P_{2})=\sigma_{m_{1}}(P)*\sigma_{m_{2}}(P)$.\smallskip
    
    2) At every $a\in M$ the model operator of $P_{1}P_{2}$ is $(P_{1}P_{2})^{a}=P^{a}_{1}P_{2}^{a}$.
\end{proposition}

\subsection{Principal symbol of transposes and adjoints}
In this subsection we shall determine the principal symbols and the model operators of transposes and adjoints of \psivdos. 

Recall that if $P\in \pvdo^{m}(M,\cE)$ then by Proposition~\ref{prop:PsiHDO.transpose-adjoint} 
its transpose $P^{t}:C^{\infty}_{c}(M,\cE^{*})\rightarrow C^{\infty}(M,\cE)$ is a \psivdo\ of order $m$ and its adjoint 
$P^{*}:C^{\infty}_{c}(M,\cE)\rightarrow C^{\infty}(M,\cE)$ is a \psivdo\ of order $\overline{m}$ (assuming $M$ endowed with a positive density and 
$\cE$ with a Hermitian metric in order to define the adjoint). 

\begin{proposition}\label{prop:PsiHDO.transpose-global}
  Let $P \in \pvdo^{m}(M,\cE)$ have principal symbol $\sigma_{m}(P)$. Then:\smallskip 
  
 1) The principal symbol of $P^{t}$ is $\sigma_{m}(P^{t})(x,\xi)= \sigma_{m}(x,-\xi)^{t}$ (this is an element of $S_{m}(\fg^{*}M,\cE*)$);\smallskip 
  
  2) If $P^{a}$ is the model operator of $P$ at $a$ then the model operator of $P^{t}$ at $a$ is the transpose operator 
  $(P^{a})^{t}: \cS_{0}(G_{x}M,\cE_{x}^{*})\rightarrow \cS_{0}(G_{x}M,\cE_{x}^{*})$.
\end{proposition}
\begin{proof}
   Let us first assume that $\cE$ is the trivial line bundle and that $P$ is a scalar operator. In a local Heisenberg chart $U \subset \Rd$ 
   we can write the distribution kernels of $P$ and $P^{t}$ in the form~(\ref{eq:PsiHDO.characterization-kernel.Heisenberg}) with 
   $K_{P}$ and $K_{P^{t}}$ in  $\cK^{\hat{m}}(\URd)$. Let $K_{P,\hat{m}}$  and $K^{t}_{P^{t},\hat{m}}$ denote the principal parts of 
   $K_{P}$ and $K_{P^{t}}$ respectively. Then the principal symbols of $P$ and $P^{t}$ are
  $\sigma_{m}(P)(x,\xi)=(K_{P,\hat{m}})^{\wedge}_{\yxi}(x,\xi)$ and $\sigma_{m}(P^{t})(x,\xi)=(K_{P^{t},\hat{m}})^{\wedge}_{\yxi}(x,\xi)$ respectively.
   Since~(\ref{eq:PsiHDO.transpose-principal-kernels}) implies  that $K_{P^{t},\hat{m}}(x,y)=K_{P,\hat{m}}(x,-y)$ and the Fourier transform commutes 
   with the multiplication by $-1$ we get
   \begin{equation}
       \sigma_{m}(P^{t})(x,\xi)=\sigma_{m}(P)(x,-\xi). 
   \end{equation}
   
  Next, for $a \in U$ let $p\in S_{m}(G_{a})$ and let $P$ be the left-invariant \psivdo\ with symbol $p$. Then the transpose $P^{t}$ is 
   such that, for $f$ and $g$ in $S_{0}(G_{a}U)$, we have 
   \begin{multline}
       \acou{P^{t}f}{g}=\acou{f}{Pv}= \acou{1}{f(x)(Pg)(x)}=\acou{1}{f(x)\acou{\check{p}(y)}{g(x.y^{-1})}}\\
       = \acou{1\otimes \check{p}(x,y)}{f(x)g(x.y^{-1})}.
   \end{multline}
   Therefore, using the change of variable $(x,y)\rightarrow (x.y^{-1},y^{-1})$ and the fact $y^{-1}=-y$ we get
   \begin{equation}
        \acou{P^{t}f}{g}=  \acou{1\otimes \check{p}(x,-y)}{f(x)g(x.y^{-1})}=\acou{1}{f(x)\acou{\check{p}(-y)}{g(x.y^{-1})}}.
       \label{eq:PsiHDO.transpose-model-operator-transpose}
   \end{equation}
    Since $\check{p}(-y)=\check{p}^{t}(y)$ with $p^{t}(\xi)=p(-\xi)$, we obtain
    \begin{equation}
        \acou{P^{t}f}{g}= \acou{1}{f(x)\acou{\check{p}^{t}(y)}{g(x.y^{-1})}}=\acou{(p^{t}*f)(x)}{g(x)}.
    \end{equation}
    Thus $P^{t}$ is the left-convolution operator with symbol $p^{t}(\xi)=p(-\xi)$. 
    
    Now, since the model operator $(P^{t})^{a}$ is the 
    left-invariant \psivdo\ with symbol $\sigma_{m}(P^{t})(a,\xi)=\sigma_{m}(P)(a,-\xi)$, it follows that $(P^{t})^{a}$ agrees 
    with the transpose of $P^{a}$. 
    
    In the general case, when $\cE$ is not the trivial bundle, we can similarly show that $P^{t}$ is a \psivdo\ of order $m$ with principal symbol 
    $\sigma_{m}(P^{t})(x,\xi)=\sigma_{m}(P)(x,-\xi)^{t}$ and such that at every point $a\in M$ its model operator at $a$ is the transpose 
    $(P^{a})^{t}$ of $P^{a}$.
\end{proof}

Assume now that $M$ is endowed with a positive density and $\cE$ with a Hermitian metric respectively and let $L^{2}(M,\cE)$ be the associated 
$L^{2}$-Hilbert space. 

\begin{proposition}\label{prop:PsiHDO.adjoint-manifold}
    Let $P \in \pvdo^{m}(M,\cE)$ have principal symbol $\sigma_{m}(P)$. Then:\smallskip

    1) The principal symbol of $P^{*}$ is $\sigma_{\bar{m}}(P^{*})(x,\xi)=\sigma_{m}(P)(x,\xi)^{*}$.\smallskip 
  
  2) If $P^{x}$ denotes the model operator of $P$ at $x\in M$ then the model operator of $P^{*}$ at $x$ is  
  the adjoint  $(P^{x})^{*}$ of $P^{x}$. 
 \end{proposition} 
 \begin{proof}
    Let us first assume that $\cE$ is the trivial line bundle, so that $P$ is a scalar operator. Moreover, since the above statements are local ones, it 
    is enough to prove them in a local Heisenberg chart $U \subset \Rd$ and we may assume  that $P$ is a \psivdo\ on $U$. 
    
    Let  $\overline{P}:C^{\infty}_{c}(U)\rightarrow C^{\infty}(U)$ be the conjugate operator of $P$, so that 
    $\overline{P}u=\overline{P(\overline{f})}$ for any  $f\in C^{\infty}_{c}(U)$. 
    By Proposition~\ref{prop:PsiVDO.characterisation-kernel2} 
    the distribution kernel of $P$ of the 
    form~(\ref{eq:PsiHDO.characterization-kernel.Heisenberg}) with $K_{P}(x,y)$ in $\cK_{\hat{m}}(\URd)$, so the kernel of $\overline{P}$ takes the form
    \begin{equation}
        k_{\overline{P}}(x,y)=\overline{k_{P}(x,y)}=|\varepsilon_{x}'|\overline{K_{P}(x,y)} \quad \bmod C^{\infty}(U\times U).
    \end{equation}
    
    Since the conjugation of distribution $K(x,y)\rightarrow \overline{K(x,y)}$ induces an anti-linear isomorphism from $\cK^{\hat{m}}(\URd)$ onto 
    $\cK^{\hat{\bar m}}(\URd)$, it follows from Proposition~\ref{prop:PsiVDO.characterisation-kernel2}  
    that $\overline{P}$ is a \psivdo\ of order $\hat{m}$ and its kernel can be 
    put into the form~(\ref{eq:PsiHDO.characterization-kernel.Heisenberg}) 
    with $K_{\overline{P}}(x,y)=\overline{K_{P}(x,y)}$. In particular, if we let $K_{P,\hat{m}}\in \cK_{\hat{m}}(\URd)$ denote 
    the leading kernel of $K_{P}$ then the leading kernel of $K_{\overline{P}}$ is $ \overline{K_{P,\hat{m}}}$. Thus $\overline{P}$ has principal symbol, 
    \begin{equation}
        \sigma_{\bar{m}}(\overline{P})(x,\xi)=[\overline{K_{P,\hat{m}}}]^{\wedge}_{\xiy}(x,\xi)=\overline{[(K_{P,\hat{m}})^{\wedge}_{\xiy}(x,-\xi)}= 
        \overline{\sigma_{m}(x,-\xi)}.   
    \end{equation}
    
    Moreover, since $\sigma_{\bar m}(\overline{P})^{\vee}_{\xiy}(x,y)= \overline{\sigma_{m}(P)^{\vee}_{\xiy}(x,y)}$ the model operator at $a\in U$ of 
    $\overline{P}$ is  such that, for any $f\in \cS_{0}(G_{a}U)$, we have 
    \begin{equation}
        (\overline{P})^{a}f(x)=\acou{ \overline{\sigma_{m}(P)^{\vee}_{\xiy}(x,y)}}{f(x.y^{-1})}= 
        \overline{\acou{\sigma_{m}(P)^{\vee}_{\xiy}(x,y)}{\overline{f(x.y^{-1}}}}= \overline{P^{a}}f(x),
    \end{equation}
    so $(\overline{P})^{a}$ agrees with $\overline{P^{a}}$. 
    
     Combining this with Proposition~\ref{prop:PsiHDO.transpose-chart} and 
    Proposition~\ref{prop:PsiHDO.transpose-global} we thus see that $\overline{P}^{t}$ is a \psivdo\ of order $\overline{m}$ such that:\smallskip 
    
    - If we put the kernel of $\overline{P}^{t}$ into the form~(\ref{eq:PsiHDO.characterization-kernel.Heisenberg}) with respect to
    $K_{\overline{P}^{t}}(x,y)\in \cK_{\hat{\bar{m}}}(\URd)$, then the leading kernel of $K_{\overline{P}^{t}}$ is 
    $K_{\overline{P}^{t},\hat{\bar{m}}}=\overline{K_{P,\hat{m}}(x,-y)}$;\smallskip 
   
   - The global principal symbol of $\overline{P}^{t}$ is 
   $\sigma_{\bar{m}}(\overline{P}^{t})=\overline{\sigma_{m}(P^{t})(x,\xi)}=\overline{\sigma_{\bar m}(P)(x,\xi)}$;\smallskip
   
   - The model operator at $a \in U$ of $\overline{P}^{t}$ is $(\overline{P}^{t})^{a}=\overline{P^{a}}^{t}=(P^{a})^{*}$.\smallskip 
   
    Now, let $d\rho(x)=\rho(x)dx$ be the smooth positive density on $U$ coming from that of $M$. The formal adjoint $P^{*}: C^{\infty}_{c}(U)\rightarrow 
    C^{\infty}(U)$ of $P$ with respect to $d\rho$ is such that%, for functions $f$ and $g$ in  $C^{\infty}_{c}(U)$, we have
    \begin{equation}
        \int_{f} \overline{Pf(x)}g(x)\rho(x)dx = \int_{U}\overline{f(x)}P^{*}g(x)\rho(x)dx, \qquad f,g\in C^{\infty}_{c}(U).
    \end{equation}
    Thus $P^{*}=\rho^{-1}\overline{P}^{t}\rho$, which shows that $P^{*}$ is a \psivdo\ of order $\bar{m}$. Moreover, as in the proof of 
    Proposition~\ref{prop:PsiHDO.transpose-chart} in Section~\ref{sec:transpose}, 
    we can prove that the kernel of $P^{*}$ can be put into the form~(\ref{eq:PsiHDO.characterization-kernel.Heisenberg}) with 
    $K_{P^{*}}(x,y)\in \cK^{\overline{m}}(\URd)$ such that 
       \begin{equation}
       K_{P^{*}}(x,y)  =\rho(x)^{-1}K_{\overline{P}^{t}}(x,y)\rho(\varepsilon_{x}^{-1}(y)) \sim \sum_{\alpha} \frac{1}{\alpha!} \rho(x)^{-1}
       \partial_{y}(\rho(\varepsilon_{x}^{-1}(y))_{|_{y=0}}K_{\overline{P}^{t}}(x,y).
    \end{equation}
In particular, the kernels $K_{P^{*}}(x,y) $ and $K_{\overline{P}^{t}}(x,y)$ agree modulo $\cK^{\hat{\bar{m}}+1}(\URd)$, hence have same leading kernel. It then 
follows that $P^{*}$ and $\overline{P}^{t}$ have same principal symbol and same model operator at a point $a\in U$, that is, we have $\sigma_{\bar 
m}(P^{*})(a,\xi)=\overline{\sigma_{m}(P)(x,\xi)}$ and $(P^{*})^{a}=(P^{a})^{*}$. 

Finally, assume that $\cE$ is a general bundle, so that the restriction of $P$ to $U$ is given by a matrix $P=(P_{ij})$ of \psivdos\ of order $r$. 
Let $h(x)\in C^{\infty}(U,GL_{r}(\C))$, $h(x)^{*}=h(x)$, be the Hermitian metric on $U\times 
\C^{r}$ coming from that of $\cE$. Then the adjoint of $P$ with respect to this Hermitian metric is $P^{*}=\rho^{-1}h^{-1}\overline{P}^{t}h\rho$. Therefore, 
in a the same way as in the scalar case we can prove that $P^{*}$ has principal symbol 
$h(x)^{-1}\overline{\sigma_{m}(P)(x,\xi)}^{t}h(x)=\sigma_{m}(P)(x,\xi)^{*}$ and its model operator at any point $a\in U$ is 
$h(a)^{-1}\overline{P^{a}}^{t}h(a)=(P^{a})^{*}$. 
\end{proof}
 
\section{Hypoellipticity, parametrices and the Rockland condition}
\label{sec:hypoellipticity}
In this section we define a Rockland condition for \psivdos\ and relate it to the invertibility of the principal symbol to get hypoellipticity 
criterions. 

First, by~\cite[Sect.~18]{BG:CHM} in a local Heisenberg chart the invertibility of the local principal symbol of a \psivdo\ 
is equivalent to the existence of a \psivdo-parametrix. Using the 
global principal symbol we can give a global reformulation of this result as follows.

\begin{proposition}\label{thm:PsiHDO.hypoellipticity}
   Let $P:C^{\infty}_{c}(M,\cE) \rightarrow C^{\infty}(M,\cE)$ be a \psivdo\ of order $m$ such that $k:=\Re m>0$. The following are equivalent:\smallskip 

   1) The principal symbol $\sigma_{m}(P)$ of $P$ is invertible with respect to the convolution product for homogeneous 
    symbols;\smallskip 
    
    2) $P$ admits a parametrix $Q$ in $\pvdo^{-m}(M,\cE)$, so that~$PQ=QP=1  \bmod \psinf(M,\cE)$.
  \end{proposition}
 \begin{proof} 
     First, it immediately follows from Proposition~\ref{prop:PsiHDO.composition2} that 2) implies 1). 
     Conversely, in a local trivializing Heisenberg chart~(\ref{eq:PsiHDO.global-local-convolution-symbols})  shows that 
     the invertibility of the global principal $\sigma_{m}(P)$ is equivalent to that of the local 
     principal symbol. Once the latter is granted  Lemma~\ref{lem:PsiHDO.asymptotic-completeness} 
     and Proposition~\ref{prop:PsiHDO.composition} allows us to carry out  in a local trivializing Heisenberg chart  
     the standard parametrix construction in a trivializing Heisenberg chart to get a parametrix for $P$ as a \psivdo\ of order $-m$ 
     (see~\cite[p.~142]{BG:CHM}). A classical partition of the unity 
     argument then allows us to produce a parametrix for $P$ in $\pvdo^{-m}(M,\cE)$.  
\end{proof}
    
When a \psivdo\ has an invertible principal symbol using the Sobolev regularity properties of its parametrices allows us to get:  

\begin{proposition}[{\cite[p.~142]{BG:CHM}}]
      Let $P:C^{\infty}_{c}(M,\cE) \rightarrow C^{\infty}(M,\cE)$ be a \psivdo\ of order $m$ whose principal symbol is invertible. 
      If  $k:=\Re m>0$ then $P$ is hypoelliptic with loss of $\frac{k}{2}$-derivatives, i.e., for any $a \in M$,  any 
   $u \in \cE'(M,\cE)$  and any $s \in \R$, we have 
   \begin{equation}
       \text{$Pu$ is $L^{2}_{s}$ near $a$} \ \Longrightarrow \  \text{$u$ is $L^{2}_{s+k}$ near $a$}.
%         \label{eq:}
   \end{equation}
   In particular, if $M$ is compact then, for any reals $s$ and $s'$, we have the estimate,
   \begin{equation}
       \|f\|_{L^{2}_{s+k}}\leq C_{ss'}(\|Pf\|_{L^{2}_{s}}+\|f\|_{L^{2}_{s'}}), \qquad u \in C^{\infty}(M,\cE).
%         \label{eq:}
   \end{equation}
\end{proposition}

\begin{remark}
    We can give sharper hypoellipticity regularity results for \psivdos\ in terms of suitably weighted Sobolev spaces (see~\cite{FS:EDdbarbCAHG}, 
    \cite{Po:CPDE1}). When $P$ is a differential operator and the Levi form is non-vanishing these results correspond to the maximal hypoellipticity 
    of $P$ as in~\cite{HN:HMOPCV}. 
\end{remark}

% Given a $H$-frame $X_{0},\ldots,X_{d}$ and a multi-order $\alpha$ we let $X^{\alpha}
% \begin{definition}
%     Let $P\in \pvdo^{m}(M,\cE)$, $m \in \N$. We say that $P$ is hypoelliptic maximal when for any $a \in M$, given a $H$-frame $X_{0},..,X_{d}$ near 
%     $a$, for any $u \in \cE'(M,\cE)$ and any $s \in \R$, we have 
%     \begin{equation}
%        \text{$Pu$ is $L^{2}_{s}$ near $a$} \ \Longrightarrow \  \text{$X^{\alpha}u$ is $L^{2}_{s}$ near $a$ for any $\alpha\in \N^{d}$ such that 
%        $\brak \alpha\leq m$},
% %         \label{eq:}
%     \end{equation}
%     where we have let $X^{\alpha}=X_{0}^{\alpha_{0}}\ldots X_{d}^{\alpha_{d}}$. 
% \end{definition}
% 
% It is easy to check that maximal hypoellipticity implies hypoellipticity with loss of $2[\frac{m}{2}]$-derivatives. Moreover, at least 
% 
% 
% When $M$ is compact, combining this with 
% the compactness of the embedding of  $H^{k/2}(M,\cE)$ into $L^{2}(M,\cE)$ we get: 
% \begin{proposition}\label{prop:PsiHDO.spectrum}
%     Suppose $M$ compact and assume that $P$ has an invertible principal symbol and a 
%     spectrum different from $\C$. Then: \smallskip 
% 
%     1) The spectrum of $P$ consists of isolated eigenvalues with finite multiplicities.\smallskip 
% 
%     2) For any $\lambda\in \op{Sp}P$ the eigenspace $\ker (P-\lambda)$ is a finite dimensional subspace of $C^{\infty}(M,\cE)$.  
% \end{proposition}

Now, let $P:C^{\infty}(M,\cE)\rightarrow C^{\infty}(M,\cE)$ be a \psido\ of order $m$ and assume that $M$ is endowed with a positive density and 
$\cE$ with a Hermitian metric.  Let $P^{a}$ be the model operator of $P$ at 
a point $a\in M$ 
and let  $\pi: G\rightarrow \cH_{\pi}$  be a (nontrivial) unitary representation of $G=G_{a}M$. We define the symbol 
$\pi_{P^{a}}$ as follows (see also~\cite{Ro:HHGRTC}, \cite{Gl:SSGMCAINGHG}, \cite{CGGP:POGD}).

Let $\cH_{\pi}^{0}(\cE_{a})$ be the subspace of $\cH_{\pi}(\cE_{a}):=\cH_{\pi}\otimes \cE_{a}$ spanned by the vectors of the form
\begin{equation}
     \pi_{f}\xi=\int_{G}(\pi_{x}\otimes 1_{\cE_{a}})(\xi \otimes f(x)) dx, %\qquad \xi \in \cH_{\pi}(\cE_{a}), f \in % \cS_{0}(G), \quad .
\end{equation}
as $\xi$ ranges over  $\cH_{\pi}$ and $f$ over $\cS_{0}(G,\cE_{a})=\cS_{0}(G)\otimes \cE_{a}$.
Then we let $\pi_{P^{a}}$ denote the (unbounded) operator of $\cH_{\pi}(\cE_{a})$ with domain 
$\cH_{\pi}^{0}(\cE_{a})$ such that
\begin{equation}
    \pi_{P^{a}}(\pi_{f}\xi)=\pi_{P^{a}f}\xi \qquad \forall f\in \cS_{0}(G,\cE_{a})\quad \forall \xi \in \cH_{\pi}.
\end{equation}

One can check that $\pi_{P^{a*}}$ is the adjoint of $\pi_{P^{a}}$ on $\cH_{\pi}^{0}$, hence 
is densely defined. Thus $\pi_{P^{a}}$ is closeable and we can let $\overline{\pi_{P^{a}}}$ denotes its closure. 

In the sequel we let $C^{\infty}_{\pi}(\cE_{a})=C^{\infty}_{\pi}\otimes \cE_{a}$, where $C^{\infty}_{\pi}\subset \cH_{\pi}$ denotes the space of 
smooth vectors of $\pi$ (i.e.~the subspace of vectors $\xi \in \cH_{\pi}$ so that $x \rightarrow \pi(x)\xi$ is smooth from $G$ to $\cH_{\pi}$).

\begin{proposition}[\cite{CGGP:POGD}]\label{PsiHDO.properties-symbol-representation}
    1) The domain of $\overline{\pi_{P^{a}}}$ always contains $C^{\infty}_{\pi}(\cE_{a})$.\smallskip 
   
   2)  If $\Re m \leq 0$ then the operator $\overline{\pi_{P^{a}}}$ is bounded.\smallskip 
   
   3) We have $\overline{(\pi_{P^{a}})^{*}}=(\overline{\pi_{P^{a}}})^{*}$.\smallskip 
   
   4) If $P_{1}$ and $P_{2}$ are \psidos\ on $M$ then $\overline{\pi_{(P_{1}P_{2})^{a}}}=\overline{\pi_{P^{a}_{1}}}\, \overline{\pi_{P^{a}_{2}}}$.  
\end{proposition}

\begin{remark}
   If $\cE_{a}=\C$ and $P^{a}$ is a differentiable operator then, as it is left-invariant, $P^{a}$ belongs to the enveloping algebra 
$\cU(\fg)$ of the Lie algebra $\fg=\fg_{a}M$ of $G$. In this case $\overline{\pi_{P^{a}}}$ coincides on $C^{\infty}_{\pi}$ with the operator 
$d\pi(P^{a})$, where $d\pi$ is the representation of $\cU(\fg)$ induced by $\pi$.  
\end{remark}

\begin{definition}
    We say that $P$ satisfies the Rockland condition at $a$ if for any nontrivial unitary irreducible representation 
$\pi$ of $G_{a}M$ the operator $\overline{\pi_{P^{a}}}$ is injective on $C^{\infty}_{\pi}(\cE_{a})$.
\end{definition}

Since $G=G_{a}M\simeq \bH^{2n+1}\times \R^{d-2n}$ with $2n=\rk \cL_{a}$, there are left-invariant vector fields $X_{0},\ldots,X_{d}$ on $G$ 
such that $X_{0},\ldots,X_{2n}$ are the left-invariant vector fields on $\bH^{2n+1}$ given by~(\ref{eq:Examples.Heisenberg-left-invariant-basis}) 
and $X_{k}=\frac{\partial}{\partial x_{k}}$ for 
$k=2n+1,\ldots,d$. Then, up to unitary 
equivalence, the nontrivial irreducible representations of $G$  are of two types:\smallskip

(i) Infinite dimensional representations $\pi^{\lambda,\xi}:G\rightarrow L^{2}(\R^{n})$ parametrized by $\lambda\in \R\setminus0$ and 
$\xi=(\xi_{2n+1},\ldots,\xi_{2n})$ such that 
\begin{gather}
    d\pi^{\lambda,\xi}(X_{0})=i\lambda|\lambda|, \qquad d\pi^{\lambda,\xi}(X_{k})=i\lambda \xi_{k}, \quad k=2n+1,\ldots,d,\\
    d\pi^{\lambda,\xi}(X_{j})=|\lambda|\frac{\partial}{\partial \xi_{j}}, 
    \qquad d\pi^{\lambda,\xi}(X_{n+j})=i\lambda \xi_{j}, \quad j=1,\ldots,n.
%     \label{eq:}
\end{gather}
Moreover, in this case we have $C^{\infty}(\pi^{\pm,\xi})=\cS(\R^{n})$.

% for $j=1,\ldots,n$ and $k=2n+1,\ldots,d$.\smallskip 
% and given by 
% \begin{equation}
%     (\pi^{(\lambda,\xi)}(x_{0},\ldots,x_{d})f)(\xi_{1},\ldots,\xi_{n})=e^{i(\lambda x_{0}+x_{1}.(\xi_{1}+ \frac{1}{2}x_{n+1})+
%     \ldots+x_{n}.(\xi_{n}+\frac{1}{2}x_{2n})+x_{2n+1}.\xi_{2n+1}+\ldots+x_{d}.\xi_{d})}f(x_{1}+\xi_{n+1},\ldots,x_{n}+\xi_{2n}),
% %     \label{eq:}
% \end{equation}

(ii) One dimensional representations $\pi^{\xi}:G\rightarrow \C$ indexed by $\xi=(\xi_{1},\ldots,\xi_{d}) \in \R^{d}\setminus 0$ such that
\begin{equation}
     d\pi^{\xi}(X_{0})=0, \qquad d\pi^{\xi}(X_{j})=i\xi_{j}, \qquad j=1,\ldots,d.
%     \label{eq:¥}
\end{equation}

 In particular, if $P=p_{m}(-iX)$ with $p \in S_{m}(G)$ then the homogeneity of the symbol $p$ implies that we have 
 $\pi^{\lambda,\xi}_{P}=|\lambda|^{m}\pi^{\pm,\xi}_{P}$ where 
 $\overline{\pi^{\pm,\xi}_{P}}=\overline{\pi^{\pm 1,\xi}_{P}}$ accordingly with the sign of $\lambda$.

 On the other hand,  for the representations in (ii) we have $\overline{\pi^{\xi}_{P}}=\pi^{\xi}_{P}=p_{m}(0,\xi_{1},\ldots,\xi_{d})$. 
 Therefore, we get: 
 
 \begin{proposition}\label{prop:Rockland.reduction}
     The Rockland condition for $P=p_{m}(-iX^{a})$, $p_{m}\in S_{m}(G_{a}M)$, holds, and only if, the following two conditions are 
     satisfied:\smallskip 

    (i) The operators $\overline{\pi^{\pm,\xi}_{P}}$, $\xi \in \R^{d-2n}$, are injective on $\cS(\R^{n})$;\smallskip 
 
   (ii) The restriction of the symbol $p_{m}$ to $\{0\}\times (\R^{n}\setminus 0)\simeq H_{a}^{*}\setminus 0$ is pointwise invertible.
 \end{proposition}
\begin{remark}
   In the case of $G_{a}M=\bH^{2n+1}$ the conditions~(i) and (ii) have those considered by Taylor~\cite{Ta:NCMA}. 
\end{remark}
 
  Next, if $P \in \pvdo^{m}(M,\cE)$ has an invertible principal symbol, hence admits a parametrix $Q\in \pvdo^{-m}(M,\cE)$, then  
 for any $a \in M$ we have $Q^{a}P^{a}$ and $(Q^{a})^{t}(P^{a})^{t}$ are equal to $1$ on $\cS_{0}(G_{a}M,\cE_{a})$ and $\cS_{0}(G_{a}M,\cE_{a}^{*})$ 
 respectively. It then follows from Proposition~\ref{PsiHDO.properties-symbol-representation} 
 that for any nontrivial irreducible unitary representation $\pi$ of $G_{a}$ the operators 
 $\overline{\pi_{P^{a}}}$ and $\overline{\pi_{(P^{a})^{t}}}$ are injective on $C^{\infty}(\pi)$. Thus $P$ and $P^{t}$ satisfy the Rockland condition 
 at every point of $M$. 
 
Conversely, if at some point $a \in M$ the model operator $(P^{a})^{t}$ satisfies the Rockland condition then $(P^{a})^{t}$ is hypoelliptic 
(see~\cite{Ro:HHGRTC}, \cite{Be:OIHGLN}, \cite{Gl:SSGMCAINGHG}; see also~\cite{HN:HGNRN3}, \cite{HN:COHHIGGLG}, \cite{CGGP:POGD}). 
The latter fact then implies that $P^{a}$ admits a 
fundamental solution $k^{a}\in \cS'(G_{a}M,\cE_{a})$ so that $P^{a}k^{a}=\delta_{a}$ (see~\cite{Fo:SEFNLG}, \cite{Ge:LTHG}, \cite{CGGP:POGD}). 
In particular, the inverse Fourier transform of $k^{a}$ yields an inverse for the symbol of $P^{a}$.

This shows that if $P$ and $P^{t}$ satisfies the Rockland condition at every point then for any $a\in M$ there exists $q_{-m}^{a}\in 
S_{-m}(\fg^{*}_{a}M,\cE_{a})$ so that $q^{a}*^{a}\sigma_{m}(P)(a,.)=\sigma_{m}(P)(a,.)*^{a}q^{a}=1$. However, it is an open problem to determine whether 
$q^{a}$ depends smoothly $a$ and so define an element of $S_{m}(\fg^{*}M,\cE)$.

The above issue is at least true in the case of sublaplacians (see~\cite{BG:CHM} and next section). As we shall now see this result can be 
extended to arbitrary \psivdos\ when the Levi form of $(M,H)$ has constant rank. More precisely, we have:

\begin{proposition}\label{thm:PsiHDO.Rockland-Parametrix}
    Suppose that the Levi form  of $(M,H)$ has constant rank and let $P \in \pvdo^{m}(M,\cE)$. Then the following are equivalent:\smallskip 
    
    (i) $P$ and $P^{t}$ satisfy the Rockland condition at every point of $M$;\smallskip 
    
   (ii) $P$ and $P^{*}$ satisfy the Rockland condition at every point of $M$;\smallskip 
    
    (iii) The principal symbol of $P$ is invertible.\smallskip 
    
   \noindent  In particular,  if (i) holds then $P$ admits a parametrix in $\pvdo^{-m}(M,\cE)$ and is hypoelliptic with loss of $\frac{1}{2}\Re m$
   derivatives.   
\end{proposition}
\begin{proof}
   We saw that (i) implies (iii)  already. The same argument shows that (ii) implies (i). Thus we only need to prove the converse statements.  
   
%     If the principal symbol of $P$ is invertible then 
%    $P$ admits a parametrix $Q\in \pvdo^{-m}(M,\cE)$ so that $QP=1 \bmod \Psi^{-\infty}(M,\cE)$. Therefore, it follows from 
%    Proposition~\ref{prop:PsiHDO.composition2}   and 
%    Proposition~\ref{PsiHDO.properties-symbol-representation} that, for any $x \in M$ and any non-trivial unitary representation $\pi$ of $G_{x}M$ we have 
% \begin{equation}
%     \overline{\pi_{Q^{x}}} \overline{\pi_{P^{x}}}= \overline{\pi_{(PQ)^{x}}}= \overline{\pi_{1}}=1,
% \end{equation}
% so that $\overline{\pi_{P^{x}}}$ is injective on $C^{\infty}(\pi)$. Hence $P^{x}$ satisfies the Rockland condition at every point $x\in M$. Similarly, as  the principal 
% symbol $\sigma_{-m}(P^{t})$ is invertible with inverse
% $\sigma_{-m}(Q^{t})$ we see that $P^{t}$ sastifies the Rockland condition at every point of $M$. 
    
    Assume that $P$ and $P^{t}$ satisfies the Rockland condition at every point. For sake of simplicity let us further assume that $\cE$ 
    is the trivial line bundle, so that $P$ and $P^{t}$ are scalar operators. We need to show that the principal symbol of $P$ is invertible. As  
    this is a local issue it is enough to prove it in a local Heisenberg chart $U \subset \Rd$, so we may assume that $P$ and $P^{t}$ are \psivdos\ 
    on $U$.  
    
    Moreover, since the Levi form has constant rank it follows from~\cite[Prop.~2.8]{Po:Pacific1} that $GM$ is a fiber bundle of Lie group with fiber 
    $G=\bH^{2n+1}\times \R^{d-2n}$, where $2n$ is the rank of the Levi form. Therefore, by considering a trivialization of this fiber bundle by means 
    of a suitable local $H$-frame (see~\cite[p.~5]{Po:Pacific1}) we may further assume that $GU$ is the trivial bundle $U\times G$. In particular,  
    the families of model operators $(P^{x})_{x \in U}$ and $((P^{t})^{x})_{x\in U}$ can be seen as smooth families of left-invariant \psivdos\ on 
    $G$ as in~\cite{CGGP:POGD}.  
    
    Now, since $P^{x}$ satisfies the Rockland condition for every $x\in U$ 
    it follows from~\cite[Thm.~5(d)]{CGGP:POGD} near every $x_{0} \in U$ there exists  an open neighborhood $V\subset U$ of $x_{0}$ and 
    a smooth family $(K^{x})_{x \in V}\subset 
    \cK_{m-d-2}(G)$ such that if for $x \in V$ we let $Q^{x}$ be the left-convolution operator with $K^{x}$ acting on $\cS_{0}(G)$ then $Q^{x}$ is a 
    left-inverse for $P^{x}$. 
    
    As $(P^{t})^{x}=(P^{x})^{t}$ satisfies the Rockland condition for every $x \in V$ the same arguments show that $(P^{x})^{t}$ is left-invertible on 
    $\cS_{0}(G)$ for any $x$ in an open neighborhood $W$ of $x_{0}$ contained in $V$. Thus $P^{x}$ is invertible with two-side inverse $Q^{x}$. 
    
    Since $Q^{x}$ is the left-invariant \psivdo\ with symbol $q_{-m}^{x}=(K^{x})^{\wedge}$ in $S_{-m}(G)$ we see that we have 
    \begin{equation}
        q_{-m}^{x}*\sigma_{m}(P)(x,.)=\sigma_{m}(P)(x,.)*q_{-m}^{x}=1 \qquad \forall x \in W.
         \label{eq:PsiHDO.inverse-sigma-m-P-x}
    \end{equation}
    In fact, as $K^{x}$ depends smoothly on $x$, that is, yields an element of $\cK_{-m}(W\times G)$, we get a symbol in $S_{-m}(W\times \fg^{*})$ by 
    letting $q_{-m}(x,\xi)=q_{-m}^{x}(\xi)$. Then~(\ref{eq:PsiHDO.inverse-sigma-m-P-x}) 
    shows that $q_{-m}$ is an inverse for $\sigma_{m}(P)$ on $W\times (\fg^{*}\setminus 0)$. 
    This shows that $\sigma_{m}(P)$ is invertible near every point of $U$, so $\sigma_{m}(P)$ is an invertible symbol. 
    
    Next, since the aforementioned result of~\cite{CGGP:POGD} remains valid \emph{mutatis standis} for systems, by working in a local trivializing 
    Heisenberg chart we can similarly show that, in the case $\cE$ is a general bundle, if $P$ and $P^{t}$ satisfy the Rockland condition at 
    every point then the the principal symbol of $P$ is invertible. Thus the assertions (i) and (ii) are equivalent. 
    
    Finally, the above arguments remain valid when we replace  the transpose $P^{t}$  by the adjoint of $P$, so the statements (ii) and (iii) are 
    equivalent.
\end{proof}

In particular, when the Levi form of $(M,H)$ has constant rank and $P$ is selfadjoint the Rockland condition for $P$ is equivalent to the 
invertibility of the principal symbol of $P$, hence is equivalent to the existence of a parametrix in $\pvdo^{-m}(M,\cE)$.  
% 
% \begin{proposition}\label{thm:PsiHDO.Rockland-Parametrix-selfadjoint}
%     Suppose that the Levi form of $(M,H)$ has constant rank and that $P$ is selfadjoint. Then the following are equivalent:\smallskip 
%     
%     (i) $P$ satisfies the Rockland condition at every point of $M$;\smallskip 
%     
%     (ii) The principal symbol of $P$ is invertible.\smallskip 
%     
%    \noindent  In particular,  if (i) and (ii) hold then $P$ admits a parametrix in $\pvdo^{-m}(M,\cE)$ and is hypoelliptic with loss of $\frac{1}{2}\Re m$
%    derivatives.   
% \end{proposition}

\section{Hypoellipticity criteria for sublaplacians}
\label{sec:sublaplacian}
The main focus of this section is on sublaplacians, which furnish several important examples of operators on Heisenberg manifolds. The scalar case was 
dealt with in~\cite{BG:CHM}, but the results were not extended to sublaplacians acting on sections of vector bundles. These extensions are necessary in 
order to deal with sublaplacians acting on forms such as the Kohn Laplacian or the horizontal sublaplacian (see next section). 

In this section, after having explained the scalar case from the point of view this paper, we extend the results to the non-scalar case.   
In particular, this will allow us to complete the treatment of the Kohn Laplacian in~\cite{BG:CHM} (see Remark~\ref{rem:Examples.Boxb} ahead). 
% this section we focus on sublaplacians, which allows us to deal with many examples. First, we reformulate the known results for scalar 
% sublaplacians, then we explain how to extend these results This extension was not carried out 
% in~\cite{BG:CHM}, but this is a necessary step in order 
% 
% examples of operators on Heisenberg manifolds are sublaplacians (see next section). In this section we explain how the results of the previous 
% sections can be used to recover the known hypoellipticity criterions for sublaplacians.  
% \section{Examples of hypoelliptic operators on Heisenberg manifolds}\label{sec:Examples}
% In this section we review from the point of view of this paper criterions insuring us the invertibility of the principal symbol for the main geometric 
% operators on Heisenberg manifolds: H\"ormander's sum of squares, Kohn Laplacian, horizontal sublaplacian and its conformal powes, contact Laplacian. 
% All these operators but the last one these operators are sublaplacians or powers of sublaplacians, so we start by dealing with general sublaplacians 
% on Heisenberg manifolds. 
% 
% We will use the following definition for sublaplacians.
\begin{definition}\label{def:Examples.sublaplacians}
A differential operator $\Delta:C^{\infty}(M,\cE)\rightarrow C^{\infty}(M,\cE)$ is a sublaplacian when, near any point $a\in M$, we can put 
$\Delta$ in the form,
%  where which $\Delta$ takes the form
 \begin{equation}
    \Delta=-(X_{1}^{2}+\ldots+X_{d}^{2})- i\mu(x) X_{0}+ \op{O}_{H}(1),
%     \sum_{j=1}^{d} b{j}(x)X_{j}+ c(x), \quad \mu(x)=(\mu_{1}(x),\ldots,\mu_{r}(x)).
    \label{eq:Heisenberg.sublaplacian.bundle}
\end{equation}
where $X_{0},X_{1},\ldots,X_{d}$ is a local $H$-frame of $TM$, the coefficient  $\mu(x)$  is a local section of $\End \cE$ 
and the notation $\op{O}_{H}(1)$  means a differential operator of Heisenberg order~$\leq 1$.
% $\mu_{1}(x),\ldots,\mu_{d}(x)$ are smooth complex-valued functions and $b_{1}(x),\ldots,b_{d}(x)$ and $c(x)$ are smooth functions with 
% values in $M_{r}(\C)$.
\end{definition}

Let us look at the Rockland condition for a \emph{scalar} sublaplacian $\Delta:C^{\infty}(M)\rightarrow C^{\infty}(M)$. Let $a\in M$ and let 
$X_{0},X_{1},\ldots,X_{d}$ be a local $H$-frame of $TM$ so that near $a$ we can write 
\begin{equation}
    \Delta=-\sum_{j=1}^{d} X_{j}^{2} - i\mu(x) X_{0}+ \op{O}_{H}(1),
     \label{eq:Heisenberg.sublaplacian.scalar}
\end{equation}
where $\mu(x)$ is a smooth function near $a$. Using~(\ref{eq:Pincipal.example}) we see that the principal symbol of $\Delta$ is 
\begin{equation}
    \sigma_{2}(\Delta)(x,\xi)=|\xi'|^{2}+\mu(a)\xi_{0}, \qquad \xi'=(\xi_{1},\ldots,\xi_{d}).
%     \label{eq:¥}
\end{equation}
In particular we have $ \sigma_{2}(\Delta)(x,0,\xi')=|\xi'|^{2}>0$ for $\xi'\neq 0$, which shows that the condition (i) of 
Proposition~\ref{prop:Rockland.reduction}  is 
always satisfied.

Let $L(x)=(L_{jk}(x))$ be the matrix of the Levi form $\cL$ with respect to the $H$-frame $X_{0},\ldots,X_{d}$, so that for $j,k=1,\ldots,d$ we have
\begin{equation}
    \cL(X_{j},X_{k})=[X_{j},X_{k}]=L_{jk}X_{0} \quad \bmod H.
%     \nonum
\end{equation}
Equivalently, if we let $g(x)$ be the metric on $H$ making orthonormal the frame $X_{1},\ldots,X_{d}$, then for any sections $X$ and $Y$ of $H$ we have 
\begin{equation}
    \cL(X,Y)=g(x)(L(x)X,Y)X_{0} \quad \bmod H.
%     \label{eq:¥}
\end{equation}

The matrix $L(x)$ is antisymmetric, so up to an orthogonal change of frame of $H$, which does not affect the form~(\ref{eq:Heisenberg.sublaplacian.scalar}), 
we may assume that $L(a)$ is in the normal form, 
\begin{equation}
   \qquad L(a)=\left( 
   \begin{array}{ccc}
       0 & D & 0  \\
       -D & 0 & 0 \\
       0 & 0 & 0
   \end{array}\right), \qquad D=\op{diag}(\lambda_{1},\ldots,\lambda_{n}), \quad \lambda_{j}>0,
     \label{eq:Sublaplacian.Levi-form.normal-form}
\end{equation}
so that $\pm i\lambda_{1},\ldots,\pm \lambda_{2n}, 0,\ldots,0$ are the eigenvalues of $L(a)$ counted with multiplicity, the multiplicity of $0$ being $d-2n$. Then 
the model vector fields $X_{0}^{a},\ldots,X_{d}^{a}$ are: 
\begin{gather}
    X_{0}^{a}=\frac{\partial}{\partial x_{0}}, \qquad X_{k}^{a}=\frac{\partial}{\partial x_{k}}, \qquad k=2n+1,..,d,\\
    X_{j}^{a}=\frac{\partial}{\partial x_{j}}-\frac{1}{2}\lambda_{j}x_{n+j}\frac{\partial}{\partial x_{0}}, \qquad 
    X_{n+j}^{a}=\frac{\partial}{\partial x_{j}}+\frac{1}{2}\lambda_{j}x_{j}\frac{\partial}{\partial x_{0}}, \qquad j=1,\ldots,n. 
\end{gather}
In terms of these vector fields the model operator of $\Delta$ at $a$ is 
\begin{equation}
    \Delta^{a}=-[(X_{1}^{a})^{2}+\ldots+(X_{1}^{a})^{2}]-i\mu(a)X_{0}^{a}.
    \label{eq:Sublaplacian.model-operator}
\end{equation}

Next, under the isomorphism $\phi:\bH^{2n+1}\times \R^{d-2n}\rightarrow G_{a}M$ given by 
\begin{equation}
    \phi(x_{0},\ldots,x_{d})=(x_{0},\lambda_{1}^{\frac{1}{2}}x_{1},\ldots,\lambda_{n}^{\frac{1}{2}}x_{n},\lambda_{1}^{\frac{1}{2}}x_{n+1}, \ldots
    \lambda_{n}^{\frac{1}{2}}x_{2n},x_{2n+1},\ldots,x_{d}), 
%     \label{eq:¥}
\end{equation}
the representations $\pi^{\pm,\xi}=\pi^{\pm 1,\xi}$, $\xi \in\{0\}^{2n}\times \R^{d-2n}$, become the representations of $G_{a}M$ such that
\begin{gather}
    d\pi^{\pm,\xi}(X_{0})=\pm i, \qquad d\pi^{\pm,\xi}(X_{k})=\pm i \xi_{k}, \quad k=2n+1,\ldots,d,\\
    d\pi^{\pm,\xi}(X_{j})=\lambda_{j}^{\frac{1}{2}}\frac{\partial}{\partial \xi_{j}}, 
    \qquad d\pi^{\pm,\xi}(X_{n+j})=\pm i\lambda_{j}^{\frac{1}{2}} \xi_{j}, \quad j=1,\ldots,n,\\
       \pi_{\Delta^{a}}^{\pm,\xi}= d\pi^{\pm,\xi}(\Delta^{a})= 
    \sum_{j=1}^{n}\lambda_{j}(-\partial_{\xi_{j}}^{2}+\xi_{j}^{2}) \pm(\xi_{2n+1}^{2}+\ldots+\xi_{d}^{2}+\mu(a)).
\end{gather}

The spectrum of the harmonic oscillator $\sum_{j=1}^{n}\lambda_{j}(-\partial_{\xi_{j}}^{2}+\xi_{j}^{2})$ is $\sum_{j=1}^{n}\lambda_{j}(1+2\N)$ and 
 all its eigenvectors belong to $\cS(\R^{n})$. Thus, the operator $ \pi_{\Delta^{a}}^{\pm,\xi}$ is injective on $\cS(\R^{n})$ if, and only if, 
  $\xi_{2n+1}^{2}+\ldots+\xi_{d}^{2}+\mu(a)$ is not $\pm \sum_{j=1}^{n}\lambda_{j}(1+2\N)$. This occurs for any $\xi \in\{0\}^{2n}\times 
  \R^{d-2n}$ if, and only if, the following condition holds
%   \begin{equation}
%  
%        \label{eq:Sublaplacian.condition.scalar}
%   \end{equation}
%  where the set $\Lambda^{a}$ is defined by 
 \begin{gather}
         \text{$\mu(a)$ is not in the singular set $\Lambda_{a}$}, 
         \label{eq:Sublaplacian.condition.scalar} \\
         \Lambda_{a}=(-\infty, -\frac12 \Tra |L(a)|]\cup [\frac12 \Tra 
  |L(a)|,\infty) \qquad \text{if $2n<d$},\\
   \Lambda_{a}=\{\pm(\frac12 \Tra |L(a)|+2\sum_{1\leq j \leq n}\alpha_{j}|\lambda_{j}|); \alpha_{j}\in \N^{d}\}\qquad \text{if $2n=d$}.
 \end{gather}
In particular, the condition~(ii) of Proposition~\ref{prop:Rockland.reduction} 
is equivalent to~(\ref{eq:Sublaplacian.condition.scalar}). Since the condition (i) is always satisfied, it follows that the 
Rockland condition for $\Delta$ is equivalent to~(\ref{eq:Sublaplacian.condition.scalar}). 

Notice also that, independently of the equivalence with the Rockland condition, 
the condition~(\ref{eq:Sublaplacian.condition.scalar}) does not depend on the choice of the $H$-frame, because 
as $\Lambda_{a}$ depends only on the eigenvalues of $L(a)$ which scale in the same way as $\mu(a)$ 
under a change of $H$-frame preserving the form~(\ref{eq:Heisenberg.sublaplacian.scalar}). 

On the other hand, since the transpose $(\Delta^{a})^{t}=(\Delta^{t})^{a}$ is given by the formula~(\ref{eq:Sublaplacian.model-operator}) 
 with $\mu(a)$ replaced by $-\mu(a)$, which has no 
effect on~(\ref{eq:Sublaplacian.condition.scalar}), we see that the Rockland condition for $(\Delta^{t})^{a}$ too is equivalent 
 to~(\ref{eq:Sublaplacian.condition.scalar}). Therefore, we have obtained: 

\begin{proposition}\label{prop:Sublaplacian.invertibility.scalar}
The Rockland conditions for $\Delta^{t}$ and $\Delta$ at $a$ are both equivalent to~(\ref{eq:Sublaplacian.condition.scalar}).
\end{proposition}

In particular, we see that if the principal symbol of $\Delta$ is invertible then the condition~(\ref{eq:Sublaplacian.condition.scalar}) holds at every 
point. As shown by Beals-Greiner the converse is true as well. The key result is the following.

\begin{proposition}[{\cite[Sect.~5]{BG:CHM}}]\label{prop:Sublaplacian.inverse.scalar}
    Let $U\subset \Rd$ be a Heisenberg chart near $a$ and set 
    \begin{equation}
        \Omega=\{(\mu,x)\in \C \times U; \ \mu \not \in \Lambda_{x}\}.
%         \label{eq:}
    \end{equation}
     Then $\Omega$ is an 
    open subset of $\C \times U$ and there exists $q_{\mu}(x,\xi)\in C^{\infty}(\Omega\times \Rdo)$ such that:\smallskip
    
   (i) $q_{\mu}(x,\xi)$ is analytic with respect to $\mu$;\smallskip
    
   (ii) We have $q_{\mu}(x,\lambda.\xi)=\lambda^{-2}q_{\mu}(x,\xi)$ for any $\lambda>0$.\smallskip
    
   (iii) For any $(\mu,x)\in \Omega$ the symbol $q_{\mu}(x,.)$ inverts $|\xi'|^{2}+i\mu\xi_{0}$ on $G_{x}U$, i.e., we have
    \begin{equation}
        q_{\mu}(x,.)*^{x}(|\xi'|^{2}+i\mu \xi_{0})=(|\xi'|^{2}+i\mu \xi_{0})*^{x}q_{\mu}(x,.)=1.
%         \label{eq:¥}
    \end{equation}
More precisely, $q_{\mu}(x,\xi)$ is obtained from the analytic continuation of the function,
\begin{gather*}
    q_{\mu}(x,\xi) =\int_{0}^{\infty}e^{-t\mu \xi_{0}}G(x,\xi,t)dt, \qquad |\Re \mu|<\frac{1}{2}\Tr |L(x)|,\\
    G(x,\xi,t)=\det{}^{-\frac{1}{2}}[\cosh (t|\xi_{0}||L(x)|)] \exp[-t\acou{\frac{\tanh (t|\xi_{0}||L(x)|)}{t|\xi_{0}||L(x)|}\xi'}{\xi'}].
%     \label{eq:¥}
\end{gather*}
%   where we have let $\xi'=(\xi_{1},\ldots,\xi_{d})$.   
\end{proposition}

 This implies that if the condition~(\ref{eq:Sublaplacian.condition.scalar}) is satisfied at every point $x \in U$ then we get an inverse $q_{-2}\in 
S_{-2}(\URd)$ for $\sigma_{2}(\Delta)(x,\xi)=|\xi'|^{2}+i\mu(x)\xi_{0}$ on $\URd$ by letting 
\begin{equation}
    q_{-2}(x,\xi)=q_{\mu(x)}(x,\xi), \qquad (x,\xi)\in \URdo.
%     \label{eq:¥}
\end{equation}
It thus follows that if~(\ref{eq:Sublaplacian.condition.scalar}) 
holds at every point of $M$ then the principal symbol of $\Delta$ is invertible near any point of $M$, hence admits 
an inverse in $S_{-2}(\fg^{*}M)$. Therefore, we get: 

\begin{proposition}\label{prop:Sublaplacian.Rockland.scalar}
A scalar sublaplacian $\Delta:C^{\infty}(M)\rightarrow C^{\infty}(M)$ has an invertible principal symbol if, and only if, it satisfies the 
condition~(\ref{eq:Sublaplacian.condition.scalar}) at every point.  
\end{proposition}

Let us now extend the above results to the case of a sublaplacian $\Delta:C^{\infty}(M,\cE)\rightarrow C^{\infty}(M,\cE)$ acting on the sections of 
the vector bundle $\cE$.

Let $a \in M$ and let $X_{0},\ldots,X_{d}$ be a local $H$-frame near $a$ with respect to which we have 
\begin{equation}
    \Delta=-\sum_{j=1}^{d} X_{j}^{2} - i\mu(x) X_{0}+ \op{O}_{H}(1),
     \label{eq:Heisenberg.sublaplacian.bundle2}
\end{equation}
where $\mu(x)$ is a smooth section of $\End \cE$. 

In a suitable basis of $\cE_{a}$ the matrix  of $\mu(a)$ is in  triangular form,
\begin{equation}
    \mu(a)=\left( 
    \begin{array}{ccc}
        \mu_{1}(a) & * & *  \\
        0 & \ddots & *  \\
        0 & 0 & \mu_{r}(a)
    \end{array}
\right).    
    %     \label{eq:}
\end{equation}
where $\mu_{1}(a),\ldots,\mu_{r}(a)$ denote the eigenvalues of $\mu(a)$ counted with multiplicity. 
Therefore, the model operator of $\Delta$ at $a$ is of the form, 
\begin{equation}
    \Delta^{a}= \left( 
    \begin{array}{ccc}
        \Delta_{1}^{a} & * & *  \\
        0 & \ddots & *  \\
        0 & 0 &\Delta_{r}^{a}
    \end{array}
\right), \qquad \Delta_{j}^{a}=-[(X_{1}^{a})^{2}+\ldots+(X_{1}^{a})^{2}]-i\mu_{j}(a)X_{0}^{a}.  
%     \label{eq:¥}
\end{equation}
It follows that $\Delta^{a}$ satisfies the Rockland condition if, and only if, so does each sublaplacian $\Delta_{j}^{a}$, $j=1,\ldots,r$. Using 
Proposition~\ref{prop:Sublaplacian.Rockland.scalar}
we then deduce that the Rockland condition $\Delta^{a}$ is equivalent to the condition, 
\begin{equation}
    \op{Sp}\mu(a) \cap \Lambda_{a}=\emptyset.
     \label{eq:Sublaplacian.condition}
\end{equation}

Notice that the same is true for the transpose $(\Delta^{a})^{t}$. Moreover, the condition~(\ref{eq:Sublaplacian.condition}) 
is independent of the choice of the basis of $\cE_{a}$ or of the $H$-frame since the condition involves $\mu(a)$ only though its eigenvalues of $\mu(a)$ 
and the latter scale in the same way as that of $L(a)$ under a change of $H$-frame preserving the form~(\ref{eq:Heisenberg.sublaplacian.bundle2}). 

Next, concerning the invertibility of the principal symbol of $\Delta$ the following extension of 
Proposition~\ref{prop:Sublaplacian.inverse.scalar} holds. 
% Moreover, as the condition involves $\mu(a)$ only via its spectrum it does dependthe choice of the 
% trivialization near $a$. Moreover, as the condition~(\ref{eq:Sublaplacian.condition.scalar}) this condition does not depend on the choice of the Heisenberg coordinates 
% either, hence makes sense intrinsically. 

\begin{proposition}\label{prop:Sublaplacian.inverse.system}
    Let $U\subset \Rd$ be a trivializing Heisenberg chart near $a$ and set
    \begin{equation}
        \Omega=\{(\mu,x)\in M_{r}(\C) \times U; \ \mu \not \in \Lambda_{x}\}. 
%     \label{eq:}
\end{equation}
Then $\Omega$ is an open subset of $M_{r}(\C) \times U$ and there exists $q_{\mu}(x,\xi)\in C^{\infty}(\Omega\times \Rdo, M_{r}(\C))$ so~that:\smallskip
    
    (i) $q_{\mu}(x,\xi)$ is analytic with respect to $\mu$;\smallskip
    
    (ii) We have $q_{\mu}(x,\lambda.\xi)=\lambda^{-2}q_{\mu}(x,\xi)$ for any $\lambda>0$.\smallskip
    
    (iii) For any $(\mu,x)\in \Omega$ the symbol $q_{\mu}(x,.)$ inverts $|\xi'|^{2}+i\mu\xi_{0}$ on $G_{x}U$, that is, we have
    \begin{equation}
        q_{\mu}(x,.)*^{x}(|\xi'|^{2}+i\mu \xi_{0})=(|\xi'|^{2}+i\mu \xi_{0})*^{x}q_{\mu}(x,.)=1.
%         \label{eq:¥}
    \end{equation}
\end{proposition}
\begin{proof}
  It is enough to prove that near point $(\mu_{0},x_{0})\in \Omega$ there exists an open neighborhood $\Omega'$ contained in $\Omega$ and 
  a function $q_{\mu}(x,\xi)\in C^{\infty}(\Omega\times \Rdo, M_{r}(\C))$ satisfying the properties (i), (ii) and (iii) on $\Omega'\times \Rdo$. 
  
   To this end observe that since $\op{Sp}\mu \subset \overline{D}(0,\|\mu\|)$ for any $\mu \in M_{r}(\C)$, we see that if we let  
 $K=\overline{B(0,\|\mu_{0}\|+1)}$ then any  $\mu \in M_{r}(\C)$ close enough to $\mu_{0}$ has its spectrum contained in $K$. 
 
 Let $\delta>0$ be small enough so that $\delta< \frac{1}{2}\op{dist}(\op{Sp}\mu_{0},\Lambda_{x_{0}})$ and define $V_{1}=\op{Sp} \mu_{0}+D(0,\delta)$ and 
 $V_{2}=\Lambda_{x_{0}}+D(0,\delta)$, so that $V_{1}$ and $V_{2}$ are disjoint open subsets of $\C$ containing $\Sp \mu_{0}$ and $\Lambda_{x_{0}}$ 
 respectively. 
 
 Notice that for any $\mu$ close enough to $\mu_{0}$ we have $\op{Sp} \mu \subset V_{1}$. Otherwise there 
 exists a sequence $(\mu_{k})_{k \geq 1} \subset M_{r}(\C)$ converging to $\mu_{0}$ and a sequence of eigenvalues 
 $(\lambda_{k})_{k \geq 1} \subset K$, $\lambda_{k }\in \Sp \mu_{k}$, such that $\lambda_{k} \not \in V_{1}$ for any $k\geq 1$. 
 Since the sequence $(\lambda_{k})_{k 
 \geq 1}$ is bounded,  we may assume that it converges to some $\lambda \not \in V_{1}$. Necessarily $\lambda$ is an eigenvalue of $\mu_{0}$, which 
 contradicts  the fact that $\lambda \not \in V_{1}$. Thus there exists $\eta_{1}>0$ so that for any $\mu \in B(\mu_{0},\eta_{1})$ 
 we have $\op{Sp} \mu \subset V_{1}$.
 
 Similarly there exists $\eta_{2}>0$ so that for any $x \in B(x_{0},\eta_{2})$ we have $\Sp |L(x)|\subset \Sp |L(x_{0})|+D(0,\delta)$, which implies 
 $\Lambda_{x}\subset \Lambda_{x_{0}}+D(0,\delta)=V_{2}$. Therefore the open set $\Omega'=B(\mu_{0},\eta_{1})\times B(x_{0},\eta_{2})$ is such that 
 for any $(\mu,x)\in \Omega'$ we have $\Sp \mu \cap \Lambda_{x}\subset V_{1}\cap V_{2}=\empty$, that is, $\Omega'$ is an open neighborhood of 
 $(\mu_{0},x_{0})$ contained in $\Omega$. %Hence $\Omega$ is open.
 
Next, let $\Gamma$ be a smooth curve of index 1 such that the bounded connected component of $\C \setminus \Gamma$ contains $V_{1}$ and its unbounded 
component contains $V_{2}$. Then we define an element of $\Hol(B(\mu_{0},\eta_{1}))\hotimes C^{\infty}(B(x_{0},\eta_{2})\times \Rdo)$ by letting 
\begin{equation}
    q_{\mu}(x,\xi) =\frac{1}{2i\pi}\int_{\Gamma}q_{\gamma}(x,\xi)(\gamma-\mu)d\gamma, \qquad (\mu,x,\xi)\in \Omega'\times \Rdo.
%     B(\mu_{0},\eta_{1})\times     B(x_{0},\eta_{2})
%     \label{eq:}
\end{equation}
This function is homogeneous of degree $-2$ with respect to $\xi$ and for any $(\mu,x)\in \Omega'$ we have
\begin{multline}
    q_{\mu}(x,.)*^{x}(|\xi'|^{2}+i\mu\xi_{0})=\frac{1}{2i\pi}\int_{\Gamma}q_{\gamma}(x,.)*^{x}(|\xi'|^{2}+i\mu\xi_{0})(\gamma-\mu)^{-1}d\gamma,\\
    \frac{1}{2i\pi}\int_{\Gamma}[(\gamma-\mu)^{-1}-iq_{\gamma}(x,.)*\xi_{0}]d\gamma=1.
\end{multline}
Similarly we have $(|\xi'|^{2}+i\mu\xi_{0})*^{x}q_{\mu}(x,.)=1$. Thus $q_{\mu}(x,\xi)$ satisfies the properties (i), (ii) and (iii) on $\Omega'\times \Rdo$. 
The proof is thus complete. 
% Since the same is true for $\Omega'$ sufficiently small open neighborhood of any other element of $\Omega$ we deduce that there exists a 
% function $q_{\mu}(x,\xi)\in C^{\infty}(\Omega\times \Rdo)$ satisfying the properties (i), (ii) and (iii).
\end{proof}

In the same way as Proposition~\ref{prop:Sublaplacian.inverse.scalar} 
in the scalar case, Proposition~\ref{prop:Sublaplacian.inverse.system} implies that when the condition~(\ref{eq:Sublaplacian.condition}) holds 
everywhere the principal symbol of $\Delta$ admits an inverse in $S_{-2}(\fg^{*}M,\cE)$. We have thus proved:
% from Proposition~\ref{prop:Sublaplacian.inverse.system} that if the condition~(\ref{eq:Sublaplacian.condition}) is satisfied for every $a \in U$ 
% then the principal symbol $\sigma_{2}(x,\xi)=|\xi'|^{2}+i\mu(x)\xi_{0}$ 
% of $\Delta$ is invertible with inverse $q_{-2}\in S_{-2}(\URd, M_{r}(\C))$ given by $q_{-2}(x,\xi)=q_{\mu(x)}(x,\xi)$. 
% Summarizing all this we have proved: 

\begin{proposition}\label{prop:Sublaplacian.Rockland-bundle}
  1) At every point $a\in M$ the Rockland conditions for $\Delta$ and $\Delta^{t}$ are equivalent to~(\ref{eq:Sublaplacian.condition}).\smallskip
    
    2) The principal symbol of $\Delta$ is invertible if, and only if, the condition~(\ref{eq:Sublaplacian.condition}) holds everywhere. 
    Moreover, when the latter occurs $\Delta$ admits a parametrix in $\pvdo^{-2}(M,\cE)$ and is hypoelliptic with loss of 1 derivative.
\end{proposition}

\section{Examples of hypoelliptic operators on Heisenberg manifolds}\label{sec:Examples}
In this section we explain how the previous results of this paper can used to deal with the hypoellipticity for the main geometric operators on 
Heisenberg manifolds: H\"ormander's sum of squares, Kohn Laplacian, horizontal sublaplacian and  contact Laplacian. In particular, the treatment 
in~\cite{BG:CHM} of the Kohn Laplacian and we establish a criterion for the invertibility of the horizontal sublaplacian, which has not been done 
before. 

% All these operators but the last one these operators are sublaplacians or powers of sublaplacians, so we start by dealing with general sublaplacians 
% on Heisenberg manifolds. 

\subsection{H\"ormander's sum of squares.}
 Let $X_{1},\ldots,X_{m}$ be (real) vector fields on a manifold $M^{d+1}$ and consider the sum of squares,
\begin{equation}
    \Delta=-(X_{1}^{2}+\ldots+X_{m}^{2}).
\end{equation}
By a celebrated theorem of H\"ormander~\cite{Ho:HSODE} the operator $\Delta$ is hypoelliptic provided that the following bracket condition is 
satisfied: the vector fields $X_{0},X_{1},\ldots,X_{m}$ together with their successive Lie brackets $[X_{j_{1}}, [X_{j_{2}},\ldots ,
X_{j_{1}}]\ldots]]$ span the tangent bundle $TM$ at every point. 

When $X_{1},\ldots,X_{m}$ span a hyperplane bundle $H$ the operator $\Delta$ is a sublaplacian with \emph{real} coefficients and the bracket condition 
reduces to $H+[H,H]=TM$, which is equivalent to the non-vanishing of the Levi form of $(M,H)$. 

More generally, given a vector bundle $\cE$, the theorem of H\"ormander holds for sublaplacians $\Delta:C^{\infty}(M,\cE)\rightarrow 
C^{\infty}(M,\cE)$ of the form
\begin{equation}
    \Delta=-(\nabla_{X_{1}}^{2}+\ldots+\nabla_{X_{m}}^{2}) +\op{O}_{H}(1),
     \label{eq:Operators.generalized-sum-of-squares}
\end{equation}
where $\nabla$ is a connection on $\cE$. In particular, if $M$ is endowed with 
a positive density and $\cE$ with a Hermitian metric, this includes the selfadjoint sum of squares,
\begin{equation}
    \Delta=\nabla_{X_{1}}^{*}\nabla_{X_{1}}+\ldots+\nabla_{X_{m}}^{*}\nabla_{X_{m}}.
\end{equation}

In fact, if $\Delta$ is a sublaplacian of the form~(\ref{eq:Operators.generalized-sum-of-squares}) then in~(\ref{eq:Heisenberg.sublaplacian.bundle})
the matrix $\mu(x)$ vanishes, so that~(\ref{eq:Sublaplacian.condition}) holds if, and only if, the Levi form does not vanish at 
$x$. Therefore, we obtain: 

\begin{proposition}
 Let $\Delta:C^{\infty}(M,\cE)\rightarrow C^{\infty}(M,\cE)$  be a (generalized) sum of squares of the form~(\ref{eq:Operators.generalized-sum-of-squares}). 
 Then:\smallskip 
 
 1) At a point $x\in M$ the operators $\Delta$ and $\Delta^{t}$ satisfies the Rockland condition if, and only if, the Levi form $\cL$ does not vanish 
 at $x$.\smallskip
 
 2) The principal symbol of $\Delta$ is invertible if, and only if, the Levi form is non-vanishing. In particular, when the latter occurs $\Delta$ 
 admits a parametrix in $\pvdo^{-2}(M,\cE)$ and  is hypoelliptic with loss of one derivative.
\end{proposition}

In particular, since the nonvanishing of the Levi form is equivalent to the  
bracket condition $H+[H,H]=TM$, we see that, the special case of Heisenberg manifolds, we recover the hypoellipticity result of~\cite{Ho:HSODE} for 
sums of squares.

\subsection{The Kohn Laplacian}
In~\cite{KR:EHFBCM} Kohn-Rossi showed that the Dolbeault complex on a bounded complex domain  induces on 
its boundary a horizontal complex of differential forms. This was later extended by Kohn~\cite{Ko:BCM} to the general setting of a CR 
manifold $M^{2n+1}$ as follows.

Let $M^{2n+1}$ be a CR manifold with CR bundle $T_{1,0}\subset T_{\C}M$ and set $T_{0,1}=\overline{T_{1,0}}$. Then the subbundle   
$H=\Re (T_{1,0}\oplus T_{0,1})\subset TM$ admits an integrable complex structure and the splitting $H\otimes \C=T_{1,0}\oplus T_{0,1}$ 
gives rise to a decomposition  $\Lambda H^{*}\otimes \C = \oplus_{0\leq p,q\leq n } \Lambda^{p,q}$, where $(p,q)$ is called a bidegree of a form with 
values in $\Lambda^{p,q}$.

Assume that $T_{\C}M$ is endowed with a Hermitian metric such 
that $T_{1,0}$ and $T_{0,1}$ are orthogonal subspaces and complex conjugation is an (antilinear) isometry. This Hermitian metric gives rise to a 
Hermitian metric on $\Lambda^{*}T^{*}_{\C}M$ with respect to which the decomposition  $\Lambda H^{*}\otimes \C = \oplus_{0\leq p,q\leq n } 
\Lambda^{p,q}$ becomes orthogonal. Let $\Pi_{p,q}:\Lambda^{*}T_{\C}^{*}M \rightarrow \Lambda^{p,q}$ be the orthogonal projection onto $\Lambda^{p,q}$. 
Then the Kohn-Rossi operator $\bar \partial_{b}:  C^{\infty}(M,\Lambda^{p,q}) \rightarrow C^{\infty}(M,\Lambda^{p,q+1})$ is given by 
\begin{equation}
    \bar\partial_{b} \eta = \Pi_{p,q+1}(d\eta),\qquad  \eta \in C^{\infty}(M,\Lambda^{p,q}). 
\end{equation}
Since the integrability of $T_{1,0}$ implies that $\bar\partial_{b}^{2}=0$, this yields chain complexes 
$\bar\partial_{b}:C^{\infty}(M,\Lambda^{p,*}) \rightarrow C^{\infty}(M,\Lambda^{p,*+1})$. 

Endowing $M$ with a smooth density $\rho>0$ we let  $\bar\partial_{b} ^{*}$ denote 
the formal adjoint of $\bar\partial_{b}$. Then the Kohn Laplacian is   
\begin{equation}
    \Box_{b}=( \bar\partial_{b}+\bar\partial^{*}_{b})^{2}=\bar\partial^{*}_{b}\bar\partial_{b} +  \bar\partial_{b}\bar\partial^{*}_{b}.
\end{equation}
The Kohn Laplacian is a sublaplacian (see~\cite[Sect.~13]{FS:EDdbarbCAHG}, \cite[Sect.~20]{BG:CHM}) and it was shown by Kohn~\cite{Ko:BCM} 
that under the condition $Y(q)$ this operator is hypoelliptic with loss of one derivative when acting on $(p,q)$-forms. 

The condition $Y(q)$ means that if for $x \in M$ we let $(r(x)-\kappa(x),\kappa(x),n-r(x))$ be the signature of $L_{\theta}$ at $x$, so that $r(x)$ is 
the rank of $L_{\theta}$ and $\kappa(x)$ the number of its negative eigenvalues, then the condition $Y(q)$ is satisfied at $x \in M$ when we have 
\begin{equation}
    q\not \in \{\kappa(x),\kappa(x)+1,\ldots,\kappa(x)+n-r(x)\}\cup \{r(x)-\kappa(x),r(x)-\kappa(x)+1,\ldots,n-\kappa(x)\}.
     \label{eq:Operators.Y(q)-condition}
\end{equation}
For instance,  when $M$ is $\kappa$-strictly pseudoconvex, the $Y(q)$-condition exactly means that we must have $q\neq \kappa$ and $q\neq n-\kappa$. 

In fact, as shown in~\cite[Sect.~21]{BG:CHM}, at every point $a\in M$ the condition $Y(q)$ is equivalent to the 
condition~(\ref{eq:Sublaplacian.condition}) for the Kohn Laplacian acting on $(p,q)$-forms. Therefore, using Proposition~\ref{prop:Sublaplacian.Rockland-bundle} we 
immediately get:  

\begin{proposition}\label{prop:Examples.Boxb}
    Let $\Box_{b}:C^{\infty}(M,\Lambda^{p,q})\rightarrow C^{\infty}(M,\Lambda^{p,q})$ be the Kohn Laplacian acting on $(p,q)$-forms. 
   
   1) At a point $x\in M$ the Rockland condition for $\Box_{b}$ is equivalent to the condition~$Y(q)$.\smallskip
   
   2) The principal symbol of $\Box_{b}$ is invertible if, and only if, the condition $Y(q)$ is satisfied at every point. In particular, when the latter 
   occurs $\Box_{b}$ admits a parametrix in $\pvdo^{-2}(M,\Lambda^{p,q})$ and is hypoelliptic with loss of one derivative.
\end{proposition}

\begin{remark}\label{rem:Examples.Boxb}
    The proof of the second part of the statement in~\cite[Sect.~21]{BG:CHM} is not quite complete, because Beals-Greiner claim that diagonalizing the 
    leading part of the Kohn Laplacian allows us to use to apply the results for scalar sublaplacians. While this property is true in case of a 
    Levi Metric (see~\cite{FS:EDdbarbCAHG}), it may fail for a general metric on $T_{\C}M$ since a smooth eigenframe needs not exists
    For instance, for the Kohn Laplacian the eigenvalues of the 
    matrix $\mu(x)$ in~(\ref{eq:Heisenberg.sublaplacian.bundle}) with respect to an orthonormal $H$-frame of $TM$ are given in terms of eigenvalues of the Levi 
    form (see~Eq.~(21.31) in~\cite{BG:CHM}), but the latter need not depend smoothly on $x$ (unless the metric on $T_{\C}M$ is a Levi metric). 
%    diagonalization  of 
%     the leading part of $\Box_{b}$ it is possible to apply Proposition~\ref{prop:Sublaplacian.invertibility.scalar}, 
%     which deals with the scalar case only. However, this not quite true because in an eigenframe 
%     the coefficient $\mu(x)$ in the form~(\ref{eq:Heisenberg.sublaplacian.bundle}) 
%     of $\Box_{b}$ depend on the eigenvalues of the Levi form, which need not be smooth functions
   Therefore, in order to deal  with the Kohn Laplacian acting on forms, we really need to use the version for sublaplacians acting on section of vector 
   bundles, provided by Proposition~\ref{prop:Sublaplacian.Rockland-bundle}, but not deal with in~\cite{BG:CHM}.
\end{remark}

\begin{remark}
The $Y(q)$-condition is only a sufficient condition for the 
 hypoellipticity of $\Box_{b}$, as the latter may be hypoelliptic when the $Y(q)$-condition fails 
(e.g.~\cite{Ko:SEPDCRM}, \cite{Ko:OMSHETCROBL}, \cite{Ni:PhD}).
\end{remark}

\subsection{The horizontal sublaplacian on a Heisenberg manifold}
Let $(M^{d+1},H)$ be a Heisenberg manifold endowed with a Riemannian metric and let $\Lambda^{*}_{\C}H^{*}=\oplus_{k=0}^{d}\Lambda^{k}_{\C}H^{*}$ 
be the (complexified) bundle of horizontal forms. Then the horizontal sublaplacian 
$\Delta_{b}:C^{\infty}(M,\Lambda^{*}_{\C}H^{*})\rightarrow C^{\infty}(M,\Lambda^{*}_{\C}H^{*})$ is given by 
\begin{equation}
    \Delta_{b}=d_{b}^{*}d_{b}+d_{b}d_{b}^{*}, \qquad d_{b}\alpha=\pi_{b}(d\alpha), 
\end{equation}
where $\pi_{b}$ denotes the orthogonal projection of $\Lambda_{\C}^{*}T^{*}M$ onto 
$\Lambda^{*}_{\C}H^{*}$. 

This operator was first introduced by Tanaka~\cite{Ta:DGSSPCM} for strictly pseudoconvex CR manifolds, but versions of this operator acting on functions 
were independently defined by Greenleaf~\cite{Gr:FESPM} and Lee~\cite{Le:FMPHI}. Moreover, it can be shown that $d_{b}^{2}=0$ if, and only if, 
the subbundle $H$ is integrable, so in general $\Delta_{b}$ is not the Laplacian of a chain complex. 

% On functions $\Delta_{b}$ is a sum of squares modulo lower order terms. When acting on horizontal forms of higher degree, that is, on sections of 
% $\Lambda^{k}_{\C}H^{*}$ with $k\geq 1$, the operator $\Delta_{b}$ is not quite a sublaplacian in the sense of 
% Definition~\ref{def:Examples.sublaplacians}. 
% More precisely, locally $\Delta_{b}$ is of the form~(\ref{eq:Heisenberg.sublaplacian.bundle})   
% with a matrix $\mu(x)$ which is diagonalizable at every point, but need to be so in a  frame of $\Lambda^{k}_{\C}H^{*}$,  
% although this occurs at least when the Levi form of $(M,H)$ has constant rank or when $H$ admits an almost complex structure, hence when 
% $M$ is a CR manifold or a contact manifold (see Appendix). 
% 
% As explained in Appendix, it is not difficult extend the results of~\cite{BG:CHM} to deal with the horizontal sublaplacian. In this case, 

As we shall see to determine under which condition the principal symbol of $\Delta_{b}$ is invertible the relevant condition to 
look at is the $X(k)$-condition below. 

\begin{definition}
    For $x\in M$ let $2r(x)$ denote the rank of the Levi form $\cL$ at $x$. Then we say that $\cL$ satisfies the condition 
$X(k)$ at $x$ when we have
\begin{equation}
    k\not \in\{r(x),r(x)+1,\ldots,d-r(x)\}.
    \label{eq:Operators.X(k)-condition}
\end{equation}
\end{definition}

For instance, the condition $X(0)$ is satisfied if, and only if, the Levi form does not vanish. Also, 
if $M^{2n+1}$ is a contact manifold or a nondegenerate CR manifold then the Levi form is everywhere nondegenerate, so $r(x)=2n$ and 
the $X(k)$-condition becomes $k\neq n$. In any case, we have:

\begin{proposition}\label{prop:Examples.horizontal-sublaplacian}
    Let $\Delta_{b}:C^{\infty}(M,\Lambda^{k}_{\C}H^{*})\rightarrow C^{\infty}(M,\Lambda^{k}_{\C}H^{*})$ be the horizontal  sublaplacian acting on 
    horizontal forms of degree $k$. 
   
   1) At a point $x\in M$ the Rockland condition for $\Delta_{b}$ is equivalent to the condition~$X(k)$.\smallskip
   
   2) The principal symbol of $\Delta_{b}$ is invertible if, and only if, the condition $X(k)$ is satisfied at every point. In particular, when the latter 
   occurs $\Delta_{b}$ admits a parametrix in $\pvdo^{-2}(M,\Lambda^{k}_{\C}H^{*})$ and is hypoelliptic with loss of one derivative.
\end{proposition}
\begin{proof}
    First, thanks to Proposition~\ref{prop:Sublaplacian.Rockland-bundle}
    we only have to check that for $k=0,\ldots,d$ at any point $a$  the condition~(\ref{eq:Sublaplacian.condition}) for 
    $\Delta_{b|_{\Lambda^{k}_{\C}H^{*}}}$ is equivalent to the condition $X(k)$. 
    
    Next, let $U\subset \Rd$ be a Heisenberg chart around $a$ together with an orthonormal $H$-frame $X_{0},X_{1},\ldots X_{d}$ of $TU$. 
    Let $g$ be the Riemannian metric of $M$. Then on $U$ we can write the Levi form $\cL$ in the form, 
    \begin{equation}
        \cL(X,Y)=[X,Y]=\acou{L(x)X}{Y}X_{0} \quad \bmod H, 
%         \label{eq:}
    \end{equation}
for some antisymmetric section $L(x)$ of $\End_{\R}H$. In particular, if for $j,k=1,\ldots,d$ we let $L_{jk}=\acou{LX_{j}}{X_{k}}$ then 
we have 
\begin{equation}
    [X_{j},X_{k}]=L_{jk}X_{0} \qquad \bmod H. 
     \label{eq:Examples.Levi-form.coefficients}
\end{equation}

Let $2n$  be the rank of $L(a)$. Since the condition~(\ref{eq:Sublaplacian.condition}) for $\Delta_{b}$ at $a$ is independent of the choice of the Heisenberg chart, 
we may assume that $U$ is chosen in such way that at $x=a$ we have $g(a)=1$ and $L(a)$ is in the normal form, 
\begin{equation}
   \qquad L(a)=\left( 
   \begin{array}{ccc}
       0 & D & 0  \\
       -D & 0 & 0  \\
       0 & 0 & 0
   \end{array}\right), \qquad D=\op{diag}(\lambda_{1},\ldots,\lambda_{n}), \quad \lambda_{j}>0,
     \label{eq:Examples.Levi-form.normal-form}
\end{equation}
so that $\pm i\lambda_{1},\ldots,\pm i\lambda_{n}$ are the nonzero eigenvalues of $L(a)$ counted with multiplicity.

Let $\omega^{1},\ldots,\omega^{n}$ be the coframe of $H^{*}$ dual to $X_{1},\ldots,X_{d}$. For a 1-form $\omega$ we let  $\varepsilon(\omega)$  denote the 
exterior product  
and $\iota(\omega)$ denote the interior product with $\omega$, that is, the contraction with the vector fields dual to $\omega$. 
For an ordered subset $J=\{j_{1},\ldots,j_{k}\}\subset 
\{1,\ldots,d\}$, so that $j_{1}<\ldots<j_{d}$, we let $\omega^{J}=\omega^{j_{1}}\wedge \ldots \wedge \omega^{j_{k}}$ (we make the convention 
that $\omega^{\emptyset}=1$). Then the forms $\omega^{J}$'s give rise to an orthonormal frame of $\Lambda^{*}_{\C}H^{*}$ over $U$. With respect to this frame 
we have
% the operators $d_{b}$ and $d_{b}^{*}$ are of the form, 
\begin{equation}
    d_{b}=\sum_{j=1}^{d}\varepsilon(\omega^{j})X_{j} \qquad \text{and} \qquad d_{b}=-\sum_{l=1}^{d}\iota(\omega^{l})X_{l}+\op{O}_{H}(1).
%     \label{eq:}
\end{equation}
Thus, %the operator  $$ is equal to
\begin{multline}
 \Delta_{b}=d_{b}^{*}d_{b}+d_{b}d_{b}^{*}=-\sum_{j,l}[\varepsilon(\omega^{j})\iota(\omega^{l})X_{j}X_{l}+ 
    \iota(\omega^{l})\varepsilon(\omega^{j})X_{l}X_{j}]+\op{O}_{H}(1) =\\
   -\frac{1}{2}\sum_{j,l=1}^{d} [(\varepsilon(\omega^{j})\iota(\omega^{l})+\iota(\omega^{l})\varepsilon(\omega^{j}))(X_{j}X_{l}+X_{l}X_{j})
   + (\varepsilon(\omega^{j})\iota(\omega^{l})-\iota(\omega^{l})\varepsilon(\omega^{j}))[X_{j},X_{l}] ]+\op{O}_{H}(1).
    \end{multline}
Combining this with~(\ref{eq:Examples.Levi-form.coefficients}) and the relations,
 \begin{equation}
   \varepsilon(\omega^{j})\iota(\omega^{l})+\iota(\omega^{l})\varepsilon(\omega^{j}) 
   =  \delta_{jl}, \qquad j,l=1,\ldots,d,
% %     \label{eq:¥}
 \end{equation}
 we then obtain
\begin{equation}
    \Delta_{b}=-\sum_{j=1}^{d}X_{j}^{2}-i\mu(x)X_{0}+\op{O}_{H}(1), \qquad 
    \mu(x)= \frac{1}{i}\sum_{j,l=1}^{d}\varepsilon(\omega^{j})\iota(\omega^{l})L_{jl}.
%     \label{eq:}
\end{equation}
In particular, thanks to~(\ref{eq:Examples.Levi-form.normal-form}) at $x=a$ we have
\begin{equation}
    \mu(a)=\frac{1}{i}\sum_{j=1}^{n}(\varepsilon(\omega^{j})\iota(\omega^{n+j})-\varepsilon(\omega^{n+j})\iota(\omega^{j}))\lambda_{j}.
%     \label{eq:¥}
\end{equation}

For $j=1,\ldots,n$ let $\theta^{j}=\frac{1}{\sqrt{2}}(\omega^{j}+i\omega^{n+j})$ and 
$\theta^{\bar{j}}=\frac{1}{\sqrt{2}}(\omega^{j}-i\omega^{n+j})$. 
Then we have:
\begin{multline}
\frac{1}{i}(\varepsilon(\omega^{j})\iota(\omega^{n+j})-\varepsilon(\omega^{n+j})\iota(\omega^{j}))= \\
\frac{-1}{2}[(\varepsilon(\theta^{j})+\varepsilon(\theta^{\bar{j}}))(\iota(\theta^{\bar{j}})-\iota(\theta^{j}))- 
(\varepsilon(\theta^{j})-\varepsilon(\theta^{\bar{j}}))(\iota(\theta^{\bar{j}})+\iota(\theta^{j}))]
 =  \varepsilon(\theta^{j})\iota(\theta^{j})-\varepsilon(\theta^{\bar{j}})\iota(\theta^{\bar{j}}).
%  
%  = (\varepsilon(\omega^{j})+i\varepsilon(\omega^{n+j}))(\iota(\omega^{j})-i\iota(\omega^{n+j}))\\
%     = \varepsilon(\omega^{j})\iota(\omega^{j})+\varepsilon(\omega^{n+j})\iota(\omega^{n+j})+
%     i(\varepsilon(\omega^{j})\iota(\omega^{n+j})-\varepsilon(\omega^{n+j})\iota(\omega^{j})).
% %         \label{eq:¥}
\end{multline}
% Similarly, we have
% \begin{equation}
%      2 = \varepsilon(\omega^{j})\iota(\omega^{j})+\varepsilon(\omega^{n+j})\iota(\omega^{n+j})-
%     i(\varepsilon(\omega^{j})\iota(\omega^{n+j})-\varepsilon(\omega^{n+j})\iota(\omega^{j})).
% %     \label{eq:¥}
% \end{equation}
Thus, 
% $\varepsilon(\theta^{j})\iota(\theta^{j})-\varepsilon(\theta^{\bar{j}})\iota(\theta^{\bar{j}})=i(\varepsilon(\omega^{j})\iota(\omega^{n+j})- 
% \varepsilon(\omega^{n+j})\iota(\omega^{j}))$, from which we get
\begin{equation}
   \mu(a)=\sum_{j=1}^{n}(\varepsilon(\theta^{j})\iota(\theta^{j})-\varepsilon(\theta^{\bar{j}})\iota(\theta^{\bar{j}}))\lambda_{j}.
    \label{eq:Examples.mu(a)}
\end{equation}

For any ordered subset $J= \{j_{1},\ldots,j_{p}\}$ of $\{1,\ldots,n\}$ we let 
\begin{equation}
    \theta^{J}=\theta^{j_{1}}\wedge \ldots \wedge \theta^{j_{p}}, \qquad \theta^{\bar{K}}=\theta^{\bar{k}_{1}}\wedge 
    \ldots\wedge\theta^{\bar{k}_{q}},
%     \label{eq:¥}
\end{equation}
Then the forms $\theta^{J}\wedge \theta^{\bar{K}}\wedge \omega^{L}$ give rise to an orthonormal frame of $\Lambda_{\C}^{*}H^{*}$ as $J$ and $K$ range 
over all the ordered subsets of $\{1,\ldots,n\}$ and $L$ over all the ordered subsets of $\{2n+1,\ldots,d\}$. For $j=1,\ldots,n$ we have 
\begin{gather}
    \varepsilon(\theta^{j})\iota(\theta^{j})(\theta^{J}\wedge \theta^{\bar{K}}\wedge \omega^{L} )= 
    \left\{ 
    \begin{array}{cl}
       \theta^{J}\wedge \theta^{\bar{K}}\wedge \omega^{L}   & \text{if $j\in J$}, \\
       0 & \text{if $j\not \in J$},
%         0 & \text{if $J=\emptyset$},
    \end{array}\right. \\
    \varepsilon(\theta^{\bar{j}})\iota(\theta^{\bar{j}})(\theta^{J}\wedge \theta^{\bar{K}}\wedge \omega^{L} )= 
    \left\{ 
    \begin{array}{cl}
       \theta^{J}\wedge \theta^{\bar{K}}\wedge \omega^{L}   & \text{if $j\in K$}, \\
        0 & \text{if $j\not \in K$}.
%         0 & \text{if $K=\emptyset$}.
    \end{array}\right.
%     \label{eq:¥}
\end{gather}
% and
% \begin{equation}
% %     \label{eq:¥}
% \end{equation}
Combining this with~(\ref{eq:Examples.mu(a)}) then gives
\begin{equation}
    \mu(a)(\theta^{J}\wedge \theta^{\bar{K}}\wedge \omega^{L}) = \mu_{J,\bar{K}}(a)\theta^{J}\wedge \theta^{\bar{K}}\wedge \omega^{L}, \quad 
        \mu_{J,\bar{K}}(a)=\sum_{j\in J} \lambda_{j}-\sum_{j\in K}\lambda_{j}.
%     \label{eq:¥}
\end{equation}
This shows that $\mu(a)$ diagonalizes in the basis of $\Lambda^{*}_{\C}H^{*}_{a}$ provided by the forms of $\theta^{J}\wedge \theta^{\bar{K}}\wedge 
\omega^{L}$ with eigenvalues given by the numbers $\mu_{J,\bar{K}}(a)$. In particular, for $k=0,\ldots,d$ we have 
\begin{equation}
    \op{Sp}\mu(a)_{|_{\Lambda^{k}H^{*}}}=\{\mu_{J,\bar{K}}; \ |J|+|K|\leq k\}. 
%     \label{eq:¥}
\end{equation}

Note that we always have $|\mu_{J,K}|\leq \sum_{j=1}^{n}\lambda_{j}$ with equality if, and only if, one the subsets $J$ or $K$ is empty and the other is 
$\{1,\ldots,n\}$, which occurs for eigenvectors in the subspace spanned by the forms $\theta^{1}\wedge \ldots\theta^{n}\wedge \omega^{L}$ and 
$\theta^{\bar{1}}\wedge \ldots\theta^{\bar{n}}\wedge \omega^{L}$  as $L$ ranges over all the subsets of $\{2n+1,\ldots,d\}$. 

Since $\lambda_{1},\ldots,\lambda_{n}$ are the eigenvalues of $|L(a)|$, each of them counted twice, if follows that the condition~(\ref{eq:Sublaplacian.condition}) 
for $\Delta_{b|_{\Lambda^{k}_{\C}H^{*}}}$ reduces to 
$\pm \sum_{j=1}^{n}\lambda_{j}\not \in \op{Sp}\mu(a)_{|_{\Lambda^{k}H^{*}}}$. 
This latter condition is satisfied if, and only if, the space $\Lambda^{k}_{\C}H^{*}_{a}$ does contain any of the forms 
$\theta^{1}\wedge \ldots\theta^{n}\wedge \omega^{L}$ and $\theta^{\bar{1}}\wedge \ldots\theta^{\bar{n}}\wedge \omega^{L}$ with $L$ subset of 
$\{2n+1,\ldots,d\}$. Therefore, the sublaplacian $\Delta_{b|_{\Lambda^{k}_{\C}H^{*}}}$ satisfies~(\ref{eq:Sublaplacian.condition}) at $a$ if, 
and only if, the integer $k$ is not between $n$ and $n+d-2n=d-n$, that is, if, and only if, the condition $X(k)$ holds at $a$. 
The proof is thus achieved. 
\end{proof}

Finally, suppose that $M$ is a CR manifold with Heisenberg structure $H=\Re(T_{1,0}\oplus T_{0,1})$ and assume that 
$T_{\C}M$ is endowed with a Hermitian metric with respect to which $T_{1,0}$ and $T_{0,1}$ are orthogonal subspaces and complex conjugation is an isometry. 
Then  we have  $d_{b}=\dbarb +\partial_{b}$, where $\partial_{b}$ denotes the conjugate of $\dbarb$, that is, $\partial_{b}\alpha =\overline{\dbarb \bar \omega}$ 
for any $\omega\in C^{\infty}(M,\Lambda_{\C}^{*}H^{*})$. Moreover, as $\dbarb \partial_{b}^{*}+\partial_{b}^{*} \dbarb= \dbarb^{*} \partial_{b}+\partial_{b} 
 \dbarb^{*}=0$ (see~\cite{Ta:DGSSPCM}) we get 
 \begin{equation}
     \Delta_{b}=\Box_{b}+ \overline{\Box}_{b},
      \label{eq:Operators.Tanaka-Kohn}
 \end{equation}
 where $\overline{\Box}_{b}$ is the conjugate of $\Box_{b}$. In particular, we see that the horizontal sublaplacian preserves the bidegree, i.e., it 
 acts on $(p,q)$-forms. 
 
  In fact, along the similar lines as that of~\cite[pp.~151--156]{BG:CHM} and of the proof of Proposition~\ref{prop:Examples.horizontal-sublaplacian} 
  above, one can show that the 
  condition~(\ref{eq:Sublaplacian.condition}) at a point $x\in M$ for $\Delta_{b}$ acting on $(p,q)$-forms is equivalent to the condition~$Y(p,q)$ below: 
  \begin{equation}
     \{ (p,q),(q,p)\}\cap \{(\kappa(x)+j,r(x)-\kappa(x)+k);\ \max(j,k)\leq n-r(x)\} =\emptyset,
%        \label{eq:¥}
  \end{equation}
  where $(r(x)-\kappa(x),\kappa(x),n-r(x))$ is the signature at $x$ of the Levi form $L_{\theta}$ associated to some non-vanishing real 1-form $\theta$ 
  anihilating $T_{1,0}\oplus T_{0,1}$. In particular, when $M$ is $\kappa$-strictly pseudoconvex the $Y(p,q)$ reduces to $(p,q)\neq (\kappa,n-\kappa)$ 
  and $(p,q)\neq (n-\kappa,\kappa)$.

% The proof essentially follows similar arguments as those involved in the proofs of Proposition~\ref{thm:PsiDO.Rockland-sublaplacian} 
% and Proposition~\ref{prop:Examples.Boxb} 
% in~\cite{BG:CHM}. We refer to Appendix for a complete proof. 

\subsection{Contact complex and the contact Laplacian}
Given an orientable contact manifold $(M^{2n+1},\theta)$ the contact complex of Rumin~\cite{Ru:FDVC} can be seen as an attempt to get on $M$ 
a complex of horizontal differential forms by forcing out the equalities $d_{b}^{2}=0$ and $(d_{b}^{*})^{2}=0$ as follows.

Let $H=\ker \theta$ and assume that $H$ is endowed with a calibrated almost complex structure $J\in \End_{\R}H$, $J^{2}=-1$, so that 
$d\theta(X,JX)=-d\theta(JX,X)>0$ for any section $X$ of $H$. We then can endow $M$ with the Riemannian metric $g_{\theta}=d\theta(.,J.)+\theta^{2}$.

In addition, let $T$ be the Reeb fields of $\theta$. Then we have 
\begin{equation}
    d_{b}^{2}=-\cL_{T}\varepsilon(d\theta)=\varepsilon(d\theta)\cL_{T},
\end{equation}
where $\varepsilon(d\theta)$ denotes the exterior multiplication by $d\theta$. 

There are two ways of modifying the space $\Lambda^{*}_{\C}H^{*}$ of horizontal forms to get a complex. The first one is to force the equality $d_{b}^{2}=0$ by
restricting the operator $d_{b}$ to $\Lambda^{*}_{2}:=\ker \varepsilon(d\theta) \cap \Lambda^{*}_{\C}H^{*}$ since this 
bundle is stable under $d_{b}$ and on there $d_{b}^{2}$ vanishes. 

The second way is to similarly force the equality $(d_{b}^{*})^{2}=0$ by restricting $d_{b}^{*}$ to  
$\Lambda^{*}_{1}:=\ker 
\iota(d\theta)\cap \Lambda^{*}_{\C}H^{*}=(\im \varepsilon(d\theta))^{\perp}\cap \Lambda^{*}_{\C}H^{*}$, where $\iota(d\theta)$ denotes the interior product 
with $d\theta$. This amounts to replace $d_{b}$ by the operator $\pi_{1}\circ d_{b}$, where $\pi_{1}$ is the orthogonal projection onto $\Lambda^{*}_{1}$.

On the other hand, since $d\theta$ is nondegenerate on $H$ the operator $\varepsilon(d\theta):\Lambda^{k}_{\C}H^{*}\rightarrow \Lambda^{k+2}_{\C}H^{*}$  is 
injective for $k\leq n-1$ and surjective for $k\geq n+1$. This implies that $\Lambda_{2}^{k}=0$ for $k\leq n$ and $\Lambda_{1}^{k}=0$ for $k\geq n+1$. 
Therefore, we only have two halves of complexes. The key observation of Rumin is that we get a full complex by connecting 
the two halves by means of the (second order) differential operator, 
\begin{equation}
    D_{R}:\Lambda_{1}^{n} \longrightarrow\Lambda_{2}^{n+1}, \qquad 
    D_{R}=\cL_{T}-d_{b}\varepsilon(d\theta)^{-1}d_{b},
\end{equation}
where $\varepsilon(d\theta)^{-1}$ is the inverse of $\varepsilon(d\theta):\Lambda^{n-1}_{\C}H^{*}\rightarrow \Lambda^{n+1}_{\C}H^{*}$. 
In other words, letting 
$\Lambda^{k}=\Lambda^{k}_{1}$ for $k=0,\ldots,n-1$ and $\Lambda^{k}=\Lambda^{k}_{2}$ for $k=n+1,\ldots,2n$, we have the  chain complex,
\begin{equation}
    C^{\infty}(M)\stackrel{d_{R}}{\rightarrow}
    \ldots \stackrel{d_{R}}{\rightarrow} C^{\infty}(M,\Lambda^{n}_{1})\stackrel{D_{R}}{\rightarrow} C^{\infty}(M,\Lambda^{n}_{2}) \stackrel{d_{R}}{\rightarrow}
    \ldots \stackrel{d_{R}}{\rightarrow} C^{\infty}(M,\Lambda^{2n}),
     \label{eq:Operators.contact-complex}
\end{equation}
where $d_{R}:C^{\infty}(M,\Lambda^{k})\rightarrow C^{\infty}(M,\Lambda^{k+1})$ is equal to $\pi_{1}\circ d_{b}$ for $k=0,\ldots,n-1$ and to $d_{b}$ for 
$k=n+1,\ldots,2n$. This complex is called the contact complex of $M$. 

The contact Laplacian $\Delta_{R}:C^{\infty}(M,\Lambda^{*}\oplus \Lambda^{n}_{*})\rightarrow C^{\infty}(M,\Lambda^{*}\oplus \Lambda^{n}_{*})$ 
is given by the formulas,
\begin{equation}
    \Delta_{R}=\left\{
    \begin{array}{ll}
        (n-k)d_{R}d_{R}^{*}+(n-k+1) d_{R}^{*}d_{R}& \text{on $\Lambda^{k}$, $k=0,\ldots,n-1$},\\
        (d_{R}^{*}d_{R})^{2}+D_{R}^{*}D_{R} &  \text{on $\Lambda^{n}_{1}$},\\
        D_{R}D_{R}^{*}+  (d_{R}d_{R}^{*})& \text{on $\Lambda^{n}_{2}$},\\
         (n-k+1)d_{R}d_{R}^{*}+(n-k) d_{R}^{*}d_{R} & \text{on $\Lambda^{k}$, $k=n+1,\ldots,2n$}.
    \end{array}\right.
\end{equation}

By comparing the contact Laplacian $\Delta_{R}$ to the horizontal sublaplacian $\Delta_{b}$ Rumin~\cite{Ru:FDVC} was able to show that on every 
degree  $\Delta_{R}$ satisfies the Rockland condition at every point. He then used results 
of Helffer-Nourrigat~(\cite{HN:HMOPCV}) to show that $\Delta_{R}$ was maximal hypoelliptic.  

Alternatively, since in the contact case both $d\theta$ 
and $\cL$ are nondegenerate on $H$, once the Rockland condition is granted 
we may directly apply Proposition~\ref{thm:PsiHDO.Rockland-Parametrix} to get: 

\begin{proposition} 1) For $k=0,\ldots,2n$ with $k\neq n$ the contact Laplacian $\Delta_{R}$ acting on sections of $\Lambda^{k}$ has an invertible 
principal symbol of degree $-2$, hence admits a parametrix in $\pvdo^{-2}(M,\Lambda^{k})$ and is 
hypoelliptic with loss of one derivative.\smallskip

2) For $k=n$ the contact Laplacian $\Delta_{R}$ acting on sections of $\Lambda^{n}_{*}$ has an invertible 
principal symbol of degree $-4$, hence admits a parametrix in $\pvdo^{-4}(M,\Lambda^{n}_{*})$ and 
is hypoelliptic with loss of two derivatives.
\end{proposition}

\section{Proof of Proposition~\ref{prop:PsiHDO.invariance}}
\label{sec.Invariance}
First, we need the lemma below. 

\begin{lemma}\label{lem:Appendix.Heisenberg.invariance}
  For $\Re m>0$ we have $\cK^{m}(\URd)\subset C^{\infty}(U)\hotimes C^{[\frac{\Re m}{2}]}(\Rd)$.  
\end{lemma}
\begin{proof}  
 Let $N=[\frac{\Re m}{2}]$ and let $\alpha$ be a multi-order such that $|\alpha|\leq N$. As $\brak\alpha \leq 2|\alpha|\leq 
-(k+d+2)$ the multiplication by $\xi^{\alpha}$ maps $S^{\hat{m}}(\URd)$ to $C^{\infty}(U)\hotimes 
L^{1}(\Rd)$. Composing it with the inverse Fourier transform with respect to $\xi$ then shows that the map 
$p\rightarrow \partial_{y}^{\alpha}\check{p}_{\xiy}$ maps $S^{\hat{m}}(\URd)$ to $C^{\infty}(U)\hotimes 
C^{0}(\Rd)$. It then follows that for any $p \in S^{\hat{m}}(\URd)$ the transform
$\check{p}_{\xiy}(x,y)$  belongs to $C^{\infty}(U)\hotimes C^{N}(\Rd)$.
     
Now, if $K \in \cK^{m}(\URd)$ then by Lemma~\ref{lem:PsiHDO.characterization.Km} 
     there exists $p\in S^{\hat{m}}(\URd)$, $\hat{m}=-(m+d+2)$, such that $K(x,y)$ is equal to $\check{p}_{\xiy}(x,y)$ modulo 
     a smooth function. Hence $K(x,y)$ belongs to  
     $C^{\infty}(U)\hotimes C^{N}(\Rd)$. The lemma is thus proved.
 \end{proof}
% \begin{proof}  
%      Let $\alpha$ be a multi-order such that $|\alpha|\leq N$. Then we have $\brak\alpha \leq 2|\alpha|\leq 
% -(k+d+2)$, so the multiplication by $\xi^{\alpha}$ maps continuously $S_{||}^{k}(\URd)$ to $C^{\infty}(U)\hotimes 
% L^{1}(\Rd)$. Composing it with the inverse Fourier transform with respect to $\xi$ then shows that the map 
% $p\rightarrow \partial_{y}^{\alpha}\check{p}_{\xiy}$ is continuous from $S_{||}^{k}(\URd)$ to $C^{\infty}(U)\hotimes 
% C^{0}(\Rd)$. Thus the map $p\rightarrow \check{p}_{\xiy}$  is continuous from $S_{||}^{k}(\URd)$ to 
% $C^{\infty}(U)\hotimes C^{N}(\Rd)$.\smallskip  
%      
%      2) Let $K \in \cK^{m}(\URd)$. Then by Lemma~\ref{lem:PsiHDO.characterization.Km} 
%      there exists $p\in S^{\hat{m}}(\URd)$, $\hat{m}=-(m+d+2)$, such that $K(x,y)$ is equal to $\check{p}_{\xiy}(x,y)$ modulo 
%      a smooth function. Since $S^{\hat m}(\URd) \subset S_{||}^{\Re \hat m}(\URd)$ it follows from the first part that $K$ is in 
%      $C^{\infty}(U)\hotimes C^{[\frac{\Re m}{2}]}(\Rd)$. Thus $\cK^{m}(\URd)$ is contained in 
%      $C^{\infty}(U)\hotimes C^{[\frac{\Re m}{2}]}(\Rd)$. 
%  \end{proof}

 We are now ready to prove Proposition~\ref{prop:PsiHDO.invariance}. Let $\tilde{U}$ be an open subset of $\Rd$ together with a hyperplane bundle 
 $\tilde{H}\subset T\tilde{U}$ and a $\tilde{H}$-frame of $T\tilde{U}$ and let $\phi:(U,H)\rightarrow (\tilde{U},\tilde{H})$ be a Heisenberg 
 diffeormorphism. Let $\tilde{P}\in \Psi_{\tilde{H}}^{m}(V)$ and set $P=\phi^{*}\tilde{P}$. We need to show that $P$ is a \psivdo\ of order $m$ on $U$. 

 First, by Proposition~\ref{prop:PsiVDO.characterisation-kernel2}
the distribution kernel of $\tilde{P}$ takes the form,   
 \begin{equation}
      k_{\tilde{P}}(\tilde{x},\tilde{y})= |{\tilde{\varepsilon}_{\tilde{x}}'}| 
      K_{\tilde{P}}(\tilde{x},-{\tilde{\varepsilon}}_{\tilde{x}}(\tilde{y})) + \tilde{R}(\tilde{x},\tilde{y}),  
\label{eq:Appendix.kP}
 \end{equation}
with $K_{\tilde{P}}(\tilde{x},\tilde{y})$ in $\cK^{\hat m}(\tilde{U}\times\Rd)$ and  $\tilde{R}(\tilde{x},\tilde{y})$ in 
$C^\infty(\tilde{U}\times \tilde{U})$. 
Therefore,  the distribution kernel  of $P=\phi^{*}P$ is given by 
\begin{equation}
    k_{P}(x,y) = |\phi'(y)|  k_{\tilde{P}}(\phi(x),\phi(y))= 
     |\varepsilon'_{x}| K(x, -\varepsilon_{x}(y)) +  \tilde{R}(\phi(x),\phi(y)),
\end{equation}
where $K$ is the distribution on $\cU=\{ (x,y)\in \URd; \  \varepsilon_{x}^{-1}(-y)\in U\} \subset \URd$ given by 
 \begin{equation}
     K(x,y)=|\partial_{y}\Phi(x,y)| K_{\tilde{P}}(\phi(x), \Phi(x,y)), 
     \qquad \Phi(x,y)= -\tilde{\varepsilon}_{\phi(x)}\circ \phi\circ \varepsilon_{x}^{-1}(-y). 
 \end{equation}

 Next, it follows from \cite[Props.~3.16, 3.18]{Po:Pacific1} that we have
\begin{equation}
    \Phi(x,y)=-\phi_{H}'(x)(-y)+\Theta(x,y)=\phi_{H}'(x)(y)+\Theta(x,y),
\end{equation}
where $\Theta(x,y)$ is a smooth function on $\cU$ with a behavior near $y=0$ of the form
\begin{equation}
    \Theta(x,y)=(\op{O(\|y\|^{3})},  \op O(\|y\|^{2}),\ldots, \op O(\|y\|^{2})).
     \label{eq:Appendix.invariance.Theta}
\end{equation}
Then a Taylor expansion around $\tilde{y}=\phi_{H}'(x)y$ gives
\begin{gather}
   K(x,y) =  \sum_{\brak\alpha<N}|\partial_{y}\Phi(x,y)| \frac{\Theta(x,y)^\alpha}{\alpha!} 
    (\partial^\alpha_{\tilde{y}}K_{\tilde{P}})(\phi(x),\phi_{H}'(x)y) + R_{N}(x,y), 
   \label{eq:PsivDO.invariance.TaylorK1}\\
     R_{N}(x,y) =  \sum_{\brak\alpha=N} 
   |\partial_{y}\Phi(x,y) | \frac{\Theta(x,y)^\alpha}{\alpha!} \int_{0}^1 (t-1)^{N-1}
    \partial^\alpha_{\tilde{y}}K_{\tilde{P}}(\phi(x),\Phi_{t}(x,y) ) dt.
    \label{eq:PsivDO.invariance.TaylorK2} 
\end{gather}
where we have let $\Phi_{t}(x,y) =\phi_{H}'(x)y +t\Theta(x,y)$.

Set $f_{\alpha}(x,y)=|\partial_{y}\Phi(x,y)| \Theta(x,y)^\alpha$. Then~(\ref{eq:Appendix.invariance.Theta})  implies 
that near $y=0$ we have 
\begin{equation}
    f_{\alpha}(x,y)=\op{O}(\|y\|^{3\alpha_{0}+2(\alpha_{1}+\ldots+\alpha_{d})})=\op{O}(\|y\|^{\frac32\brak\alpha}).
     \label{eq:Appendix.behavior-f-alpha}
\end{equation}
Thus all the homogeneous components of degree $< \frac32\brak\alpha$ in  
the Taylor expansion for $f_{\alpha}(x,y)$ at $y=0$ must be zero. Therefore, we can write
\begin{equation}
    f_{\alpha}(x,y)= \sum_{\frac32\brak\alpha \leq \brak\beta < \frac{3}{2}N} \!\! f_{\alpha\beta}(x) \frac{y^\beta }{\beta !}+ 
    \sum_{\brak\beta\dot{=}\frac{3}{2}N}r_{N\alpha\beta}(x,y)y^{\beta},
     \label{eq:Appendix.behavior-Theta-alpha}
\end{equation}
where we have let $f_{\alpha\beta}(x)=\partial_{y}^{\beta}f_{\alpha}(x,0)$, the functions $r_{M\alpha\beta}(x,y)$ are in $C^{\infty}(\cU)$ 
and the notation $\brak 
\beta\dot{=}\frac{3}{2}N$ means that $\brak 
\beta$ is equal to $\frac{3}{2}N$ if $\frac{3}{2}N$ is an integer and to $\frac{3}{2}N+\frac{1}{2}$ otherwise. Thus,
\begin{equation}
    K(x,y)= \sum_{\brak\alpha<N} \sum_{\frac32\brak\alpha \leq \brak\beta < \frac{3}{2}N}\!\! K_{\alpha\beta}(x,y) + 
    \sum_{\brak\alpha<N} R_{N\alpha}(x,y) + R_{N}(x,y),
     \label{eq:Appendix.expansion-K}
\end{equation}
where we have let
\begin{gather}
     K_{\alpha\beta}(x,y)= f_{\alpha\beta}(x) y^\beta 
    (\partial^\alpha_{\tilde{y}}K_{\tilde{P}})(\phi(x),\phi_{H}'(x)y),\\ R_{N\alpha }(x,y) = 
    \sum_{\brak\beta\dot{=}\frac{3}{2}N} r_{M\alpha\beta}(x,y)y^\beta (\partial^\alpha_{\tilde{y}}K_{\tilde{P}})(\phi(x),\phi_{H}'(x)y). 
\end{gather}

As in the proof of Proposition~\ref{prop:PsiVDO.characterisation-kernel2} 
the smoothness of $\phi_{H}'(x)y$ and the fact that $\phi_{H}'(x)(\lambda.y)=\lambda.\phi_{H'}(x)y$ 
for any $\lambda 
\in \R$  imply that $K_{\alpha\beta}(x,y)$ belongs to $\cK^{\hat m -\brak\alpha+\brak\beta}(U\times\Rd)$. 
Notice that if  $\frac32\brak\alpha \leq \brak\beta \dot{=}\frac{3}{2}N$ then 
$\Re \hat m-\brak\alpha+\brak\beta\geq \Re\hat m +\frac{1}{3}\brak\beta\geq \Re\hat m + \frac{N}{2}$. 
It thus follows from Lemma~\ref{lem:Appendix.Heisenberg.invariance} that, for any integer $J$, the remainder term $R_{N\alpha}$ is 
in $C^{J}(\URd)$ as soon as $N$ is large enough.  

Let $\pi_{x}:\URd \rightarrow U$ denote the projection on the first coordinate. In the sequel we will say that a distribution $K(x,y)\in \cD'(\URd)$ is properly 
supported with respect to $x$ when $\pi_{x|_{\supp K}}$ is a proper map, i.e.~ for any compact $L \subset U$ the set $\supp K \cap (L\times \Rd)$ is compact.

In order to deal with the regularity of $R_{N}(x,y)$ in~(\ref{eq:PsivDO.invariance.TaylorK2}) 
we need the lemma below. 
\begin{lemma}\label{lem:Appendix.Theta-alpha}
  There exists a  function $\chi_{n}\in  C^{\infty}_{c}(\cU)$ properly supported with respect to $x$ such that   
  $\chi(x,y)=1$ near $y=0$ and, for any multi-order $\alpha$, we can write   
 \begin{equation}
      \chi(x,y)\Theta(x,y)^{\alpha}=\sum_{\brak \beta\dot{=}\frac{3}{2}\brak \alpha} \theta_{\alpha\beta}(t,x,y) \Phi_{t}(x,y)^{\beta} 
      \label{eq:Appendix.claim}
\end{equation}
where the functions $\theta_{\alpha\beta}(t,x,y)$ are in $C^{\infty}([0,1]\times U\times \Rd)$.
\end{lemma}
\begin{proof}[Proof of the lemma]
   Let $U'$ be a relatively compact open subset of $U$ and let $(t_{0},x_{0})\in [0,1]\times U'$. 
   Since $\Phi_{t_{0}}(x_{0},0)=0$ and $\partial_{y}\Phi_{t_{0}}(x_{0},0)= \phi'_{H}(x_{0})$ is invertible
  the implicit function theorem implies that there 
  exist an open interval $I_{x_{0}}$ containing $t_{x_{0}}$, an open subset $U_{x_{0}}$ of $U$ containing $x_{0}$, open subsets 
   $V_{x_{0}}$ 
   and $\tilde{V}_{x_{0}}$ of $\Rd$ containing $0$ and a smooth map $\Psi_{x_{0}}(t,x,\tilde{y})$ from 
   $I_{x_{0}}\times U_{x_{0}}\times \tilde{V}_{x_{0}}$ to 
   $V_{x_{0}}$ such that 
    $U_{x_{0}}\times V_{x_{0}}$ is contained in $\cU$ and  for any $(t,x,y)$ in $I_{x_{0}}\times U_{x_{0}}\times V_{x_{0}}$ and any 
    $\tilde{y}$ in $\tilde{V}_{x_{0}}$ we have 
    \begin{equation}
        \tilde{y}=\Phi_{t}(x,y) \Longleftrightarrow y=\Psi_{x_{0}}(t,x,\tilde{y}). 
         \label{eq:Appendix.definition-Psi}
    \end{equation}
   
    Since $[0,1]\times \overline{U'}$ is compact we can cover it by finitely many products $I_{x_{k}}\times U_{x_{k}}$, $k=1,..,p$, with 
    $(t_{k},x_{k})\in [0,1]\times \overline{U'}$. In particular, the sets $I=\cup_{k}I_{k}$ and $U''=\cup U_{k}$ are open neighborhoods of $I$ and 
    $\overline{U'}$ respectively. Thanks to~(\ref{eq:Appendix.definition-Psi}) 
    we have $\Psi_{x_{k}}=\Psi_{x_{l}}$ on 
    $(I_{x_{k}}\times U_{x_{k}}\times V_{x_{k}})\cap (I_{x_{l}}\times U_{x_{l}}\times V_{x_{l}})$. Therefore, 
    setting $V= \cap_{k}V_{k}$ and $\tilde{V}=\cap_{k}\tilde{V}_{k}$ we have $U''\times V\subset \cU$ and there exists a smooth map $\Psi$ from 
    $I\times U''\times \tilde{V}$ such that for any $(t,x,y)$ in $I\times U''\times V$ and any 
    $\tilde{y}$ in $\tilde{V}$ we have 
    \begin{equation}
        \tilde{y}=\Phi_{t}(x,y) \Longleftrightarrow y=\Psi(t,x,\tilde{y}). 
    \end{equation}
    
   Furthermore, as $\partial_{\tilde{y}} \Psi(t,x,0)=[\partial_{y}\Phi_{t}(x,0)]^{-1}=\phi_{H}'(x)^{-1}$ and 
   $\phi_{H}'(x)^{-1}(\lambda.y)=\lambda \phi_{H}'(x)^{-1}(y)$ for any $\lambda \in \R$, the function $\Theta(x,\Psi(t, x,\tilde{y}))$ 
   behaves near 
   $\tilde{y}=0$ as in~(\ref{eq:Appendix.invariance.Theta}). Therefore, as 
   in~(\ref{eq:Appendix.behavior-f-alpha})--(\ref{eq:Appendix.behavior-Theta-alpha})  for any multi-order 
   $\alpha$ we can write 
   \begin{equation}
       \Theta(x,\Psi(t, x,\tilde{y}))^{\alpha}= \sum_{\brak\beta \dot{=} \frac32 \brak\alpha} 
          \tilde{\theta}_{\alpha\beta}(t,x,\tilde{y}) \tilde{y}^{\beta}, \quad \tilde{\theta}_{\alpha\beta}(t,x,\tilde{y})\in C^{\infty}(I\times 
   U''\times \tilde{V}).
   \end{equation}
 Setting $\tilde{y}=\Phi_{t}(x,y)$ then gives 
  \begin{equation}
       \Theta(x,y)^{\alpha}=\sum_{\brak  \beta\dot{=}\frac{3}{2}\brak\alpha} \theta_{\alpha\beta}(t,x,y) 
       \Phi_{t}(x,y)^{\beta}, \quad \theta_{\alpha\beta}(t,x,y)\in C^{\infty}(I\times U''\times V).
       \label{eq:Appendix.lemma.form-Theta-alpha}
  \end{equation}
   
   All this allows us to construct locally finite coverings $(U'_{n})_{n\geq 0}$ and $(U''_{n})_{n\geq 0}$ of $U$ by relatively compact open 
   subsets in such way that, for each integer $n$, the open $U''_{n}$ contains $\overline{U'_{n}}$ and there exists an open $V_{n}\subset \Rd$ containing $0$
so that, for any multiorder $\alpha$, on $[0,1]\times U''_{n}\times V_{n}$ we have   
     \begin{equation}
       \Theta(x,y)^{\alpha}= \sum_{\brak \beta\dot{=}\frac{3}{2}\brak\alpha} \theta_{\alpha\beta}^{(n)}(t,x,y) \Phi_{t}(x,y)^{\beta}, \qquad 
       \theta_{\alpha\beta}^{(n)}(t,x,y)\in C^{\infty}([0,1]\times U''_{n}\times V_{n}).
      \label{eq:Appendix.lemma.form-Theta-alpha-n}
   \end{equation}

   For each $n$ let $\varphi_{n}\in C^{\infty}_{c}(U_{n}'')$ be such that $\varphi_{n}=1$ on $U_{n}'$ and  let 
  $\psi_{n}\in C^{\infty}_{c}(V_{n})$ be such that $\psi_{n}=1$ on a neighborhood $V'_{n}$ of $0$. Then we   
  construct a locally  finite family $(\chi_{n})_{n\geq 0} \subset C^{\infty}_{c}(\cU)$ as follows: for $n=0$ we set  
  $\chi_{0}(x,y)=\varphi_{0}(x)\psi_{0}(y)$ and for $n\geq 1$ we let
\begin{equation}
   \chi_{n}(x,y)= (1-\varphi_{0}(x)\psi_{0}(y))\ldots (1-\varphi_{n-1}(x)\psi_{n-1}(y))\varphi_{n}(x)\psi_{n}(y) . 
\end{equation}
  Then $\chi =  \sum \chi_{n}$ is  a well defined smooth function on $\URd$ 
  supported on $\cup_{n\geq 0} (U_{n}''\times V_{n}) \subset \cU$, hence properly supported with respect to $x$. Also, as $\chi(x,y)=1$ on each product 
  $U_{n}'\times V_{n}'$ we see that $\chi(x,y)=1$ on a neighborhood of $U \times \{0\}$. In addition, 
  thanks to~(\ref{eq:Appendix.lemma.form-Theta-alpha-n}) 
  on $[0,1]\times\cU$ we have 
  \begin{equation}
      \chi(x,y)\Theta(x,y)^{\alpha}=\sum_{n\geq 0} \chi_{n}(x,y)\Theta(x,y)^{\alpha}=  \sum_{\brak \beta\dot{=}\frac{3}{2}\brak\alpha}  
      \theta_{\alpha\beta}(t,x,y) \Phi_{t}(x,y)^{\beta},
  \end{equation}
where the functions $\theta_{\alpha\beta}(t,x,y) :=\sum_{n} \chi_{n}(x,y) \theta_{\alpha\beta}^{(n)}(t,x,y)$ are in 
$C^{\infty}([0,1]\times\URd)$. The lemma is thus proved.
 \end{proof}

Let us go back  to the proof of Proposition~\ref{prop:PsiHDO.invariance}. Thanks to~(\ref{eq:PsivDO.invariance.TaylorK2}) 
and~(\ref{eq:Appendix.claim})  on $U\times \R$ we have 
\begin{equation}
    \chi(x,y) R_{N}(x,y)=\sum_{\brak\alpha=N}\sum_{\brak\beta\dot{=}\frac{3}{2}N} \int_{0}^{1} 
    r_{\alpha\beta}(t,x,y)(\tilde{y}^{\beta}\partial^\alpha_{\tilde{y}}K_{\tilde{P}})(\phi(x),\Phi_{t}(x,y))dt, 
\end{equation}
for some functions $r_{\alpha\beta}(t,x,y)$ in $C^{\infty}([0,1\times U\times \Rd)$. Since  
$\tilde{y}^{\beta}\partial^\alpha_{\tilde{y}}K_{\tilde{P}}$ belongs to $\cK^{\hat{m}+N/2}(\tilde{U}\times \Rd)$ it follows from 
Lemma~\ref{lem:Appendix.Heisenberg.invariance} that, for any integer $J\geq 0$, as soon as 
$N$ is taken large enough $\tilde{y}^{\beta}\partial^\alpha_{\tilde{y}}K_{\tilde{P}}$ is in $C^{J}(\tilde{U}\times \Rd)$ and so 
$\chi(x,y) R_{N}(x,y)$ is in $C^{J}(\URd)$.

On the other hand, set $K_{P}(x,y)=\chi(x,y)K(x,y)=\sum \chi_{n}(x,y)K(x,y)$. Since $\chi(x,y)$ is supported in $\cU$ and is properly supported with 
respect to $x$ this defines a distribution on $U\times \Rd$. Moreover, using~(\ref{eq:Appendix.expansion-K}) 
we get
\begin{equation}
    K_{P}(x,y)= \sum_{\brak \alpha<N} \sum_{\frac32\brak\alpha \leq \brak\beta 
    < \frac{3}{2}N} K_{\alpha\beta}(x,y)+\sum_{j=1}^{3}R_{N}^{(j)},
        \label{eq:Appendix.expansion-KP1}
 \end{equation}   
 where the remainder terms $R_{N,z}^{(j)}$, $j=1,2,3$ are given by 
 \begin{gather}
R_{N}^{(1)}=\chi(x,y)R_{N}(x,y),\quad  R_{N}^{(2)}= \sum_{\brak \alpha <N} \chi(x,y) R_{N\alpha}(x,y),\\
   R_{N}^{(2)}(x,y)=  \sum_{\brak \alpha<N} \sum_{\frac32\brak\alpha \leq \brak\beta 
    <\frac{3}{2}N}  (1-\chi(x,y))K_{\alpha\beta}(x,y).
\end{gather}
Each term $K_{\alpha\beta}(x,y)$ belongs to $\cK^{\hat 
m-\brak\alpha+\brak \beta}(\URd)$ and, as $\hat{m}+\brak\beta-\brak\alpha=\hat{m}+j$ 
and $\frac32\brak\alpha \leq \brak\beta$ imply $\brak\alpha\leq 2j$ and $\brak\beta \leq \frac{4}{3}j$, 
in the r.h.s.~(\ref{eq:Appendix.expansion-KP1}) there are only  finitely 
many such distributions in a given space $\cK^{\hat m+j}(\URd)$ as $\alpha$ and $\beta$ range over all multi-orders such that 
$\frac32\brak\alpha \leq \brak\beta$.    

Furthermore, the reminder term $R_{N}^{(3)}$ is smooth and the other remainder terms $R_{N}^{(j)}$, $j=1,2$, are in $C^{J}(\URd)$ as soon as $N$ is 
large enough.  Thus, 
\begin{equation}
    K_{P}(x,y)  \sim\sum_{\frac32\brak\alpha \leq \brak\beta} K_{\alpha\beta}(x,y), 
    \label{eq:Appendix.asymptotic-expansion-KP'}
\end{equation}
which implies that $K_{P}$ belongs to $\cK^{\hat m}(\URd)$ and satisfies~(\ref{eq:Appendix.asymptotic-expansion-KP'}).

Finally, from~(\ref{eq:Appendix.kP}) and the very definition of $\Phi(x,y)$ on $U\times U$, we deduce that the distribution kernel of $P$ 
is equal to
\begin{multline}
      |\varepsilon_{x}'| K_{P}(x,-\varepsilon_{x}(y)) + [1-\chi(x,\varepsilon_{x}(y))]  
       |\tilde{\varepsilon}'_{\phi(x)}| K_{\tilde{P}}(\phi(x),-\tilde{\varepsilon}_{\phi(x)}\circ\phi(y))+ \tilde{R}(\phi(x),\phi(y))\\
       =  |\varepsilon_{x}'| K_{P}(x,-\varepsilon_{x}(y)) \quad \bmod C^{\infty}(U\times U).
\end{multline}
Combining this with Proposition~\ref{prop:PsiVDO.characterisation-kernel2} and the fact that $K_{P}(x,y)$ 
satisfies~(\ref{eq:Appendix.asymptotic-expansion-KP'}) completes the proof of Proposition~\ref{prop:PsiHDO.invariance}.

\section{Proof of Proposition~\ref{prop:PsiHDO.transpose-chart}}\label{sec:transpose}
Let $P:C^{\infty}_{c}(U)\rightarrow C^{\infty}(U)$ be a  \psivdo\ of order $m$ and let us show that its transpose operator 
$P^{t}:C^{\infty}_{c}(U)\rightarrow C^{\infty}(U)$ is a \psivdo\ of order $m$. 
By Proposition~\ref{prop:PsiVDO.characterisation-kernel2} the distribution kernel of $P$ is of the form, 
\begin{equation}
    k_{P}(x,y)=|\varepsilon_{x}'|K_{P}(x,-\varepsilon_{x}(y)) +R(x,y),
\end{equation}
with $K_{P}(x,y)$ in $\cK^{m}(\URd)$ and $R(x,y)$ in $C^{\infty}(U\times U)$. Thus the distribution kernel of $P^{t}$ can be written as 
\begin{equation}
    k_{P^{t}}(x,y)=k_{P}(y,x)=|\varepsilon_{y}'|K_{P}(y,-\varepsilon_{y}(x)) +R(y,x)= |\varepsilon_{x}'|K(x,-\varepsilon_{x}(y)) +R(y,x),
     \label{eq:Appendix.kPt-K}
\end{equation}
where $K$ is the distribution on the open $\cU=\{(x,y); \ \varepsilon_{x}^{-1}(-y)\in U\}$ given by
\begin{equation}
    K(x,y)=|\varepsilon_{x}'|^{-1}|\varepsilon_{y}'|K_{P}(\varepsilon_{x}^{-1}(-y), -\varepsilon_{\varepsilon_{x}^{-1}(-y)}(x)).
\end{equation}

\begin{lemma}
    On $\cU$ we have
    \begin{equation}
        \varepsilon_{\varepsilon_{x}^{-1}(-y)}(x)=y-\Theta(x,y),
    \end{equation}
 where $\Theta: \cU \rightarrow \Rd$ is a smooth map with a behavior near $y=0$ of the form~(\ref{eq:Appendix.invariance.Theta}). 
\end{lemma}
\begin{proof}
    Let $(x,y)\in \cU$ and $Y\in G_{x}U$ let $\lambda_{y}(Y)=y.Y$, that is,  $\lambda_{y}$ is the left multiplication by $y$ on $G_{x}U$. 
    Then by~\cite[Eq.~(3.32)]{Po:Pacific1} for $Y$ small enough we have
    \begin{equation}
        \lim_{t\rightarrow 0} \varepsilon_{x}\circ \varepsilon_{\varepsilon_{x}^{-1}(t.-y)}^{-1}(t.Y)=\lambda_{-y}(Y)=\lambda_{y}^{-1}(Y).
    \end{equation}
    Since $\varepsilon_{x}\circ \varepsilon_{\varepsilon_{x}^{-1}(t.-y)}^{-1}(t.Y)$ is a smooth function of $(t,Y)$ near $(0,0)$, it 
    follows from the implicit function theorem that for $Y$ small enough we have 
    \begin{equation}
        \lim_{t \rightarrow 0} t^{-1}.\varepsilon_{\varepsilon_{x}^{-1}(t.-y)} \circ \varepsilon_{x}^{-1}(Y)=\lambda_{y}(Y).
    \end{equation}
    In particular, for $Y=0$ we get 
    \begin{equation}
         \lim_{t \rightarrow 0} t^{-1}.\varepsilon_{\varepsilon_{x}^{-1}(t.-y)}(x)=y.
         \label{eq:Appendix.inverse-asymptotic}
    \end{equation}
    
    Now, the function $ \varepsilon_{\varepsilon_{x}^{-1}(-y)}(x)$ depends smoothly on $(x,y)\in \cU$, so~(\ref{eq:Appendix.inverse-asymptotic}) 
    allows us to put it into the form,
    \begin{equation}
        \varepsilon_{\varepsilon_{x}^{-1}(-y)}(x)=y-\Theta(x,y),
    \end{equation}
 where $\Theta=(\Theta_{0},\ldots,\Theta_{d})$ is smooth map from $\cU$ to $\Rd$ with a behavior near $y=0$ of the form
 \begin{equation}
     \Theta_{0}(x,y)=\op{O}(|y_{0}|^{2}+|y_{0}||y|+|y|^{3}), \qquad \Theta_{j}(x,y)=\op{O}(|y|^{2}), \quad j=1,\ldots,d. 
 \end{equation}
 In particular, near $y=0$ the map $\Theta$ has a behavior of the form~(\ref{eq:Appendix.invariance.Theta}). %The lemma is thus proved. 
\end{proof}

Next, a Taylor expansion around $(\varepsilon_{x}^{-1}(-y),-y)$ gives
\begin{gather}
    K(x,y)=\sum_{\brak\alpha<N} |\varepsilon_{x}'|^{-1}|\varepsilon_{y}'|\frac{\theta(x,y)^{\alpha}}{\alpha!} 
    (\partial_{y}^{\alpha}K_{P})(\varepsilon_{x}^{-1}(-y),-y) +R_{N}(x,y),\\
    R_{N}(x,y)= \sum_{\brak\alpha=N}|\varepsilon_{x}'|^{-1}|\varepsilon_{y}'|\frac{\theta(x,y)^{\alpha}}{\alpha!} \int_{0}^{1}(1-t)^{N-1} 
    \partial_{y}^{\alpha}K_{P})(\varepsilon_{x}^{-1}(-y),\Phi_{t}(x,y)),
    \label{eq:AppendixB.RN}
\end{gather}
where we have let $\Phi_{t}(x,y)=-y+t\Theta(x,y)$. 

Let $a_{\alpha}(x,y)= 
|\varepsilon_{x}'|^{-1}|\varepsilon_{y}'|\frac{\theta(x,y)^{\alpha}}{\alpha!}$. Then because of~(\ref{eq:Appendix.invariance.Theta}) 
the same arguments used to prove~(\ref{eq:Appendix.behavior-Theta-alpha}) show 
that there exist functions $r_{N\alpha}(x,y)\in C^{\infty}(\cU)$, $\brak \beta \dot{=}\frac{3}{2}N$, such that 
\begin{equation}
    a_{\alpha}(x,y)=\sum_{\frac{3}{2}\brak\alpha \leq \brak \beta < \frac{3}{2}N}a_{\alpha\beta}(x)y^{\beta} + 
    \sum_{\brak \beta \dot{=}\frac{3}{2}N}r_{N\alpha}(x,y)y^{\beta},
\end{equation}
where we have let $a_{\alpha \beta}(x)=\frac{1}{\beta!}\partial^{\beta}a_{\alpha}(x,0)$. Thus,
\begin{gather}
    K(x,y)= \sum_{\brak\alpha<N} \sum_{\frac{3}{2}\brak\alpha \leq \brak \beta < \frac{3}{2}N}   f_{\alpha\beta}(x)y^{\beta} 
    (\partial_{y}^{\alpha}K_{P})(\varepsilon_{x}^{-1}(-y),-y)  + 
    \sum_{\brak\alpha<N}R_{N\alpha}(x,y) +R_{N}(x,y),\\
    R_{N\alpha}(x,y)= \sum_{\brak \beta \dot{=}\frac{3}{2}N}r_{N\alpha}(x,y)y^{\beta} (\partial_{y}^{\alpha}K_{P})(\varepsilon_{x}^{-1}(-y),-y).
    \label{eq:AppendixB.RNalpha}
\end{gather}

Next, a further Taylor  expansion shows that $  (\partial_{y}^{\alpha}K_{P})(\varepsilon_{x}^{-1}(-y),-y)$ is equal to
\begin{equation}
\sum_{|\gamma|<N}\frac{1}{\gamma!}(\varepsilon_{x}^{-1}(-y)-x)^{\gamma} 
  (\partial^{\gamma}_{x}\partial_{y}^{\alpha}K_{P})(x,-y) + \sum_{|\gamma|=N} \int_{0}^{1}(1-t)^{N-1} 
  (\partial^{\gamma}_{x}\partial_{y}^{\alpha}K_{P})(\varepsilon_{t}(x,y),-y),
\end{equation}
where we have let $\varepsilon_{t}(x,y)=x+t(\varepsilon_{x}^{-1}(-y)-x)$. Since $\varepsilon_{x}^{-1}(-y)-x$ is polynomial in $y$ of degree 2 and 
vanishes for $y=0$, we can write
\begin{equation}
    \frac{1}{\gamma!}(\varepsilon_{x}^{-1}(-y)-x)^{\gamma} =\sum_{|\gamma|\leq |\delta \leq 2|\gamma|} a_{\gamma\delta}(x)y^{\delta}, \qquad 
    b_{\gamma\delta}(x)=\frac{1}{\gamma!\delta!}[\partial_{y}(\varepsilon_{x}^{-1}(-y)-x)^{\gamma}] _{y=0}.
\end{equation}
Therefore, we can put $K(x,y)$ into the form, 
\begin{equation}
    K(x,y)= \sum_{\alpha,\beta,\gamma,\delta}^{(N)}
    K_{\alpha\beta\gamma\delta}(x,y) + \sum_{\brak\alpha<N} \sum_{\frac{3}{2}\brak\alpha \leq \brak \beta < \frac{3}{2}N} R_{N\alpha\beta}(x,y)+
    \sum_{\brak\alpha<N} R_{N\alpha}(x,y)  +R_{N}(x,y),
\end{equation}
where the first summation goes over all the multi-orders $\alpha$, $\beta$, $\gamma$ and $\delta$ such that $\brak\alpha<N$, $\frac{3}{2}\brak\alpha \leq 
\brak \beta < \frac{3}{2}N$ and $|\gamma|\leq |\delta| \leq 2|\gamma|<2N$, and we have let
\begin{gather}
     K_{\alpha\beta\gamma\delta}(x,y)=a_{\alpha\beta\gamma\delta}(x) y^{\beta+\delta}  (\partial^{\gamma}_{x}\partial_{y}^{\alpha}K_{P})(x,-y), 
  \qquad  f_{\alpha\beta\gamma\delta}(x)=a_{\alpha\beta}(x) b_{\gamma\delta}(x),\\
    R_{N\alpha \beta}(x,y)= \sum_{|\gamma|=N}\sum_{N\leq |\delta|\leq 2N} a_{\alpha\beta\gamma\delta}(x)y^{\beta+\delta} \int_{0}^{1}(1-t)^{N-1} 
  (\partial^{\gamma}_{x}\partial_{y}^{\alpha}K_{P})(\varepsilon_{t}(x,y),-y).
  \label{eq:AppendixB.RNalphabeta}
\end{gather}

Now, the distribution $y^{\beta}K_{P}(x,-y)$ belongs to $\cK^{\hat m-\brak \alpha+\brak \beta}(\URd)$. In particular, 
if $\frac{3}{2}\brak \alpha \leq \brak \beta\dot{=}\frac{3}{2}N$ then $\Re \hat{m}-\brak\alpha+\brak \beta\geq \Re \hat{m}+\frac{1}{3}\brak \beta \geq  \Re 
\hat{m}+\frac{1}{2}N$. Therefore, for any given integer $J$  Lemma~\ref{lem:Appendix.Heisenberg.invariance} tells us that
$y^{\beta}K_{P}(x,-y)$ is in $C^{J}(\URd)$ as soon as $N$ is large enough. 
It follows that all the remainder terms 
$R_{N\alpha}(x,y)$, $\brak \alpha<N$, belong to $C^{J}(\cU)$ for $N$ large enough. 

Similarly, if $\frac{3}{2}\brak \alpha \leq \brak \beta $ and $|\gamma|=N\leq |\delta|\leq 2N$ then $\Re \hat{m}-\brak\alpha +\brak \beta+\brak 
\delta\geq \Re \hat{m}+\brak{\delta}\geq \Re \hat{m}+\frac{1}{2}N$, so using Lemma~\ref{lem:Appendix.Heisenberg.invariance} 
we see that 
$y^{\beta+\delta} (\partial^{\gamma}_{x}\partial_{y}^{\alpha}K_{P})(x,-y)$ is in $C^{J}(\URd)$ for $N$ large enough. It then follows that for $N$ 
large enough the remainder terms $R_{N\alpha\beta}(x,y)$  with $\brak \alpha<N$ and $\frac{3}{2}\brak \alpha \leq \brak \beta \dot{=}\frac{3}{2}N$ are all 
in $C^{J}(\cU)$ as soon as $N$ is chosen large enough. 

In order to deal with the last remainder term $R_{N}(x,y)$ notice that, along the same lines as that of the proof of 
Lemma~\ref{lem:Appendix.Theta-alpha},  
one can show that there exists a $\chi \in C^{\infty}(\cU)$ properly supported with respect to $x$ such that 
  $\chi(x,y)=1$ near $y=0$ and, for any multi-order $\alpha$, we can write   
 \begin{equation}
      \chi(x,y)\Theta(x,y)^{\alpha}=\sum_{\brak \beta\dot{=}\frac{3}{2}\brak \alpha} \theta_{\alpha\beta}(t,x,y) \Phi_{t}(x,y)^{\beta},
\end{equation}
where the functions $\theta_{\alpha\beta}(t,x,y)$ are in $C^{\infty}([0,1]\times U\times \Rd)$. Then we can put $ \chi(x,y) R_{N}(x,y)$ into the form,
\begin{equation}
    \chi(x,y) R_{N}(x,y)= \sum_{\brak\alpha=N}\sum_{\brak \beta\dot{=}\frac{3}{2}\brak \alpha} 
    |\varepsilon_{x}'|^{-1}|\varepsilon_{y}'| \int_{0}^{1}r_{N\alpha\beta}(t,x,y)
    (y^{\beta}\partial_{y}^{\alpha}K_{P})(\varepsilon_{x}^{-1}(-y),\Phi_{t}(x,y)),
\end{equation}
for some functions $r_{N\alpha\beta}(t,x,y)$ in $C^{\infty}([0,1]\times \URd)$. As $(y^{\beta}\partial_{y}^{\alpha}K_{P})$ is in  
$\cK^{\hat{m}-\brak \alpha+\brak\beta}(\URd)$ and we have $\Re \hat{m}-\brak\alpha+\brak\beta \geq \Re \hat{m}+\frac{1}{2}\brak\alpha=\Re 
\hat{m}+\frac{1}{2}N$, using Lemma~\ref{lem:Appendix.Heisenberg.invariance} we see 
that for $N$ large enough $ \chi(x,y) R_{N}(x,y)$ is in $C^{J}(\cU)$, so is in $C^{J}(\URd)$ since $ \chi(x,y) R_{N}(x,y)$ is a properly supported  with respect to 
$x$.  

Let $K_{P^{t}}(x,y)=\chi(x,y)K(x,y)$. This defines a distribution on $\URd$ since $\chi$ is properly supported with respect to $x$. Moreover, we have
\begin{equation}
  K_{P^{t}}(x,y)= \sum_{\alpha,\beta,\gamma,\delta}^{(N)}
    K_{\alpha\beta\gamma\delta}(x,y) + \sum_{j=1}^{4}R_{N}^{(j)}(x,y), 
\end{equation}
where the remainder terms $R_{N}^{(j)}$, $j=1,\ldots,4$, are given by
\begin{gather}
    R_{N}^{(1)}=\chi(x,y)R_{N}(x,y), \quad 
    R_{N}^{(2)}= \sum_{\brak\alpha<N}\chi(x,y)R_{N\alpha}(x,y),\\
    R_{N}^{(3)}= \sum_{\brak\alpha<N} \sum_{\frac{3}{2}\brak\alpha \leq \brak \beta < \frac{3}{2}N}\chi(x,y)R_{N\alpha\beta}(x,y),\quad 
    R_{N}^{(4)}= \sum_{\alpha,\beta,\gamma,\delta}^{(N)}(1-\chi(x,y))K_{\alpha\beta\gamma\delta}(x,y).
\end{gather}

Note that $K_{\alpha\beta\gamma\delta}(x,y)$ belongs to $\cK^{\hat{m}-\brak\alpha+\brak\gamma+\brak\delta}(\URd)$ and there are finitely many terms 
of a given order as $\alpha$, $\beta$, $\gamma$ and $\delta$ ranges over all the multi-orders such that $\frac{3}{2}\brak\alpha \leq 
\brak \beta$ and $|\gamma|\leq |\delta| \leq 2|\gamma|$. 

On the other hand, the remainder term $R_{N}^{(4)}$ is smooth and it follows from th observations above that the other remainder terms are in $C^{J}(\URd)$ as soon  as 
$N$ is large enough. Thus,
\begin{equation}
    K_{P^{t}}(x,y) \sim \sum_{\frac{3}{2}\brak\alpha \leq \brak \beta} \sum_{|\gamma|\leq |\delta| \leq 2|\gamma|}K_{\alpha\beta\gamma\delta}(x,y),
     \label{eq:Appendix.expansionKPt}
\end{equation}
which incidentally shows that $K_{P^{t}}(x,y)$ belongs to $\cK^{\hat{m}}(\URd)$. 

Finally, thanks to~(\ref{eq:Appendix.kPt-K}) we can put the kernel of $P^{t}$ into the form,
\begin{multline}
    k_{P^{t}}(x,y)=|\varepsilon_{x}'|K_{P}(x,-\varepsilon_{x}(y))+|\varepsilon_{x}'|[(1-\chi)K](x,-\varepsilon_{x}(y))+R(y,x)\\ 
    = |\varepsilon_{x}'|K_{P}(x,-\varepsilon_{x}(y)) \quad \bmod C^{\infty}(U\times U).
\end{multline}
It then follows from Proposition~\ref{prop:PsiVDO.characterisation-kernel2} that $P^{t}$ is a \psivdo\ of order $m$. Moreover, working out the expression for 
$K_{\alpha\beta\gamma\delta}$ shows that the asymptotics expansion~(\ref{eq:Appendix.expansionKPt}) 
reduces to~(\ref{eq:PsiHDO.transpose-expansion-kernel}). The proof of Proposition~\ref{prop:PsiHDO.transpose-chart} 
is thus achieved. 
% \appendix 
% 
% \section*{Appendix: Condition $X(k)$ and hypoellipticity of the horizontal sublaplacian}
% 
% In this appendix we give a proof of Proposition~\ref{prop:Examples.horizontal-sublaplacian}. 

\end{document}